\documentclass[11pt,letterpaper]{amsart}
\usepackage[utf8]{inputenc}
\usepackage{amssymb,amsthm, graphics, comment,bm}
\usepackage{mathtools}
\usepackage{amsthm}
\usepackage{diagmac2}
\usepackage{textcomp}
\usepackage{color}
\usepackage[T1]{fontenc}
\usepackage[alphabetic]{amsrefs}
\usepackage[all,cmtip,color,matrix,arrow]{xy}
\usepackage{yfonts}
\usepackage{tikz}
\usepackage{tikz-cd}
\usepackage{stmaryrd}
\usepackage{mathrsfs}
\usepackage{enumerate}
\usepackage{hyperref}
\hypersetup{
	colorlinks=true,
	linkcolor=blue,
	filecolor=magenta,      
	urlcolor=cyan,
}
\usepackage[mathscr]{euscript}
\usepackage[all,cmtip]{xy}
\usepackage{cleveref}
\usepackage{fullpage}

\newcommand{\bA}{\mathbb{A}}

\newcommand{\bG}{\mathbb{G}}

\newcommand{\bN}{\mathbb{N}}
\newcommand{\bP}{\mathbb{P}}
\newcommand{\bQ}{\mathbb{Q}}

\newcommand{\bZ}{\mathbb{Z}}

\newcommand{\sC}{\mathcal{C}}
\newcommand{\sD}{\mathcal{D}}

\newcommand{\oH}{\operatorname{H}}

\newcommand{\sP}{\mathscr{P}}
\newcommand{\sR}{\mathscr{R}}
\newcommand{\sS}{\mathscr{S}}

\newcommand{\sX}{\mathcal{X}}

\newcommand{\cA}{\mathcal{A}}
\newcommand{\cB}{\mathcal{B}}
\newcommand{\cC}{\mathcal{C}}
\newcommand{\cD}{\mathcal{D}}
\newcommand{\cE}{\mathcal{E}}
\newcommand{\cF}{\mathcal{F}}
\newcommand{\cG}{\mathcal{G}}
\newcommand{\cH}{\mathcal{H}}
\newcommand{\cI}{\mathcal{I}}

\newcommand{\cL}{\mathcal{L}}
\newcommand{\cM}{\mathcal{M}}
\newcommand{\cN}{\mathcal{N}}
\newcommand{\cO}{\mathcal{O}}
\newcommand{\cP}{\mathcal{P}}
\newcommand{\cQ}{\mathcal{Q}}
\newcommand{\cR}{\mathcal{R}}
\newcommand{\cS}{\mathcal{S}}
\newcommand{\cT}{\mathcal{T}}
\newcommand{\cU}{\mathcal{U}}
\newcommand{\cV}{\mathcal{V}}

\newcommand{\cX}{\mathcal{X}}
\newcommand{\cY}{\mathcal{Y}}
\newcommand{\cZ}{\mathcal{Z}}

\newcommand{\spec}{\operatorname{Spec}}

\newcommand{\sY}{\mathscr{Y}}

\newcommand{\Sym}{\operatorname{Sym}}
\newcommand{\Hom}{\operatorname{Hom}}

\newcommand{\Proj}{\operatorname{Proj}}

\newcommand{\Pic}{\operatorname{Pic}}

\newcommand{\Gm}{\bG_{{\rm m}}}
\newcommand{\Gmdelta}{\bG_{{\rm m,\Delta}}}
\newcommand{\GL}{\operatorname{GL}}

\newcommand{\Aut}{\operatorname{Aut}}
\newcommand{\PGL}{\operatorname{PGL}}

\newcommand{\Id}{\operatorname{Id}}

\newcommand{\Bl}{\operatorname{Bl}}

\newcommand{\bmu}{\bm{\mu}}
\newcommand{\dm}{{\rm DM}}
\newcommand{\id}{\rm{id}}

\newcommand{\Mcal}{\mathcal{M}}
\newcommand{\Mbar}{\overline{\mathcal{M}}}
\newcommand{\Mtilde}{\widetilde{\mathcal{M}}}
\newcommand{\Ctilde}{\widetilde{\mathcal{C}}}
\newcommand{\Hass}{{\rm H}}
\newcommand{\git}{{\rm GIT}}
\DeclareMathOperator{\EB}{EB}
\newcommand{\CY}{{\rm CY}}

\newtheorem{theorem}{Theorem}[section]
\newtheorem{Teo}[theorem]{Theorem}
\newtheorem*{Teo*}{Theorem}
\newtheorem{Lemma}[theorem]{Lemma}

\newtheorem{Cor}[theorem]{Corollary}

\newtheorem{Prop}[theorem]{Proposition}
\newtheorem*{Ques*}{Question}

\theoremstyle{definition}

\newtheorem*{Oss'}{Remark}

\newtheorem{Def}[theorem]{Definition}
\newtheorem*{Def*}{Definition}
\newtheorem{Notation}[theorem]{Notation}
\newtheorem{Assumption}[theorem]{Assumption}

\newtheorem{Remark}[theorem]{Remark}
\newtheorem{Conjecture}[theorem]{Conjecture}

\theoremstyle{definition}
\newtheorem{example}[theorem]{Example}

\begin{document}
\title[Stable maps to quotient stacks with a properly stable point]{Stable maps to quotient stacks with a properly stable point}
	\author[A. Di Lorenzo]{Andrea Di Lorenzo}
	\address[A. Di Lorenzo]{University of Pisa, Italy}
	\email{andrea.dilorenzo@unipi.it}
	\author[G. Inchiostro]{Giovanni Inchiostro}
	\address[G. Inchiostro]{University of Washington, Seattle, Washington, USA}
	\email{ginchios@uw.edu}
	\maketitle
	\begin{abstract}
{We compactify moduli stacks of maps from curves to a large class of quotient stacks $[W/G]$ admitting projective good moduli spaces. Our construction extends quasimap-theoretic compactifications and relies on a new birational operation for algebraic stacks, which we call an \emph{extended weighted blow-up}. As applications, we construct compactifications of certain moduli spaces of fibrations in log Calabi--Yau pairs, of maps to stacks of the form $[W/\Gm^n]$, of maps to the GIT moduli stack of $2n$ points on $\bP^1$, and of maps to the GIT moduli stack of plane cubics. In the appendix, we use extended weighted blow-ups to give a modular proof of a conjecture of Hassett.}
	\end{abstract}
 \section{Introduction}
{Once a moduli space of curves has been constructed, a natural next question is to construct moduli spaces of curves equipped with additional geometric data, such as families of varieties over the curve, $G$-torsors, or related structures. When this additional data is parametrized by an algebraic stack $\cX$, the problem can be formulated as the construction of a moduli space of maps $(C,f\colon C\to \cX)$ from smooth curves to $\cX$. A central question in moduli theory is whether such moduli spaces admit compactifications whose boundary points retain a meaningful modular interpretation. We briefly recall several cases in which this problem has been addressed for maps from smooth curves to an algebraic stack $\cX$:}

 \begin{enumerate}
     \item for $\cX$ a projective variety by Kontsevich \cites{Kon, FP},
     \item for $\cX$ a Deligne--Mumford stack with a projective coarse moduli space by Abramovich and Vistoli \cite{AV_compactifying},
     \item for $\cX=\cB \Gm$, $\cX=[X/\Gm]$ with $X$ projective or $\cX=\cB \GL_n$ by Caporaso, Frenkel-Teleman-Tolland, Cooper and Pandharipande, and more recently Halpern-Leistner and Herrero \cites{Cap, FTT,cooper2022git, Pand, halpern2025quantum},
     \item for $\cX=[\spec(A)/G]$ and under the assumption that there is a character $\theta:G\to \Gm$ such that
$[\spec(A)(\Bbbk_\theta)^{ss}/G]$ is a Deligne--Mumford stack by Cheong, Ciocan-Fontanine, Kim and Maulik via quasimap theory \cites{ciocan2010moduli, ciocan2014stable, orbifold_qmap_theory},
\item Halpern--Leistner and Fernández Herrero treat maps from a fixed smooth curve to a global quotient in  \cite{HLH23}.
\end{enumerate}
In this paper we address this problem for maps from curves to a broad class of algebraic stacks, described below in increasing order of generality. Our strategy has two parts: we first construct compactifications when the target stack $\cX$ contains a proper Deligne–Mumford open substack $\cX_{\dm}\subseteq \cX$, with some additional mild assumptions; this is achieved in \Cref{teo_intro_cX_contains_open_proper_dm}, and in more detail in \Cref{thm:quasimaps are dm}. Then, we enlarge $\cX\hookrightarrow \widetilde{\cX}$ a general target $\cX$ with a properly stable point, so that to guarantee that such an enlargement contains a proper Deligne–Mumford open substack $\widetilde{\cX}_{\dm}\subseteq \widetilde{\cX}$ as above; this step leads to \Cref{thm_intro_simplified}. The boundary objects are \textit{stable quasimaps}, introduced in \Cref{def_intro_qmaps}: they are maps from nodal (twisted) curves, which generically map to the proper Deligne–Mumford locus, subject to a stability condition.

\medskip

Let $\cX=[W/G]$ be a quotient stack, where $G$ is reductive, admitting a projective good moduli space $\pi\colon \cX \to X$ with an ample line bundle $L_X$. Assume that there is an open substack $\cX_\dm\subseteq \cX$ that is proper and Deligne--Mumford. Combining the definition of quasimaps to affine GIT quotients \cite{ciocan2014stable} and \cite{dli1}*{Definition 2.2}, we make the following definition.

\begin{Def}\label{def_intro_qmaps}
    A stable $n$-marked \emph{quasimap} of genus $g$ is a pair \[(\phi:\sC\longrightarrow \cX, \Sigma_1,\ldots, \Sigma_n) \] such that
    \begin{enumerate}
        \item the pair $(\sC;\Sigma_1,\ldots,\Sigma_n)$ is a twisted, $n$-marked curve of genus $g$ as in \cite{AV_compactifying},
        \item the morphism $\phi:\sC\to \cX$ is representable and maps the generic points, the marked gerbes $\Sigma_i$ and the nodal locus of $\sC$ to $\cX_\dm \subseteq \cX$;
        \item the line bundle \[\omega_\cC(\Sigma_1+\ldots+\Sigma_n)\otimes \phi^*\pi^*L_X^{\otimes 3}\] is nef, and
        \item if $\sD\subseteq \sC$ is an irreducible component such that $\omega_\cC( \Sigma_1+\ldots + \Sigma_n)\otimes \phi^*\pi^*L_X^{\otimes 3}$ has degree zero on $\cD$, then its coarse moduli space $D$ is isomorphic to $\bP^1$ and there are two geometric points $p_1,p_2\in \sD$ such that $\phi(p_1)$ and $\phi(p_2)$ are not isomorphic points of $\cX$.
    \end{enumerate}
\end{Def}  
Recall that a twisted curve is a mild stack-theoretic generalization of a nodal curve. As a first approximation, one may think of a stable quasimap as a morphism 
$
\phi\colon\cC\to \cX$ from a twisted nodal curve $\cC$ such that $\phi^{-1}(\cX\smallsetminus \cX_{\dm})$ consists of finitely many smooth, unmarked points of $\cC$. These points are the base points of the quasimap $\phi$; conditions (3) and (4) instead are stability conditions.

\medskip

The following statement is a streamlined version of our main theorem, stated under the simplifying assumption that $\cX$ is smooth. 
\begin{Teo}\label{teo_intro_cX_contains_open_proper_dm}
Within the setup above, stable $n$-marked quasimaps of genus $g$ and class $\beta$ form a stack $\cQ_{g,n}(\cX,\cX_\dm,\beta)$ which is proper, Deligne--Mumford and carries a perfect obstruction theory.
\end{Teo}

Smoothness is imposed here only to simplify the statement: it suffices to assume that $\cX_{\dm}$ agrees with the relative semistable locus for a line bundle on $\cX$. 
Moreover, in the full version of our theorem, as can be found in \Cref{sec:epsilon qmaps}, we construct a proper Deligne-Mumford stack of $\epsilon$-stable quasimaps (\Cref{thm:epsilon quasimaps are dm}), of which the stack above is a particular case obtained by setting $0<\epsilon \ll 1$ (\Cref{thm:quasimaps are dm}). We expect this broader notion of $\epsilon$-stability to be useful in enumerative geometry and plan to pursue this direction in future work.

The theorem above has been used as a key ingredient in \cites{ISZ,IZ1,IZ2}, for producing compact moduli spaces of fibrations in log Calabi--Yau pairs; we will discuss further applications in \Cref{sec:first_applications}. 

\smallskip

Suppose now that $\cX$ does not contain a proper Deligne--Mumford open substack, but that there exists a dense open subset $U\subset X$ such that $\mathcal{U}=\cX\times_XU$ is Deligne--Mumford. When $\cX$ is irreducible, the existence of such an open subset $U$ is equivalent to the existence of a properly stable point; for a GIT quotient, this amounts to the stable locus being nonempty. In this situation, we construct in \Cref{thm_you_can_embed_a_Stack_in_one_with_open_dm_by_taking_deformations_to_nc} an algebraic stack $\widetilde{\mathcal{X}}$ with an open embedding $\mathcal{X}\subseteq \widetilde{\mathcal{X}}$ inducing an isomorphism on good moduli spaces, such that $\widetilde{\mathcal{X}}$ contains a proper open Deligne--Mumford stack $\widetilde{\mathcal{X}}_\dm$ compactifying $\mathcal{U}$. In particular, applying \Cref{teo_intro_cX_contains_open_proper_dm} to $\widetilde{\mathcal{X}}$ leads to one of our main theorems.
\begin{Teo}\label{thm_intro_simplified}
Let $\cX=[W/G]$ be a quotient stack with $G$ reductive, with a projective good moduli space $X$, and admitting a dense open subset $U\subseteq X$ such that $\mathcal{U}=\cX\times_XU$ is Deligne--Mumford. Then the proper Deligne--Mumford stack $\cQ_{g,n}(\widetilde{\cX},\widetilde{\cX}_\dm,\beta)$ compactifies the stack of maps $\phi\colon C\to \cX$ of class $\beta$ from a smooth, $n$-marked curve $C$ of genus $g$ whose generic point is mapped to $\mathcal{U}$.
The boundary of $\cQ_{g,n}(\widetilde{\cX},\widetilde{\cX}_\dm,\beta)$ consists of stable, $n$-marked quasimaps with target $\widetilde{\cX}$. 
\end{Teo}

The construction of the enlargement $\widetilde{\cX}$ of $\cX$ is \textit{explicit}, using \cite{ER21}: the proper Deligne--Mumford open substack $\widetilde{\cX}_{\dm}$ is obtained from the canonical reduction-of-stabilizers procedure of Edidin–Rydh \cite{ER21}, which generalizes Kirwan’s partial desingularization \cite{kirwan}. Specifically, we define \textit{extended weighted blow-ups}, which are a generalization of weighted blow-ups of \cite{QR}, and
the enlargement $\widetilde{\cX}$ of $\cX$ is obtained by performing a sequence of extended weighted blow-ups. The properties of the enlargement $\cX\subseteq \widetilde{\cX}$ are presented in \Cref{thm_you_can_embed_a_Stack_in_one_with_open_dm_by_taking_deformations_to_nc}.

\subsection{First applications} \label{sec:first_applications} 
We provide some applications of our main results. 
\begin{enumerate}
\item When $\cX$ is a certain moduli of boundary-polarized Calabi--Yau pairs. This application leads to a moduli stack parametrizing fibrations in log Calabi--Yau pairs.
    \item when $\cX$ admits a morphism $\cX\to \cB\Gm^n$ representable by Deligne--Mumford stacks. For example, all stacks of the form $[W/\Gm^n]$ admitting a projective good moduli space.
    \item when $\cX$ is the GIT compactification of the moduli stack of smooth plane cubics. In this example one can still give a modular interpretation of a morphism $\phi:C\to \widetilde{\cX}$, in terms of plane cubics, by composing $\phi$ with the projection $\widetilde{\cX}\to \cX$ of \Cref{thm_you_can_embed_a_Stack_in_one_with_open_dm_by_taking_deformations_to_nc}.
\end{enumerate}

\subsubsection{Applications to moduli of fibrations in log Calabi--Yau pairs}\label{subsectuon_k_dim_1} One of the main applications of \Cref{thm_intro_simplified} is for constructing certain moduli of varieties. 
 While substantial progress has been made in constructing moduli spaces for varieties of general type \cites{KSB,al1,kollarmodulibook}, and similarly for Fano varieties \cites{blum2019uniqueness,alper2020reductivity,blum2021properness,blum2022openness,xu2020uniqueness,MR4445441}, the situation for varieties of intermediate Kodaira dimension is less clear.
 
 The main application of \Cref{teo_intro_cX_contains_open_proper_dm} is \Cref{teo_intro_kodairadim1}: the
 construction of a compact moduli stack of certain pairs $(Y,cD)$ of Kodaira dimension one, which are fibrations $f\colon (Y,cD)\to C$ with fibers in a moduli space of log Calabi--Yau pairs. The quasimap condition and the stability condition of \Cref{def_stable_qmap} are translated in terms of the birational geometry of the corresponding fibered variety: labels (Sing), (Stab) and (Num) in the statement of \Cref{teo_intro_kodairadim1} stand for singularity condition, stability condition and numerical condition.
 
 \medskip
 
 We denote by $\cX:=\cD\cP^{\CY}_m$ the moduli space of \cites{blum2024good} parametrizing boundary polarized Calabi--Yau pairs $(S,cD)$ of dimension at most 2 admitting a $\bQ$-Gorenstein smoothing, and with $c<1$. 
\begin{Teo}\label{teo_intro_kodairadim1}
    The algebraic stack $\cX$ contains a proper Deligne--Mumford open substack $\cX_\dm$, denoted in \cite{blum2024good} by $\cD\cP^{\operatorname{KSBA}}_m$, such that, if $\cX\smallsetminus \cX_{\dm}$ has codimension at least two, the  assumptions of \Cref{teo_intro_cX_contains_open_proper_dm} apply for the inclusion $\cX_\dm\subseteq \cX$. In particular:
    \begin{itemize}
        \item the stack $\cQ_g(\cX,\cX_\dm,\beta)$ compactifies the moduli of fibrations $\pi\colon(Y,cD)\to C$ with fibers in $\cX$ such that:
        \begin{itemize}
            \item[(Sing)] the curve $C$ is smooth and the generic fiber of $\pi$ is has klt singularities,
            \item[(Stab)] either $\omega_C$ is ample, or not all the fibers of $\pi$ are $S$-equivalent,
            \item[(Num)] the family $\pi$ comes from a map $C\to \cX$ of class $\beta$ from a curve of genus $g$.
        \end{itemize}
        \item the boundary of $\cQ_g(\cX,\cX_\dm,\beta)$ parametrizes families $\pi:(\sY,c\sD)\to \sC$ of pairs in $\cX$ with fibers in $\cX$, fibered over a twisted curve $\sC$, such that:
        \begin{itemize}
            \item[(Sing)] the set $\Delta:=\{p\in \sC:(\sY_p,(c+\epsilon)\sD_p)$ does \underline{not} have semi-log-canonical singularities for any $0<\epsilon \ll 1\}$ is a finite union of smooth points of $\sC$,
            \item[(Stab)] if $\sR\subseteq \sC$ is an irreducible component such that $\deg(\omega_\sC |_\sR)< 0$, then not all the fibers of $\pi|_\sR\colon(\sY|_\sR,c\sD|_\sR)\to \sR$ are $S$-equivalent; whereas if $\deg(\omega_\sC |_\sR)=0$, then not all fibers of $\pi|_\sR$ are isomorphic, and
            \item[(Num)] the family $\pi$ comes from a map $\sC\to \cX$ of class $\beta$ from a twisted curve of genus $g$.
        \end{itemize}
    \end{itemize}
\end{Teo}
\begin{Remark}Being $S$-equivalent is a weaker condition than being isomorphic; one can consult \cite{ascher2023moduli}*{Definition 6.7} for the definition.
    Moreover, from \cite{ascher2023moduli}*{Proposition 14.17}, one can rephrase condition (S) in point (1) of \Cref{teo_intro_kodairadim1} by asking that either $\omega_C$ is ample, or the moduli part in the canonical bundle formula for $\pi$ has positive degree. 
\end{Remark}
\begin{Remark}
    As stated, the fibrations $\pi\colon (Y,cD)\to C$ in part (1) of \Cref{teo_intro_kodairadim1}
    must have \textit{all} fibers that are in $\cD\cP^{\CY}_m$ and the generic fiber in $\cD\cP^{\operatorname{KSBA}}_m$. In \Cref{remark_stefano_some_fibers_are_not_slc} we explain how to use the pointed version of \Cref{teo_intro_cX_contains_open_proper_dm}
     to treat the case when one has a morphism $\pi'\colon (Y,D)\to C$ over a smooth curve $C$ such that $D$ does not have any vertical component and only the generic fiber of $\pi$ is in $\cD\cP^{\operatorname{KSBA}}_m$.
\end{Remark}
\subsubsection{Applications to target stacks of the form $[W/\Gm^n]$}
Our second main application is for target stacks of the form $[W/\Gm^n]$. 
In \Cref{teo_toric_stack_has_open_dense_dm} we show that no enlargement is needed: one may take 
one can take $\cX=\widetilde{\cX}$, as we prove that any algebraic stack of the form $[W/\Gm^n]$ which satisfies the assumptions of \Cref{thm_intro_simplified} already contains a proper Deligne--Mumford stack. Combined with \Cref{teo_intro_cX_contains_open_proper_dm}, this gives a proper moduli stack of maps to algebraic stacks of the form $[W/\Gm^n]$, admitting a projective good moduli space, and a properly stable point.

\subsubsection{Applications to fibrations in plane cubics} In this section we give an example of how to use the \Cref{thm_you_can_embed_a_Stack_in_one_with_open_dm_by_taking_deformations_to_nc} and \Cref{teo_intro_cX_contains_open_proper_dm} to compactify the space of maps to the GIT moduli space of plane cubics. Recall that this GIT moduli space is defined as $ [\bP(\oH^0(\cO_{\bP^2}(3)))^{ss}/\PGL_3]$, where $\bP(\oH^0(\cO_{\bP^2}(3)))^{ss}$ is the open in $\bP(\oH^0(\cO_{\bP^2}(3)))$ parametrizing plane cubics with at most nodal singularities. This stack has neither a proper Deligne–Mumford open substack nor, as far as we know, a modular enlargement containing one.

Nevertheless, the properties of extended weighted blow-ups allow us to give a detailed description of the objects on the boundary of the space of quasimaps to an appropriate extended blow-up of $[\bP(\oH^0(\cO_{\bP^2}(3)))^{ss}/\PGL_3]$.

\subsection{Extended weighted blow-ups}

To prove \Cref{thm_intro_simplified}, the key step is to construct an enlargement $\cX\subseteq \widetilde{\cX}$ such that the assumptions of \Cref{teo_intro_cX_contains_open_proper_dm} hold. For doing so, we introduce a new birational transformation for algebraic stacks, which we call \emph{extended blow-up}. On a first approximation, an extended weighted blow-up is an enlargement of an algebraic stack which simultaneously contains the original stack and its weighted blow-up as open substacks, while preserving the good moduli space. Applying this construction multiple times, we prove the following result.
\begin{Teo}\label{thm_you_can_embed_a_Stack_in_one_with_open_dm_by_taking_deformations_to_nc}
    Let $\cX$ be an algebraic stack with a good moduli space $p:\cX\to X$, and with an open dense $U\subseteq X$ such that $\cX\times_XU$ is Deligne--Mumford. Then there is an algebraic stack $\widetilde{\cX}$ with an open embedding $i:\cX\hookrightarrow \widetilde{\cX}$ such that:
    \begin{enumerate}
        \item $\widetilde{\cX}$ has a good moduli space which is isomorphic to $X$ via the inclusion $i$,
        \item there is a morphism $\pi:\widetilde{\cX}\to \cX$ which is an isomorphism over $(p\circ \pi)^{-1}(U)$,
        \item the morphism $\pi\circ i$ is isomorphic to the identity,
        \item there is a line bundle $\cL_{DM}$ on $\widetilde{\cX}$ such that the open substack $\widetilde{\cX}(\cL_{DM})^{ss}_X \subseteq \widetilde{\cX}$ given by the $\cL_{DM}$-semistable locus over $X$ is Deligne--Mumford and proper over $X$,
        \item if $\cX$ is a global quotient by a reductive group, then $\widetilde{\cX}$ can be chosen to be a global quotient by a reductive group.
    \end{enumerate}
\end{Teo}
In particular, from \Cref{teo_intro_cX_contains_open_proper_dm}, if $\cX$ is a global quotient and one fixes integers $g$, $n$ and a class $\beta$, there is a moduli space $\cQ_{g,n}(\widetilde{\cX},\widetilde{\cX}(\cL_{DM})^{ss}_X,\beta)$ of stable quasimaps to $\widetilde{\cX}$ which is a proper Deligne--Mumford stack, and compactifies the smooth quasimaps to $\mathcal{\cX}$, thus yielding \Cref{thm_intro_simplified}. 


\medskip

Finally, in the appendix we establish a criterion (see \Cref{prop:criterion}) for telling when a morphism of algebraic stacks is an extended weighted blow-up. For $\cX$ the GIT compactification of the moduli stack of $2n$ distinct points on $\bP^1$, we use this criterion to show that a \textit{modular} enlargement $\cX\subseteq \widetilde{\cX}$ is an extended weighted blow-up,
and we use it to give a modular proof of a conjecture of Hassett:
\begin{Conjecture}\label{conjecture_hassett}
    Let $\overline{\cM}_{0,(\frac{1}{n}  + \epsilon,\ldots,\frac{1}{n}  + \epsilon)}$ be the moduli space of weighted $2n$-pointed stable curves of genus 0 and with weights $\frac{1}{n}+\epsilon$. Let $\widetilde{M}_{0,2n}$ the coarse moduli space of $[\overline{\cM}_{0,(\frac{1}{n}  + \epsilon,\ldots,\frac{1}{n}  + \epsilon)}/S_{2n}]$ where the action permutes the $2n$ points. Let $\bP(\oH^0(\bP^1,\cO_{\bP^1}(2n)))/\!\!/\PGL_2$ the GIT moduli space of $2n$ unordered points on $\bP^1$, linearized via $\cO(2)$. Then there is a map $\widetilde{M}_{0,2n}\to \bP(\oH^0(\bP^1,\cO_{\bP^1}(2n)))/\!\!/\PGL_2$ which is a weighted blow-up.
\end{Conjecture}
\Cref{conjecture_hassett} was proven in \cite{moduli_pointed_curves_via_git} by showing that the corresponding statement for ordered points is true: there is a map $\overline{\cM}_{0,(\frac{1}{n}  + \epsilon,\ldots,\frac{1}{n}  + \epsilon)}\to (\bP^1)^n/\!\!/\PGL_2$ that is a weighted blow-up (in this case, it is the Kirwan desingularization). The conjecture was then proved by observing that the previous map is $S_{2n}$-equivariant, and the $S_{2n}$-quotient induces $\widetilde{M}_{0,2n}\to \bP(\oH^0(\bP^1,\cO_{\bP^1}(2n)))/\!\!/\PGL_2$ which is a weighted blow-up.
It is natural to wonder if one can construct the map above from the modular properties of the quotient stack $[\bP(\oH^0(\bP^1,\cO_{\bP^1}(2n)))(\cO(2))^{ss}/\PGL_2]$, which we will denote by $\cC^{\git}_{2n}$, directly. Specifically, the question is whether one can construct a moduli stack of rational $2n$-marked curves $\Mcal$, having $\widetilde{M}_{0,2n}$ as coarse moduli space, and with a map $\cM\to \cC^{\git}_{2n}$ which induces the desired weighted blow-up on good moduli spaces. This is the content of the next theorem, proved in \Cref{appendix}.

\begin{Teo}\label{thm_intro_hassett}
Consider pairs $(P,D)$ where:
        \begin{itemize}
            \item either $P\cong \bP^1$ or $P$ is a twisted conic, namely a twisted curve that is the nodal union of two root stacks of $\bP^1$ at $0$, glued along the $\cB\bmu_2$ gerbe, 
            \item $D$ is supported on the smooth locus of $P$, has degree $2n$ and each point of $P$ has multiplicity at most $n$ in $D$ (in other terms, $(P,\frac{1}{n}D)$ is semi log-canonical), and
            \item the divisor $\omega_P(\frac{1}{n}D)$ is nef.
        \end{itemize}
        Then:
        \begin{enumerate}
            \item there is an algebraic stack $\Ctilde^{\CY}_{2n}$ which is a moduli stack for pairs as above,
            \item there is a morphism $\Phi:\Ctilde^{\CY}_{2n}\to \cC^{\git}_{2n}$ which is an extended weighted blow-up with center the polystable but not stable point of $\cC^{\git}_{2n}$,
            \item if we denote by $X$ the good moduli space of $\cC^{\git}_{2n}$, there is a line bundle $\cL$ on $\Ctilde^{\CY}_{2n}$ such that the open substack $\Ctilde^{\CY}_{2n}(\cL)^{ss}_X\subseteq\widetilde{\cC}^\CY_{2n}$ is a proper Deligne--Mumford stack and its coarse moduli space is a weighted blow-up of $C_{2n}^{\git}$ at the unique polystable but not stable point, and
            \item the pairs in $\Ctilde^{\CY}_{2n}(\cL)^{ss}_{X}$ are pairs $(P,D)$ as above such that each point of $P$ has multiplicity at most $n-1$ in $D$. 
        \end{enumerate}
        In particular, the weighted blow-up constructed in \cite{moduli_pointed_curves_via_git} is induced on good moduli spaces by the morphism $\Ctilde_{2n}^{\CY}(\cL)^{ss}_{X}\hookrightarrow\Ctilde_{2n}^{\CY}\xrightarrow{\Phi}\cC^{\git}_{2n}$.
        \end{Teo}
        \Cref{remark_why_we_need_twisted_nodes} explains why one needs to consider twisted curves rather than schematic curves.

 \subsection{Organization of the paper}

 In  \Cref{section_background}, we collect preliminary results on good moduli spaces and relative GIT. In \Cref{sec:epsilon qmaps}, we define stable and $\epsilon$-stable quasimaps and prove \Cref{teo_intro_cX_contains_open_proper_dm}. \Cref{section_examples_qmaps} contains applications to fibrations in log Calabi–-Yau pairs and to quotient stacks by tori. In \Cref{section_extended_blowups} we introduce extended weighted blow-ups and prove \Cref{thm_you_can_embed_a_Stack_in_one_with_open_dm_by_taking_deformations_to_nc}.  \Cref{section_examples_extended} applies the construction to the GIT quotient stack of plane cubics. 
 
 Finally, \Cref{appendix} gives a criterion for recognizing extended weighted blow-ups and applies it to prove \Cref{thm_intro_hassett}.

 \subsection{Acknowledgements} We are thankful to Luca Battistella, Dori Bejleri, Harold Blum, Daniel Halpern-Leistner, Felix Janda, Nick Kuhn, Denis Nesterov and Dhruv Ranganathan for helpful conversations. We are thankful to Stefano Filipazzi and Andres Fernandez Herrero for their detailed feedback on an earlier version of this manuscript. We are particularly thankful to Andres Fernandez Herrero for suggesting an improvement on \Cref{lemma_one_can_check_rep_at_poli_pts}.
 
 \subsection{Conventions.} We work over a field, denoted by $k$, which is of characteristic 0 and algebraically closed. Unless otherwise stated, all the stacks have affine diagonal, are of finite type over $k$, and the stabilizers of the points are reductive (but a priori not connected).
 We denote the set of characters (resp. cocharacters) of a group $G$ by $\mathbf{X}(G)$ (resp. $\mathbf{X}(G)^*)$, and the natural pairing between them by \[<\cdot, \cdot >:\text{ }\mathbf{X}(G)^*\times \mathbf{X}(G)\to \bZ. \] 
 Given a character $\chi\in \mathbf{X}(G)$, we denote by $\Bbbk_\chi$ the representation of $G$ with character $\chi$.
 We will denote by $\bmu_n$ the group of $n$-th roots of unity of $k$ and by $\Theta:=[\bA^1/\Gm]$. Unless otherwise stated, all curves are assumed to be connected and proper.
 \section{Good moduli spaces and relative GIT}\label{section_background}
In this section we report a few results on algebraic stacks, good moduli spaces and geometric invariant theory which will be useful for the rest of the manuscript. The first new tool we produce is \Cref{prop_stable_locus_of_pi*L^n_otimes_G}, which roughly speaking allows us to express a sequence of relative semistable loci over a base $X$ as one single relative semistable locus directly over $X$. Although technical, this proposition turns out to be useful in the sequel. 

The next tools, namely \Cref{lemma_git_for_ample_linearization_is_contained_in_the_affine_git_for_the_cone} and \Cref{lemma_cone_construction_over_X}, allow us to perform some standard reductions when dealing with quotient stacks and relative semistable loci. In \Cref{sec:relative semistable loci} we show why the hypothesis that an open substack is a relative semistable locus is not restrictive as it may seem.
 \subsection{Relative GIT}\label{subsection_relative_GIT}
 The goal of this subsection is to prove \Cref{prop_stable_locus_of_pi*L^n_otimes_G}, which is the (cohomologically) affine version of \cite{ER21}*{Proposition 3.18}, and is used a few times throughout the paper. For doing so, we begin by recalling the following definition.
 
 \begin{Def}[\cite{Alp}*{\S 11}]
     Let $p:\cX\to S$ be a morphism from an algebraic stack $\cX$ to a scheme $S$, and let $\cL$ be a line bundle on $\cX$. We define the $\cL$-\textit{semistable locus} of $\cX$ over $S$, denoted by $\cX(\cL)^{ss}_S$, as the locus of points $x\in \cX$ such that there exist $U\subset S$ an open neighbourhood of $p(x)$ and a section $s\in \oH^0(U, p_*(\cL^{\otimes n}))=H^0(p^{-1}U,\cL^{\otimes n})$ for some $n>0$ such that
     \begin{enumerate}
         \item the section $s$ does not vanish at $x$, and
         \item the induced morphism $p^{-1}(U)_s \to U$ is cohomologically affine.
     \end{enumerate} 
 \end{Def}
\begin{Def}\label[Def]{def_unstable_locus}
    Given a good moduli space $p\colon \cX\to X$ and a line bundle $\cL$ on $\cX$, the \textit{$\cL$-unstable locus} over $X$ is defined over an affine open subset $U\subseteq X$ as the vanishing locus of the ideal generated by the sections of $\cL(p^{-1}U)$.
\end{Def}
\begin{Def}
    Assume that a group $G$ acts on an affine scheme $\spec(R)$, let $\chi\colon G \to \Gm$ be a character of $G$ and let $f\in R$. We say that $f$ is semi-invariant for the character $\chi$ if for every $g\in G$ we have $g*f = \chi(g)f$. 
\end{Def}
 \begin{Remark}
     As observed in \cite{ER21}*{Remark 3.16} and the paragraph that follows, the previous definition can be extended to the case in which $S$ is an algebraic stack, as being cohomologically affine commutes with flat descent for morphisms of stacks with quasi-affine diagonal from \cite{Alp}*{Proposition 3.9.(vii)}. More explicitly, for every flat morphism $S'\to S$, inducing $g:\cX':=\cX\times_SS'\to \cX$ we have that $g^{-1}(\cX(\cL)^{ss}_S) = \cX'(g^*\cL)^{ss}_{S'}$.
 \end{Remark}
 \begin{Remark}
     Let $p:\cX\to S$ be a morphism from an algebraic stack $\cX$ to a scheme $S$, and let $\cL$ be a line bundle on $\cX$. Then from \cite{Alp}*{Theorem 11.5} the stack $\cX(\cL)^{ss}_S$ admits a good moduli space over $S$.
 \end{Remark}

\begin{Lemma}\label{lm:character bundle}
    Let $\cX=[\spec(A)/G]$ be an algebraic stack with good moduli space $X=\spec(A^G)$ and $\cL$ a line bundle on $\cX$. Then there exists an affine $G\times \Gm$-scheme $\spec(A')$ over $\spec(A)$ such that
    \begin{enumerate}
        \item there is an isomorphism $\cX\simeq[\spec(A')/G\times \Gm]$, and
        \item the line bundle $\cL$ is isomorphic to $\Bbbk_\chi$ for some character $\chi\colon G\times \Gm\to\Gm$.
    \end{enumerate}
\end{Lemma}
\begin{proof}
    Let $\cY\to \cX$ be the $\Gm$-torsor determined by $\cL$. Then the cartesian product $\cY \times_{\cX} \spec(A)$ is (1) a $\Gm$-torsor over $\spec(A)$ and (2) a $G$-torsor over $\cY$. From (1) we deduce that $\cY \times_{\cX} \spec(A)$ is affine over $\spec(A)$, hence isomorphic to $\spec(A')$ for some $A$-algebra $A'$; from (2) we deduce that $\cY=[\spec(A')/G]$, hence $\cX=[\spec(A')/G\times \Gm]$.
    By construction, given the character $\chi\colon G\times \Gm \to \Gm$ coming from the second projection, we have that $\cL$ is the pull-back of $\Bbbk_{\Id}$.
\end{proof}
Let $\Theta:=[\bA^1/\Gm]$, with geometric points $1$ and $0$.
Recall that given a stack $\cX$ and a line bundle $\cL$ on $\cX$, for any geometric point $p$ in $\cX$ and any morphism $\lambda\colon\Theta\to\cX$ which maps $1\to p$, we can define
\[\mu^{\cL}(p,\lambda)=-\text{weight of }(\lambda^*\cL)|_{[0/\Gm]}.\]
For $\cX=[\spec(A)/G]$ with $G$ reductive, we have \cite{Alp21}*{Proposition 6.9.1}
\[ \Hom(\Theta,\cX) = \left\{ x \in \spec(A)\text{, }\lambda\colon\Gm\to G \text{ such that there exists }\lim_{t\to 0}\lambda(t)\cdot x\right\}.\]
Moreover, for a point $x\in \spec(A)$ which is mapped to $p$, a character $\chi:G\to\Gm$ and $\lambda:\Gm\to G$ a cocharacter such that $\lim_{t\to 0}\lambda(t)\cdot p = y \in \spec(A)$, we have
\[ \mu^\chi(x,\lambda) = \mu^{\Bbbk_\chi}(p,f)\]
where $f(1)=p$ and $f(0)=q$ is the image of $y$ via $\spec(A)\to\cX$. As it is common, we will refer to both these functions as Hilbert-Mumford functions.
\begin{Lemma}\label{lemma_relative_GIT_for_gms_can_be_checked_using_HM}
     Let $\cX$ be an algebraic stack with a good moduli space $X$, a line bundle $\cL$ and $x\in \cX$ a point. Then the following conditions are equivalent:
     \begin{enumerate}
         \item $x\in \cX(\cL)^{ss}_X$, and
         \item for every morphism $\lambda:\Theta\to \cX$ sending $1\mapsto p$, the weight of $\lambda^*\cL$ on $\cB\Gm$ is non-positive.
     \end{enumerate}
     
     Moreover, if $\cX=[\spec(A)/G]$ and $\cL=\Bbbk_\chi$ for a $G$-character $\chi$, the conditions above are equivalent to the existence of $f\in A$ which is $\chi^d$-semi-invariant for some $d>0$, and which does not vanish on any point $q\to \spec(A)\to [\spec(A)/G]$ isomorphic to $p$.
 \end{Lemma}
 \begin{proof}
     Since conditions (1) and (2) are \'etale local on $X$, from \cite{luna}*{Theorem 4.12} we can assume that $\cX=[\spec(B)/\Gamma]$ and $X=\spec(B^\Gamma)$. 
     By \Cref{lm:character bundle}, we can assume that $\cL\simeq\Bbbk_\chi$ comes for a $\Gamma$-character $\chi$.
Now, a point $x$ is $\Bbbk_\chi$-semistable in $[\spec(B)/\Gamma]$ if and only if the corresponding point $p\in\spec(B)$ is $\chi$-semistable. The desired statement follows from the Hilbert-Mumford criterion for line bundles linearized by characters on affine rings \cite{Hos}*{Proposition 2.5}. 
 \end{proof}
Recall that, given a reductive group $G$, we can always define a norm $|\cdot |$ on $\Hom(\Gm,G)$ which is invariant under conjugation \cite{Hos}*{\S 2.2}.
 In what follows, when we talk of a \emph{normalized} Hilbert-Mumford function on a quotient stack $[W/G]$ with $G$ reductive, we are implicitly fixing a norm $|\cdot |$ on the lattice of cocharacters of $G$. 
This defines an element $\oH^4(\cB G,\bQ)$, hence an element $b\in \oH^4([W/G],\bQ)$ via pullback. Given $f\colon \Theta \to \cX$, we set $|f|$ to be $f^*b\in \oH^4(\Theta,\bQ)\simeq\bQ$.
 
\begin{Remark}\label{remark_norm_on_maps_from_theta_has_same_info_as_norm_from_1ps}
    If $\lambda\colon \Gm \to G$ such that $\lim_{t\to 0}\lambda(t)\cdot p=q$, the norm of the induced $f\colon\Theta\to [W/G]$ which maps $1\mapsto p$ and $0\mapsto q$ coincides with the previously defined norm of $\lambda$ \cite{HalInst}*{Example 4.1.17}. Moreover, from \cite{HalInst}*{Example 4.1.17}, for each $f\colon\Theta\to [W/G]$ there is a one parameter subgroup $\lambda\colon\Gm\to G$ such that $\lim_{t\to 0}\lambda(t)p = q$ and $\frac{\mu^{\Bbbk_\chi}(p,f)}{|f|} = \frac{\mu^{\chi}(p,\lambda)}{|\lambda|}$.
\end{Remark} 
 \begin{Lemma}\label{lm:bounded}
     Let $G$ be a reductive group and let $\cX=[\spec(A)/G]$ be a quotient stack. Let $\chi\colon G\to \Gm$ be a character.
     Then:
     \begin{enumerate}
         \item given an unstable point $p$, there exists $f_{\min}\colon\Theta\to\cX$ such that
         \[\frac{\mu^{\Bbbk_\chi}(p,f_{\min})}{|f_{\min}|} = \inf_{f(1)\simeq p} \left(\frac{\mu^{\Bbbk_\chi}(p,f)}{|f|}\right);\]
         \item  if we denote by $\cX^c$ the points $p\in \cX$ such that there is $f:\Theta\to \cX$ such that $f(1)\simeq p$ and $f(0)\not \simeq p$, the function \[\cX^c\to \mathbb{R},\text{ } p\mapsto \inf\left(\frac{\mu^{\Bbbk_\chi}(p,f)}{|f|}\text{ such that } f:\Theta\to \cX, f(1)\simeq p\text{ and }f(1)\not \simeq f(0)\right)\] takes finitely many values
     \end{enumerate}
 \end{Lemma}
 \begin{proof}
     If $G$ is a reductive group acting on an affine scheme $\spec(A)$ and $\chi\colon G\to\Gm$ a character, for any unstable point $x$ there exists a cocharacter $\lambda_{min}:\Gm\to G$ such that the normalized Hilbert-Mumford function
    \[ (x,\lambda) \longmapsto \frac{\mu^{\chi}(x,\lambda)}{|\lambda|} \]
    reaches its minimum value.
    Moreover, the set of minimal values that the normalized Hilbert-Mumford function reaches over the set of points $p\in \spec(A)$ whose $G$-orbit is not closed is finite. These statements follow from \cite{Hos}*{\S 2.2}, especially the discussion after \cite{Hos}*{Lemma 2.13}; we sketch the salient steps. 
    
    First observe that from \Cref{remark_norm_on_maps_from_theta_has_same_info_as_norm_from_1ps}
    it suffices to consider morphisms $\Theta\to \cX$ which come from one parameter subgroups. If $\lambda:\Gm\to G$ is a one parameter subgroup and $g\in G$, then $$\frac{\mu^{\chi}(x,\lambda)}{|\lambda|} = \frac{\mu^{\chi}(gx,g\cdot\lambda\cdot g^{-1})}{|g\cdot\lambda\cdot g^{-1}|},$$ where we denoted by $g\cdot\lambda\cdot g^{-1}$ the one parameter subgroup $\Gm\to G$ that sends $t\mapsto g\lambda(t)g^{-1}$.
        All one parameter subgroups are contained in a maximal torus of $G$ and all maximal tori are conjugated, so it suffices to check the desired statement for one parameter subgroups contained in a fixed maximal torus of $T$ of $G$. We can $T$-equivariantly embed $\spec(A)$ in $\bA^n$, so it suffices to check the desired statement for a $T$-action on $\bA^n$, which we can further assume to be diagonal as each $T$-representation of $T$ is diagonalizable. In this case, one can check the desired statement directly.

    To conclude, observe that for any $\Bbbk_\chi$-unstable point $p$ in $\cX$ and any point $x$ of $\spec(A)$ mapping to $p$, and any cocharacter $\lambda$ inducing $f\colon\Theta\to \cX$, we have $\mu^{\chi}(x,\lambda)=\mu^{\Bbbk_\chi}(p,f)$ and $|\lambda|=|f|$.
 \end{proof}

\begin{Prop}\label{prop_stable_locus_of_pi*L^n_otimes_G}
    For $i=1,2$, let $\cX_i$ be an algebraic stack with the same good moduli space $X$. Suppose that there is a morphism $\pi:\cX_2\to \cX_1$ inducing an isomorphism on good moduli spaces and such that the natural map $\cO_{\cX_1}\to \pi_*\cO_{\cX_2}$ is an isomorphism. Let $\cL$ be a line bundle on $\cX_1$ and let $\cM$ be a line bundle on $\cX_2$. Set $\cU_1:=\cX_1(\cL)^{ss}_{X}$ and denote $\cU_1\to U_1$ its good moduli space. Let $\cU_2:=(\pi^{-1}(\cU_1))(\cM)^{ss}_{U_1}$.
    
    Then, for $m\gg 0$, there is an equality $\cU_2 = \cX_2(\pi^*\cL^{\otimes m}\otimes \cM)^{ss}_{X}$. In other terms, the locus in $\cX_2$ which maps to $\cU_1$ and is $\cM$-semistable over $U_1$, agrees with the $\pi^*\cL^{\otimes m}\otimes \cM$-semistable locus of $\cX_2$ over $X$. 
\end{Prop}
\Cref{prop_stable_locus_of_pi*L^n_otimes_G} is the affine version of \cite{ER21}*{Proposition 3.18}, with a similar proof. The only exception is Step 2, which adapts some of the arguments in \cite{reichstein1989stability}*{Proof of Theorem 2.1}.
\begin{proof}
We will use several times the fact that the map $\pi$ is cohomologically affine, since it induces an isomorphism on good moduli spaces \cite{Alp}*{Proposition 3.13}.

\textbf{Step 1.} We show that the desired statement follows from the case $\cX_2=\cX_1$.

First observe that $\xi:\pi^{-1}(\cX_1(\cL)^{ss}_X)\to \cX_1(\cL)^{ss}_X\to U_1$ is still a good moduli space as it is cohomologically affine since it is a composition of cohomologically affine morphisms, and $\cO_{U_1}\to \xi_*\cO_{\pi^{-1}
(\cX_1(\cL)^{ss}_X)}$ is an isomorphism.  It suffices to check that \begin{equation}\label{eq_2}\cX_2(\pi^*\cL)^{ss}_X = \pi^{-1}(\cX_1(\cL)^{ss}_X).\end{equation} Indeed, if this is the case, and if we know that the desired statement holds in the case $\cX_2=\cX_1$, then \[
\cU_2 = (\pi^{-1}(\cX_1(\cL)^{ss}_X))(\cM)^{ss}_{U_1} = (\cX_2(\pi^*\cL)^{ss}_X)(\cM)^{ss}_{U_1} =\cX_2((\pi^*\cL)^{\otimes m}\otimes\cM)^{ss}_X
\]
where the first equality is the definition of $\cU_2$, the second one is \Cref{eq_2} together with the fact that $\xi$ is still a good moduli space, the third one is the desired statement in the case $\cX_1=\cX_2$.

It suffices that we show the two inclusions in \Cref{eq_2} set theoretically, since both stacks are open substacks of $\cX_2$. If $x_1\in \cX_1(\cL)^{ss}_X$ then, if $p:\cX_1\to X$ is the good moduli space, there is $W\subseteq X$ an open subscheme and $s\in p_*(\cL^{\otimes n})(W)$ such that $s$ does not vanish at $x_1$ and the locus $p^{-1}(W)_s\subseteq \cX_1$ where $s$ does not vanish is cohomologically affine over $W$. Since $\cX_2\to\cX_1$ is cohomologically affine, also \[(p\circ\pi)^{-1}(W)_{\pi^*s} = \pi^{-1}(p^{-1}(W)_s)\] is cohomologically affine over $W$, and $\pi^*s$ is also a section of $\pi^*\cL$. Since $\pi^*s$ does not vanish along $\pi^{-1}(x_1)$ we have that $\pi^{-1}(x_1)\subseteq \cX_2(\pi^*\cL)^{ss}_X$.

For the other inclusion, let $x_2\in \cX_2(\pi^*\cL)^{ss}_X$. Then there is $W$ an open neighbourhood of $p(\pi(x_2))$ and a section $\widetilde{s}\in p_*\pi_*\pi^*\cL(W)$ such that $(p\circ\pi)^{-1}(W)_{\widetilde{s}}\to W$ is cohomologically affine, and $\widetilde{s}$ does not vanish at $x_2$. Observe that, since $\pi_*\cO_{\cX_2}=\cO_{\cX_1}$, we have that
$\widetilde{s}=\pi^*s$ for a section $s\in p_*\cL(W)$, and \[(p\circ\pi)^{-1}(W)_{\widetilde{s}}=(p\circ\pi)^{-1}(W)_{\pi^*s} = \pi^{-1}(p^{-1}(W)_s)\]
where the first equality follows since $\widetilde{s}=\pi^*s$ and the second equality follows by definition of non-vanishing locus.
We need to show that $p^{-1}(W)_s\to W$ is cohomologically affine. So let $\cF\to \cG$ be a surjective morphism of quasi-coherent sheaves on $p^{-1}(W)_s$, we need to show that $p_*\cF\to p_*\cG$ remains surjective. For doing so, it suffices to check that the natural map \[\cH\to\pi_*\pi^*\cH\] is an isomorphism for any quasi-coherent sheaf $\cH$. Indeed, if this is the case then $\pi^*\cF\to \pi^*\cG$ will be surjective, which implies that $p_*\pi_*\pi^*\cF\to p_*\pi_*\pi^*\cG$ will be surjective since $p\circ\pi:(p\circ\pi^{-1})(W_s)\to W$ is cohomologically affine, and it will agree with $p_*\cF\to p_*\cG$. 

Consider a smooth atlas $Z\to p^{-1}(W)_s$ from an affine scheme $Z$. If we denote by \[\cY:=\pi^{-1}(p^{-1}(W)_s)\times_{p^{-1}(W)_s}Z,\] then the second projection $\mu:\cY\to Z$ is a good moduli space. Indeed, it is cohomologically affine since $\cX_2\to \cX_1$ is cohomologically affine and being cohomologically affine commutes with base change for a quasi-affine morphism \cite{Alp}*{Proposition 3.10.(vi)}. The morphism $Z\to p^{-1}(W)_s$ is quasi-affine (in fact, affine) as $Z$ is affine and $p^{-1}(W)_s$ has affine diagonal. Moreover $\mu_*\cO_\cY=\cO_{Z}$ since $\pi_*\cO_{\cX_2}=\cO_{\cX_1}$. Then \cite{Alp}*{Proposition 4.5} applies. This finishes Step 1. 
\color{black}

Therefore we can assume that $\cX_1=\cX_2$ and we will adopt the following notation:
\begin{itemize}
    \item we will denote $\cX_1=\cX_2$ by $\cX$, with the good moduli space morphism $p:\cX\to X$,
    \item we will not write $\pi$ any more as $\pi=\Id$ when $\cX_1=\cX_2$.
\end{itemize}

Let $\cV_m:=\cX(\cL^{\otimes m}\otimes \cM)^{ss}_{X}$.  The proof is in three steps, but first we perform some reductions. To start, observe that our claim is \'etale local on the good moduli space $X$, hence:
\begin{enumerate}
    \item we can assume that $\cX=[\spec(A)/G]$ from \cite{AHREtale};
    \item we can assume that the line bundles $\cL$ (respectively $\cM$) come from a character $\theta$ of $G$ (a character $\chi$ of $G$) by applying \Cref{lm:character bundle} twice;
    \item semistability with respect to $\cL$ and $\cM$ can be checked using the Hilbert-Mumford criterion by \Cref{lemma_relative_GIT_for_gms_can_be_checked_using_HM}; the same conclusion holds for $\cG_m := \cL^{\otimes m}\otimes \cM$;
    \item by \Cref{lm:bounded}, for $\cN=\cL, \cM$, we can assume that for every $\cN$-unstable point $q$ there exists $\lambda_{\min}$ which minimizes $\mu^{\cN}(q,\lambda)/|\lambda|$; furthermore, the set of such minimal values is finite.
\end{enumerate}
\textbf{Step 2.} We prove that there is $m_0$ such that, for every $m\ge m_0$, if $q\in \cX\smallsetminus \cU_1$, then $q\notin \cV_m$. 
For this, it is enough to prove that there exists an $m_0$ such that, for every $m\geq m_0$ and $q\notin \cU_1$, there exists a $\lambda\colon\Theta\to\cX$ which maps $1\mapsto q$ and $\lambda(0)\not\simeq q$ such that  $\mu^{\cG_m}(q,\lambda) < 0$. 

The idea is to factor $\mu^{\cG_m}(q,\lambda) = \mu^{\cL^{\otimes m}\otimes \cM}(q,\lambda) = m\mu^{\cL}(q,\lambda) + \mu^{ \cM}(q,\lambda)$, pick a $\lambda$ which makes $\mu^{\cL}(q, \lambda) $ negative, and pick $m$ big enough to cancel the contribution of $\mu^{ \cM}(q,\lambda)$. In order to perform this strategy uniformly on $\cX\smallsetminus \cU_1$, we use the normalized Hilbert-Mumford function.

For each point $q\not\in \cU_1$ there is a morphism $\lambda\colon\Theta\to \cX$ such that $\lambda(1)=q\not\simeq\lambda(0)$. Then the results from \Cref{lm:bounded} apply to points in $\cX\smallsetminus\cU_1$ and the normalized Hilbert-Mumford function $\frac{\mu^{\cG_m}(q,\lambda)}{|\lambda|}$. So let
 
\begin{align*}d=\max_{q \notin \cU_1} \left\{ \inf_{\lambda\colon 1 \to q} \frac{\mu^{\cL}(q,\lambda)} {|\lambda|} \right\},\text{ } e=\max_{q\notin \cU_1} \left\{\inf_{\lambda\colon 1 \to q} \frac{\mu^{\cM^{\vee}}(q,\lambda)}{|\lambda|} \right\}
    =\max_{q\notin \cU_1} \left\{\sup_{\lambda\colon 1 \to q} \frac{\mu^{\cM}(q,\lambda)}{|\lambda|} \right\}.
\end{align*}
Observe that with the reductions that we performed, we have $d<0$ and $e<\infty$.

Let $m_0>0$ be an integer such that $m_0\cdot d + e < 0$. For every $q\not \in \cU_1$, let $\lambda_{\min}\colon\Theta\to \cX$ be a morphism which maps $1\mapsto q$ and such that the normalized Hilbert-Mumford function in $q$ reaches its minimum. So
\begin{align*}
    \frac{\mu^{\cG_m}(q,\lambda_{\min})}{|\lambda_{\min}|} = \frac{\mu^{\cL^{\otimes m}}(q,\lambda_{\min)}}{|\lambda_{\min}|} + \frac{\mu^{\cM}(q,\lambda)}{|\lambda|} \le md + e <0
\end{align*}
for every $m\geq m_0$. This shows that $q\notin \cV_m$.
From now on we assume $m\geq m_0$, and we use $\cV$ to denote $\cV_m$.

    \noindent
    \textbf{Step 3.} We show that if $q\in \cU_1\smallsetminus\cU_2$ then $q\notin \cV$.

    Indeed, as $\cU_1\to \cU_1\to U_1$ is cohomologically affine and the line bundle is a character of a group, we can check if a point in $\cU_1$ lies in $\cU_2$ using the Hilbert-Mumford criterion. In particular, if $q\in \cU_1\smallsetminus\cU_2$ there is a morphism $\lambda:\Theta\to \cU_1$ sending $1\mapsto q$ such that $\mu^{\cM}(q,\lambda)<0$. Since $\lambda(\Theta)\subseteq \cU_1$, there is an invariant section of $\cL$ not vanishing on $\lambda(\Theta)$. As for $\lambda(\Theta)$ to be in $\cU_1$ the point $q$ cannot be strictly stable, we have $\mu^{{\cL}}(q,\lambda)=0$. Now the desired statement follows by computing the Hilbert-Mumford function: 
    \[\mu^{\cL^{\otimes m}\otimes \cM}(q,\lambda) = m\mu^{\cL}(q,\pi\circ\lambda) + \mu^{\cM}(q,\lambda) = \mu^{\cM}(q,\lambda)<0.\]

    \noindent
    \textbf{Step 4.} We prove that if a point $q\in \cU_2$, then it is also in $\cV$. This concludes the proof.
    
    We need to produce a section of $\cL^{\otimes m}\otimes \cM$ which does not vanish at $q$. Let $p:\cX\to X$ be the good moduli space morphism. 
    Since $q\in \cU_2\subset\cU_1$, up to shrinking the base, there is a section $s$ of $p_*\cL$ which does not vanish at $q$. As $U_1=\Proj(\bigoplus_n p_*(\cL^{\otimes n}))$, there is an open affine subset of the form $U_{1,s}$ such that $q$ maps to $U_{1,s}$. 
    As $\cU_2$ is the relative stable locus for $\cM=\Bbbk_\chi$,
    and since $\cX$ is a global quotient, from \Cref{lemma_relative_GIT_for_gms_can_be_checked_using_HM} there is $f\in \oH^0(\cX_{s},\cM)$ such that $f(q)\neq 0$.
    From \cite{Liu_AG_book}*{\S 5.1.25}, which is stated for schemes but holds for algebraic stacks as well, there is a section of $\cL^{\otimes m}\otimes\cM$ on the whole $\cX$ lifting $f\otimes s^n$. This section will not vanish at $q$, as desired. 
 \end{proof}
 \subsection{Reduction results for quotient stacks}
 The following lemma is the standard tool for passing from projective GIT to affine GIT. \begin{Lemma}\label{lemma_git_for_ample_linearization_is_contained_in_the_affine_git_for_the_cone}
    Let $W$ be a quasiprojective scheme with an action of $G$, and $L$ a $G$-linearized ample line bundle. For $d>0$, consider the character \[\chi:G\times \Gm\to \Gm, \text{ }(g,t)\mapsto t^{-d}.\]
    This induces a $G\times \Gm$-equivariant line bundle $\Bbbk_\chi$, and up to replacing $L$ with $L^{\otimes m}$ for $m>0$,
    \[
    [W(L)^{ss}/G] \subseteq [\spec(\bigoplus_{i\ge 0} \oH^0(L^{\otimes i}))(\Bbbk_{\chi})^{ss}/G\times \Gm].
    \]
    Moreover, if $W$ is either affine or projective, the inclusion above is an isomorphism, and if $W=\spec(A)$ is affine we also have 
    \[
    [\spec(A)(L)^{ss}/G] \cong [\spec(\bigoplus_{i\in \bZ} \oH^0(L^{\otimes i}))(\Bbbk_{\chi})^{ss}/G\times \Gm].
    \]
\end{Lemma}
\begin{proof} As the semistable locus of $\Bbbk_\chi$ agrees with the one of $\Bbbk_\chi^{\otimes k}$ for every $k\ge 1$, we assume $d=1$.
    Without loss of generality, we can assume that $L$ is generated in degree one and very ample, and also the ring $R_L=\oplus_{n\geq 0} \oH^0(W,L^{\otimes n})^G$ is generated in degree one.
    
    Let $\pi:P_L\to W$ be the $\Gm$-torsor associated to the $G$-linearized line bundle $L$ on $W$, so that there is a $G$-action on $P$ for which $\pi$ is $G$-equivariant. Let $\Bbbk_\chi$ be the $\Gm$-linearized line bundle on $P_L$ associated to the character $\chi$. Then $\pi^*L=\Bbbk_\chi$, and we can regard the latter as a $G\times\Gm$-linearized line bundle. 
    
    We claim that $P_L(\Bbbk_\chi)^{ss}=\pi^{-1}(W(L)^{ss})$, where the first semistable subset is taken with respect to the $G\times \Gm$-action, and the second semistable subset with respect to the $G$-action. This is equivalent to saying that $[W(L)^{ss}/G]\simeq [P_L(\Bbbk_\chi)/G\times \Gm]$.

    The claim follows from the equality $\pi^*\colon \oH^0(W,L)^{G} \simeq \oH^0(P_L,\Bbbk_\chi)^{G\times \Gm}$, which follows from the fact that $\pi^*\colon \oH^0(W,L) \simeq \oH^0(P_L,\Bbbk_\chi)^{\Gm}$ is an isomorphism of $G$-representations.

    The same argument applies to $\Proj(R_L)$ with the line bundle $\cO(1)$, where
    \[ R_L:= \oplus_{n\geq 0} \oH^0(W,L^{\otimes n}).\]
    In this particular case the $\Gm$-torsor over $\Proj(R_L)$ is isomorphic to $\spec(R_L)\smallsetminus\{0\}$. From \cite{stacks-project}*{\href{https://stacks.math.columbia.edu/tag/01PZ}{Tag 01PZ}} and \cite{stacks-project}*{\href{https://stacks.math.columbia.edu/tag/01PZ}{Tag 01Q1}}, the induced morphism $\iota:W\to \Proj(R_L)$ is an open embedding with $\iota^*\cO(1)\simeq L$ and we have a cartesian diagram
    \[
    \begin{tikzcd}
        P_L \ar[r, "j"] \ar[d, "\pi"] & \spec(R_L)\smallsetminus\{0\} \ar[d, "\rho"] \\
        W \ar[r, "\iota"] & \Proj(R_L)
    \end{tikzcd}
    \]
    with $\rho^*\cO(1)\simeq \Bbbk_\chi$. By cartesianity, also $j$ is an open embedding.

    We claim that any point $x$ in $P_L(\Bbbk_\chi)^{ss}$ is also $\Bbbk_\chi$-semistable in $\spec(R_L)$. For this, we have to study the restriction homomorphism
    \[ j^*\colon \oH^0(\spec(R_L),\Bbbk_\chi)^{G\times\Gm} \longrightarrow \oH^0(P_L,\Bbbk_\chi)^{G\times \Gm}.\]
    Observe that
    \[
        \oH^0(\spec(R_L),\cO)\simeq R_L = \oplus_{n\geq 0} \oH^0(W,L^{\otimes n}),
    \]
    hence the composition
    \begin{align*}
        \oH^0(W,L) \hookrightarrow \oH^0(\spec(R_L)\smallsetminus\{0\},\Bbbk_\chi)^{G\times \Gm} \overset{j^*}{\longrightarrow} \oH^0(P_L,\Bbbk_\chi)^{G\times\Gm} \simeq \oH^0(W,L)^G
    \end{align*}
    is necessarily surjective. From this it follows that if there is a $G\times \Gm$-invariant section of $\Bbbk_\chi$ over $P_L$ which does not vanish in $x$, then there is also a $\Gm\times G$-invariant section of $\Bbbk_\chi$ over $\spec(R_L)\smallsetminus\{0\}$ which does not vanish over $x$. This proves that
    \[
    [W(L)^{ss} /G] \simeq [P_L(\Bbbk_\chi)^{ss} / G\times\Gm ] \overset{\text{open}}{\hookrightarrow} [\spec(R_L)(\Bbbk_\chi)^{ss}/G\times\Gm],
    \]
    since the origin is never semistable for this linearization.

    If $W$ is projective or affine, the open embedding $\iota$ is an isomorphism, which implies that $P_L\simeq \spec(R_L)\smallsetminus\{0\}$, hence
    \[ [W(L)^{ss}/G] \simeq [\spec(R_L)(\Bbbk_\chi)^{ss}/G\times\Gm].\]
    If $W$ is affine, we have that also $P_L$ is affine, hence $P_L\simeq \spec(\cO_{P_L}(P_L))$. By definition $P_L$ is the relative spectrum of the $\cO_W$-algebra $\oplus_{n\in \bZ} L^{\otimes n}$, hence
    \[ \oH^0(P_L,\cO_{P_L}) \simeq \oH^0(W, \oplus_{n\in \bZ} L^{\otimes n}) \simeq \oplus_{n\in \bZ} \oH^0(W,L^{\otimes n}),\]
    from which we deduce $[P_L(\Bbbk_\chi)^{ss}/G]\simeq [\spec(\oplus_{n\in \bZ} \oH^0(W,L^{\otimes n})(\Bbbk_\chi)^{ss})/G\times\Gm]$.
\end{proof}
The following lemma allows us to pass from relative GIT for a line bundle to relative GIT for a character, which is better behaved, as we can use the affine Hilbert-Mumford criterion, see \Cref{lemma_relative_GIT_for_gms_can_be_checked_using_HM}.
\begin{Lemma}\label{lemma_cone_construction_over_X}
    Let $\cX=[W/G]$ be a quotient stack by a reductive group $G$, admitting a good moduli space $\cX\to X$ and let $\cL$ be a line bundle on $\cX$. Let $\chi:G\times \Gm\to \Gm$ be the character $(g,t)\mapsto t^{-1}$.
    Then there is a scheme $C$ with an action of $G\times \Gm$ such that
\begin{itemize}
    \item there is an isomorphism $[C(\Bbbk_\chi)^{ss}_X/G\times \Gm]\cong \cX(\cL)^{ss}_X$, 
    \item  there is an isomorphism $[C/G\times \Gm]\cong \cX$, and 
    \item the inclusion $C(\Bbbk_\chi)^{ss}_X\subseteq C$ induces $\cX(\cL)^{ss}_X\to \cX$.
\end{itemize}
    In particular, for every Zariski cover $\spec(A)\to X$, 
        \[C(\Bbbk_\chi)^{ss}_X\times_X\spec(A)=(C\times_X\spec(A))(\Bbbk_\chi)^{ss}.\]
\end{Lemma}
\begin{proof}

    With the same notation of the proof of \Cref{lemma_git_for_ample_linearization_is_contained_in_the_affine_git_for_the_cone}, we can pick $C=P_\cL$ the $\Gm$-torsor associated to the line bundle $\cL$, so that $[C/G\times\Gm]\simeq\cX$. 

    For every Zariski open cover $\spec(A) \to X$, we have that $W\times_X \spec(A)$ is affine, because $\pi:W\to X$ is affine as $W\to [W/G]$ is a $G$-torsor and $[W/G]\to X$ is cohomologically affine. In particular, this implies that the pullback of $\cL$ to $W\times_X\spec(A)$ is ample from \cite{stacks-project}*{\href{https://stacks.math.columbia.edu/tag/0890}{Tag 0890}} and the fact that, over an affine scheme, each coherent module is globally generated. 

    Furthermore, we can assume that $\cX(\cL)^{ss}_{X}\times_X \spec(A)$ is isomorphic to $[(W\times_X \spec(A))(L)^{ss} /G]$, where $L$ denotes the $G$-equivariant line bundle on $W\times_X \spec(A)$ which descends to the pullback of $\cL$: indeed, by definition of relative semistable locus, for every point $p\in \cX(\cL)^{ss}_{X}\times_X \spec(A)$, up to shrinking the base we can always assume that there exists a section of the pullback of $\cL$ that does not vanish on $p$.
    
    By \Cref{lemma_git_for_ample_linearization_is_contained_in_the_affine_git_for_the_cone} we get that $[(W\times_X \spec(A))(L)^{ss} /G]$ is isomorphic to $[(C\times_X \spec(A))(\Bbbk_\chi)^{ss}/G\times\Gm]$.
    On the other hand, by definition, the latter is equal to $[C(\Bbbk_\chi)^{ss}_X \times_X \spec(A) / G\times\Gm]$ as desired.
\end{proof}
\subsection{Relative semistable loci}\label{sec:relative semistable loci}
We end this subsection giving a few conditions under which an open substack of an algebraic stack is the semistable locus of a line bundle.
\begin{Lemma}\label{lemma_condition_when_X_dm_is_the_ss_locus_for_a_vb_when_X_is_normal}
    Let $\cX$ be an algebraic stack with a good moduli space $\cX\to X$, and let $\cX'\subseteq \cX$ be a dense open substack that admits a good moduli space $X'$ projective over $X$. Let $\cL$ be a line bundle on $\cX$ which, once restricted to $\cX'$, descends to a line bundle $L_{X'}$ on $X'$, that is ample over $X$. Assume finally that the complement of $\cX'$ in $\cX$ has codimension at least 2, and that $\cX$ is normal. Then $\cX'=\cX(\cL)^{ss}_X$.
\end{Lemma}
\begin{proof}The statement is \'etale local over $X$ so we can assume that $\cX=[\spec(A)/G]$ and $X=\spec(A^G)$, so $L'$ is ample.
We will use the following chain of equalities: \begin{equation}\label{equality_sections}
    \oH^0(\cX,\cL^{\otimes m}) = \oH^0(\cX',\cL^{\otimes m}|_{\cX'}) = \oH^0(X',L_{X'}^{\otimes m})
\end{equation} where the first equality follows since $\cX$ is normal and the complement of $\cX'$ has codimension at least 2.

    The inclusion $\cX'\subseteq \cX(\cL)^{ss}_X$ can be proved as follows. For every $p\in \cX'$ there is an $m\gg 0$ and a section of $L_{X'}^{\otimes m}$ which does not vanish at $p$. On the other hand from \Cref{equality_sections}, such a section extends to a section $s\in\oH^0(\cX,\cL^{\otimes m})$.
    The locus where $s$ does not vanish is of the form $[\spec(A)\smallsetminus V(\pi^*s)/G]$ where $\pi:\spec(A)\to \cX$ is the quotient map. Since $\spec(A)\smallsetminus V(\pi^*s)$ is affine, also $[\spec(A)\smallsetminus V(\pi^*s)/G]$ is cohomologically affine. This proves $\cX_\dm\subseteq \cX(\cL)^{ss}_X$.

    From \cite{Alp}*{Theorem 11.5}, there is a morphism \[\xi:\cX(\cL)^{ss}_X\to \Proj\left(\bigoplus \oH^0(\cX,\cL^{\otimes m})\right)\] and from \Cref{equality_sections} and since $L_{X'}$ is ample, $\Proj(\bigoplus \oH^0(\cX,\cL^{\otimes m}))= X'$. In particular, as $\cX'\to X'$ is surjective, also $\cX(\cL)^{ss}_X\to X'$ is surjective, so $\xi$ is a good moduli space from \cite{Alp}*{Theorem 11.5}. Both $\cX'\to X'$ and $\cX(\cL)^{ss}_X\to X'$ are cohomologically affine, the stack $\cX(\cL)^{ss}$ is normal and the complement of $\cX'$ in $\cX(\cL)^{ss}_X$ is at least 2. So $\cX'=\cX(\cL)^{ss}_X$ as desired.
\end{proof}
\begin{Cor}\label{cor_when_X_dm_is_ss_locus_for_a_lb}
    Let $\cX$ and $\cX'$ be as in \Cref{lemma_condition_when_X_dm_is_the_ss_locus_for_a_vb_when_X_is_normal}, except that $\cX'$ is assumed to be Deligne--Mumford and $\cX$ is not assumed to be normal. If $\cX(\cL)^{ss}_X$ is dense, we have $\cX'=\cX(\cL)^{ss}_X$.
\end{Cor}
Observe that if the properly stable locus $\{x\in X:\pi^{-1}(x)$ is Deligne--Mumford$\}$ is dense then $\cX(\cL)^{ss}_X$ is dense, where $\pi:\cX\to X$ is the good moduli space map.
\begin{proof}The statement is  \'etale local over $X$ so we can assume that $\cX=[\spec(A)/G]$ and $X=\spec(A^G)$, so $L'$ is ample. Observe that $\cX(\cL)^{ss}_X\subseteq \cX'$. Indeed,
consider $\nu:\cX^n\to \cX$ be the normalization. Then $\nu^{-1}(\cX(\cL)^{ss}_X)\subseteq \cX^n(\nu^*\cL)^{ss}_X$ since $\nu $ is affine, and $\cX^n(\nu^*\cL)^{ss}_X = (\cX')^n = \nu^{-1}(\cX')$ where the first equality follows from \Cref{lemma_condition_when_X_dm_is_the_ss_locus_for_a_vb_when_X_is_normal} (the assumpions of \Cref{lemma_condition_when_X_dm_is_the_ss_locus_for_a_vb_when_X_is_normal} are still satisfied as $\nu$ is the normalization). Then $\cX(\cL)^{ss}_X =\nu(\nu^{-1}(\cX(\cL)^{ss}_X))\subseteq \nu(\nu^{-1}(\cX'))=\cX'$. So $\cX(\cL)^{ss}_X$ is an open substack of $\cX'$, so it is Deligne--Mumford, and from \Cref{lemma_git_for_ample_linearization_is_contained_in_the_affine_git_for_the_cone} the coarse moduli space of $\cX(\cL)^{ss}_X$ is projective, hence $\cX(\cL)^{ss}_X=\cX'$.
\end{proof}

 \subsection{Luna slice}\label{subsection_luna_slice}
 For the convenience of the reader, we briefly recall some classical results on GIT, which are well-known.
    Let $X$ be a scheme of finite type with an action of a reductive group $G$, and let $U\subseteq X$ be the semistable locus with respect to a $G$-equivariant line bundle $\cL$ on $X$. If $x\in U$ is a closed point with stabilizer $G_x$, we have a cartesian diagram
    \[
    \xymatrix{
    [\spec(A)/G_x]\ar[r]\ar[d] & [U/G]\ar[d] \\ 
    \spec(A^{G_x})\ar[r]^-\alpha  & U/\!\!/_{\cL} G,
    }
    \]
    where $\alpha$ is an \'{e}tale neighborhood of the image of $x$ in the good moduli space \cite{Alp21}*{Corollary 6.7.3}. If $x$ is a smooth point with tangent space $T_x$, there is an action of $G_x$ on $T_x$. 
    Assume that the orbit of $x$ is closed: then from Matsushima’s Theorem the group $G_x$ is reductive, so we can decompose the $G_x$-representation \[T_x=T_x(G\cdot x)\oplus N_x\] where $T_x(G\cdot x)$ is the tangent space at $x$ of the orbit $Gx$ of $x$, and $N_x$ is $G_x$-invariant. Then it easily follows from  \cite{Alp21}*{Theorem 6.7.5, "Luna’s \'{E}tale Slice Theorem"} applied to $[\spec(A)/G_x]$ that we have a cartesian diagram
    \[
    \xymatrix{[N_x/G_x]\ar[d] & [\spec(A)/G_x]\ar[l]\ar[d] \\
    N_x/\!\!/G_x & \spec(A^{G_x}) \ar[l]_-\beta ,}
    \]
    where $\beta$ is an \'{e}tale neighborhood of $0\in N_x$
\begin{Def}
    With the notations above, we define $[N_x/G_x]$ to be a \textit{Luna slice} for $[U/G]$ at $x$.
\end{Def}
In other terms, the map $[N_x/G_n]\to N_x/\!\!/G_x$ gives a local description of $\pi:[U/G]\to U/\!\!/_\cL G$ in a neighbourhood of $\pi(x)$.

\section{Quasimaps and $\epsilon$-stability}\label{sec:epsilon qmaps}
In this section we construct the stack of stable quasimaps and $\epsilon$-stable quasimaps with target a quotient stack $\cX$ with a projective good moduli space $X$ and an open dense substack $\cX_\dm$ which is Deligne-Mumford and proper; see \Cref{assumptions_epsilon_qmaps} for the precise hypotheses on the stack.

We do so by first showing in \Cref{prop_ideal_sheaf_unstable_locus_pulls_back}
that the target stack $\cX$ can be embedded as an open substack of the quotient stack $[\spec(R)/G\times \Gm]$, where $\spec(R)$ is an affine scheme of finite type and $G$ is a reductive group. This embedding enjoys some special properties: there are characters $\rho$ and $\xi$ such that $\cX$ is equal to the GIT semistable locus of $\spec(R)$ with respect to the linearization $\Bbbk_\xi$, and $\cX_\dm$ is the semistable locus with respect to $\Bbbk_{\xi^N\otimes \rho}$ for any $N\gg 0$. 

This embedding allows us to relate our stack of quasimaps with target $\cX$ to the well studied quasimaps to $[\spec(R)/G\times\Gm]$ of \cite{ciocan2014stable}, linearized by the character $\xi^N\otimes \rho$. We show that when picking $N\gg 0$, the $\epsilon$-stability condition for the classical quasimaps to $[\spec(R)/G\times\Gm]$ with respect to $\xi^N\otimes \rho$ forces the $\epsilon$-stable quasimaps to factor through the open substack $\cX$. Once this property is secured, is easy to see that the two stacks of $\epsilon$-stable quasimaps, those with target $\cX$ and the classical ones, are isomorphic (\Cref{thm:relative quasimaps are quasimaps}).
Geometric properties for $\epsilon$-stable quasimaps to $\cX$ are then deduced from the known properties in the classical case (\Cref{thm:epsilon quasimaps are dm}).

An additional argument yields that also the stack of stable quasimaps has the same geometric properties as the classical one, and it is a limit of the $\epsilon$-stable one for $0<\epsilon\ll 1$ (\Cref{thm:quasimaps are dm}).
\subsection{The stack of $\epsilon$-stable quasimaps}\label{subsection_epsilon_stable}
From now on we will adopt the following assumption.
\begin{Assumption}\label{assumptions_epsilon_qmaps}
    We will denote by $\cX$ an algebraic stack having a projective good moduli space $\cX\to X$, with an ample line bundle $L$ on $X$. Further, we will assume that:
    \begin{enumerate}
        \item there exists a line bundle $\cL$ on $\cX$ such that $\cX_\dm := \cX(\cL)^{ss}_X$ is Deligne--Mumford, where the latter denotes the relative $\cL$-semistable locus of $\cX$ over the good moduli space $X$, and 
        \item there is an isomorphism $\cX\cong [W/G]$ for a scheme $W$ and a reductive group $G$.
    \end{enumerate}
    Up to replacing $\cL$ with $p^*L^{\otimes m}\otimes \cL$, we can further assume that $\cL$ descends to a very ample line bundle on the coarse moduli space of $\cX(\cL)^{ss}_X$. 
\end{Assumption}
\begin{Remark}\label{rmk_if_W_is_Qfactorial_then_any_open_is_the_ss_locus_of_a_linebundle}
    Condition (1) is automatically satisfied if $\cX$ is smooth and contains a proper Deligne--Mumford open substack by \cite{Alp}*{Theorem 11.14.}, and the same proof goes through if $\cX=[W/G]$ and $W$ is $\bQ$-factorial. See also \Cref{cor_when_X_dm_is_ss_locus_for_a_lb} for other cases when Condition (1) holds.
\end{Remark}
The following definitions already appeared \cite{orbifold_qmap_theory}, albeit only for quotient stacks $\cX=[W/G]$ with $W$ affine.
\begin{Def}\label{def:quasimap}
    An $n$-marked \emph{quasimap} of genus $g$ over $S$ is  a pair \[(\phi\colon \sC\to \cX, \Sigma_1,\ldots, \Sigma_n)\] such that
    \begin{enumerate}
        \item $(\sC \to S,\Sigma_1,\ldots,\Sigma_n)$ is a twisted, $n$-marked curve of genus $g$ as in \cite{AV_compactifying};
        \item the morphism $\phi:\sC\to \cX$ is representable and maps the generic points, the marked gerbes and the nodal locus of $\sC\to S$ to $\cX_\dm\subseteq \cX$.
    \end{enumerate}
    If $S=\spec(k)$, we will just say an $n$-marked quasimap of genus $g$.
\end{Def}   
We will denote the coarse moduli space of a twisted curve $\sC$ by $C$, and the image of the $\Sigma_i$ in $C$ by $p_i$.
\begin{Def}\label{def_beta}
    Given a quasimap $(\phi\colon \sC \to \cX,\Sigma_1,\ldots,\Sigma_n)$ over a geometric point, the \emph{class} of the quasimap is the homomorphism
    \[\phi_*[\sC] \colon \Pic(\cX) \longrightarrow \bQ,\quad L \longmapsto \deg(\phi^*L).\]
\end{Def}

To define $\epsilon$-stability, we need the notion of length of a point; see also \cite[7.1.2]{ciocan2014stable}.
\begin{Def}
    Given a morphism $\phi\colon \cC\to \cX$ from a twisted curve $\cC$, given a smooth point $x\in \cC$ we define the \textit{length} of $x$ with respect to $\cL$ as the
    length at $x$ of the cokernel of 
    \[
    \phi^{-1}(I)\to \cO_\cC
    \]
    where $I$ is the ideal sheaf of the $\cL$-unstable locus of $\cX$ over $X$. The length of a point $x\in \cC$ is denoted by 
    $\ell_{(\cX,\cL)}(x)$. We will drop any reference to the pair $(\cX,\cL)$ when they are clear.
\end{Def}
\begin{Remark}
    The definition above works with our assumption that, in particular, $\cL$ descends to a relatively very ample line bundle for $X_\dm \to X$. One can define the notion of length also when $\cL$ is relatively ample but not very ample, just as in \cite[Definition 7.1.1]{ciocan2014stable}. On the other hand, up to rescaling $\cL$, one can always assume that this is the case, as explained in \cite[loc. cit.]{ciocan2014stable}.
\end{Remark}
Recall that, as any line bundle $L$ on $\cC$ can be regarded as a $\mathbb{Q}$-line bundle on $C$, we have well defined notions of ampleness and degree.
\begin{Def}\label{def_epsilon_stabe_qmaps_for_us}
    Fix a real number $\epsilon>0$. A quasimap $(\phi:\sC\to \cX, \Sigma_1,\ldots, \Sigma_n) $ is $\epsilon$-stable if
    \begin{enumerate}
        \item[(Stab)] the line bundle $\omega_\sC(\Sigma_1+\ldots+\Sigma_n)\otimes \phi^*(p^*L^3\otimes \cL^{\epsilon})$ is ample, and 
        \item[(Sing)] for every $x\in \cC$ such that $\phi(x)\not\in \cX_\dm$, we have
        \[
        \ell(x)<\frac{1}{\epsilon}.
        \]
    \end{enumerate}
\end{Def}
The abbreviations (Stab) and (Sing) stand for singularity condition and stability condition.

\begin{Remark}
    Stability can also be defined for $\cL$ that descends to a relatively ample but not very ample line bundle on $X_\dm$: given $m$ such that $\cL^m$ descends to a relatively very ample line bundle, we simply set $\epsilon$-stable quasimaps for $\cL$ to be $\epsilon/m$-stable quasimaps for $\cL^m$; see also \cite[\S 7.1.4]{ciocan2014stable}. In general, replacing $\cL$ by a multiple of it does not change the $\epsilon$-stable quasimaps, granted that one rescales $\epsilon$ accordingly. 
    
    The power $3$ in condition (Stab) can be replace with any $n\ge 3$.
\end{Remark}
\begin{Remark}\label{rmk_pull_back_of_L_is_ample}
    Since $\cL$ is assumed to descend to a very ample line bundle on $\cX_\dm$, it is immediate to check that the pull-back $\phi^*\cL$ via a quasimap $\phi\colon \sC\to \cX$ is nef.
\end{Remark}
\begin{Def}
We define the category $\cQ_{g,n}^\epsilon(\cX,\cL,\beta)$ of $\epsilon$-stable quasimaps to $\cX$ with respect to $\cL$ of class $\beta$ as the following category fibered in groupoids over the big \'{e}tale site of schemes over $k$:
\begin{enumerate}
    \item for any scheme $S$, the objects are the $n$-marked quasimaps to $\cX$ of genus $g$ over $S$ whose geometric fibers are $\epsilon$-stable and of class $\beta$, and
    \item the morphisms $(\phi\colon\sC\to\cX,\Sigma) \longrightarrow (\phi'\colon\sC'\to\cX,\Sigma')$ over $S$ consist of an isomorphism $\alpha\colon\sC\simeq \sC'$ such that $\alpha^{-1}(\Sigma'_i)=\Sigma_i$ and an isomorphism $\gamma\colon\phi'\circ\alpha \simeq \phi$.
\end{enumerate}
\end{Def}
The following is one of the two main theorems of this section:
\begin{Teo}\label{thm:epsilon quasimaps are dm}
    The stack $\cQ_{g,n}^\epsilon(\cX,\cL,\beta)$ is a proper Deligne--Mumford stack.
\end{Teo}

Our strategy is to reduce to the classical case $\cX=[\spec(R)/G\times\Gm]$ of \cite{ciocan2014stable} as follows. We first prove that we can enlarge the inclusion $\cX_\dm\hookrightarrow \cX$ as
\[
\cX_\dm\hookrightarrow \cX\hookrightarrow [\spec(R)/G\times\Gm]
\]
for a certain reductive group $G$, and for a certain ring $R$, such that that $\cX$ is the semistable locus for a character $\xi\colon G\to \Gm$ and $\cX_\dm$ is the the semistable locus for another character $\xi^N\otimes \rho$, up to picking $N$ large enough. This first step is achieved in \Cref{prop_ideal_sheaf_unstable_locus_pulls_back}.

Then, we show that for a given $\epsilon$, one can pick $N$ large enough so that any $\epsilon$-quasimap with target $[\spec(R)/G\times\Gm]$, and with respect to the character $\xi^N\otimes \rho$, must map to $\cX$. In other terms, if a quasimap was not to land in $\cX$, then one can increase $N$ so that such a quasimap would have a basepoint of length greater than $\epsilon$.
In the process we will prove that if an $\epsilon$-stable quasimap with respect to the character $\xi^N\otimes \rho$ and with target $[\spec(R)/G\times\Gm]$ lands in $\cX$, then it is also $\epsilon$-stable in the sense of \Cref{def_epsilon_stabe_qmaps_for_us}. This will conclude the argument.
\begin{Remark}\label{rmk_L_is_given_by_a_line_bundle}
    From \Cref{lemma_cone_construction_over_X} one can assume, up to replacing:
    \begin{enumerate}
        \item the scheme $W$ with another scheme affine over $X$,
        \item the group $G$ with $G\times \Gm$, and 
        \item $\cL$ with a multiple of it, 
    \end{enumerate}
 that the line bundle $\cL$ is given by $\Bbbk_\rho$ for a certain character $\rho\colon G\to \Gm$. From now on, we will assume that this is the case. 
\end{Remark}
\subsection{Construction of the embedding}
We begin by embedding $\cX$ in an algebraic stack of the form $[\spec(R)/G\times\Gm]$.
\begin{Prop}\label[Prop]{prop_ideal_sheaf_unstable_locus_pulls_back}
    There is an affine scheme of finite type $\spec(R)$ endowed with an action of $G\times\Gm$ and two characters $\xi$ and $\rho$ of $G\times\Gm$, such that for every $N$ sufficiently large, we have:
    \begin{enumerate}
        \item $\cX$ is an open substack of $[\spec(R)/G\times\Gm]$ which agrees with the $ \Bbbk_\xi$-semistable locus of $[\spec(R)/G\times\Gm](\Bbbk_\xi)$,
        \item $\cX_\dm$ agrees with the $\Bbbk_{\xi^N\otimes \rho}$-semistable locus of $[\spec(R)/G\times\Gm]$, and
        \item\label{condition_in_itemize_unstable_locus_pulls_back}  up to replacing $\cL$ and $L$ with some positive multiples of them, we can arrange the ideal sheaf of the $\Bbbk_{\xi^N\otimes\rho}$-unstable locus in $[\spec(R)/G\times\Gm]$ to restrict to the ideal sheaf of the $\cL$-unstable locus of $\cX$ over $X$.
    \end{enumerate}
\end{Prop}
We start with the following preparatory lemma.
\begin{Lemma}\label[Lemma]{lm:embedding}
    Let $\cX=[W/G]$ be a quotient stack with a projective good moduli space $\cX\to X$, such that $G$ is reductive. Then there is a projective scheme $\overline{W}$ with a $G$ action and an ample $G$-linearized line bundle $M$ on $\overline{W}$ such that there are two open embeddings \[\cX\hookrightarrow [\overline{W}/G]\hookrightarrow \left[\spec\left(\bigoplus_n \oH^0(\overline{W},M^{\otimes n})\right)/G\times \Gm\right].\] Moreover, we have that $\cX = [\overline{W}(M)^{ss}/G]$, i.e. the open substack $\cX$ of $[\overline{W}/G]$ coincides with the $M$-semistable locus of $\overline{W}$. Finally, the character \begin{center}$\xi\colon G\times \Gm\to \Gm$, $(g,t)\mapsto t^{-1}$\end{center}is such that \[\cX=\left[\spec\left(\bigoplus_n \oH^0(\overline{W},M^{\otimes n})\right)/ G\times\Gm\right](\Bbbk_{\xi})^{ss}.\]
\end{Lemma}
The proof consists of two parts. In the first two steps we argue that, as $\cX$ is a global quotient with a projective good moduli space, then $\cX$ comes from GIT, i.e. there is a projective $G$-variety $\overline{W}$ and a $G$-linearized ample line bundle $\cL$ such that $W=\overline{W}(\cL)^{ss}$.
The proof in \cite{teleman2000quantization}*{Lemma 6.1} goes through in the previous setting. This is also (essentially) the argument used in \cite{AHH}*{proof of Proposition 5.11}. To conclude the argument we use \Cref{lemma_git_for_ample_linearization_is_contained_in_the_affine_git_for_the_cone}.
\begin{proof}
We will denote by $p\colon W\to X$ the composition $W\to [W/G]=\cX\to X$.

    \underline{Step 1}. We embed $G$-equivariantly $W$ into a vector bundle $F\to X$.
    
    As $L$ is very ample, there are global sections $s_0,\ldots, s_n$ of $L$ such that $X= \bigcup_i D(s_i)$ with $p^{-1}(D(s_i))=\spec(A_i)$ and $D(s_i)=\spec(A_i^G)$. This follows since $W\to X$ is affine.

    As $p\colon W\to X$ is of finite type, for each $i$ there exist elements $f_{i1},\ldots, f_{i\ell}$ of $A_i$ such that \[A_i= A_i^G[f_{i1},\ldots, f_{i\ell}].\] Take $k$ sufficiently large so that for each $i$ and $j$, we have that $f_{ij}\otimes p^*s_i^k$ extends to global sections of $p^*L^{\otimes k}$ on $W$. 
    As $G$ is reductive, the $G$-representation $\oH^0(W,p^* L^{\otimes k})$ splits into finite dimensional, irreducible $G$-representations. Let $V\subset \oH^0(W,p^* L^{\otimes k})$ be a finite dimensional $G$-representation containing $f_{ij}\otimes p^*s_i^k$ for all $i,j$.

    We claim that there is a surjective morphism $\Sym^\bullet (V\otimes L^{-k}) \to p_*\cO_W$ of $\cO_X$-algebras. Granted this, by taking the relative spectrum over $X$ we obtain a closed, $G$-equivariant embedding $W \hookrightarrow F$ over $X$, with $F\to X$ a vector bundle.
    
    To define the morphism of $\cO_X$-algebras, we only have to define a map $V\otimes L^{-k} \to p_*\cO_W$: locally on the $D(s_i)$, every section of $L^{-k}$ is of the form $g\cdot s^{-k}_i$ with $g$ an element of $A_i^G$, hence we define the morphism by sending an element \[f_{ij}\otimes p^*s_i^k \otimes g\cdot s_i^{-k} \mapsto gf_{ij}.\]One can check that this morphisms glue. 
    
    Surjectivity of $\Sym^\bullet (V\otimes L^{-k}) \to p_*\cO_W$ can be checked locally, and by construction the polynomial algebra $A_i^G[f_{i1},\ldots, f_{i\ell}]$ surjects onto $A_i$.

    \underline{Step 2}. Take the closure $W\subseteq \overline{W}$ of the composition $W\hookrightarrow F \hookrightarrow \bP_{X}(F\oplus\cO_X)$, and let $E$ be the hyperplane at infinity. Let $H:=\cO_{\overline{W}}(E)$. We argue that \[\overline{W}(\overline{p}^*L^{\otimes N} \otimes H)^{ss}=W\] for $N\ge n_0$ for some $n_0$. 
    
    Let $M:=\overline{p}^*L^{\otimes N} \otimes H$.    
From the Hilbert-Mumford criterion, as any map $\Theta\to [W/G]$ factors via the reduced structure $[W^{red}/G]$, up to replacing $W$ with $W^{red}$ we can assume that $W$ is reduced. Observe that the algebraic stack $[W^{red}/G]$ still admits a projective good moduli space, which is a closed subscheme of $X$. Now the desired statement follows since the semistable locus of $\overline{W}$ \textit{over} $X$ for the line bundle $\cO_{\bP_X(F\oplus \cO_X)}(E
)|_{\overline{W}}$ is $W$; this is proved in \cite{AHH}*{proof of Proposition 5.11, second paragraph}. The $G$-semistable locus for $X$ with respect to the ample line bundle $L$ is $X$, as the action of $G$ on $X$ is trivial. The desired statement now follows from \cite{ER21}*{Proposition 3.18}.

    \underline{Step 3}. We conclude the argument using \Cref{lemma_git_for_ample_linearization_is_contained_in_the_affine_git_for_the_cone}, possibly replacing $M$ with $M^{\otimes N}$ for $N\gg 0$ as in \textit{loc. cit.}

    Take the affine cone $\spec(\oplus_k \oH^0(\overline{W}, M^k))$, so that we have the embedding \[[\overline{W}/G] \hookrightarrow [\spec(\oplus_k \oH^0(\overline{W}, M^k))/G\times \Gm].\] By \Cref{lemma_git_for_ample_linearization_is_contained_in_the_affine_git_for_the_cone}, we have that \[[\spec(\oplus_k \oH^0(\overline{W}, M^k))(\Bbbk_\xi)^{ss}/G\times \Gm] = [W/G] = \cX.\]
\end{proof}
\begin{Remark}
    Observe that, with the notations of \Cref{prop_ideal_sheaf_unstable_locus_pulls_back}, the line bundle $\Bbbk_\xi$ restricts to a positive multiple of $p^*L$ on $W$ and to a positive multiple of what we denoted by $M$ in Step 2 of the proof of \Cref{prop_ideal_sheaf_unstable_locus_pulls_back} on $\overline{W}$. 
\end{Remark}

\begin{Lemma}\label[Lemma]{lm:embedding 2}
With the assumptions and notations of \Cref{lm:embedding}, there is $n_0$ such that for every $N\ge n_0$ we have 
\[
\cX_\dm = [\spec(\oplus_k \oH^0(\overline{W}, M^k))/G\times \Gm](\Bbbk_{\xi^N\otimes \rho})^{ss}
\]
where $\xi\colon G\times \Gm\to \Gm$ sends $(g,t)\mapsto t^{-1}$ and we still denote the composition $G\times \Gm\xrightarrow{\pi_1} G\xrightarrow{\rho} \Gm$ by $\rho$.
\end{Lemma}
\begin{proof}Let $R:=\oplus_k \oH^0(\overline{W}, M^k)$.
    We apply \Cref{prop_stable_locus_of_pi*L^n_otimes_G} with $\cX_1=\cX_2=[\spec(R)/G\times \Gm]$. More explicitly, there is a character (namely, the character $\xi$ of \Cref{lm:embedding}) such that \[[\spec(R)/G\times \Gm](\Bbbk_\xi)^{ss}=\cX.\] Similarly, by \Cref{rmk_L_is_given_by_a_line_bundle}, there is a character $\rho$ such that $\cX_\dm = \cX(\Bbbk_\rho)^{ss}_X$. From \Cref{prop_stable_locus_of_pi*L^n_otimes_G} \[[\spec(R)/G\times \Gm](\Bbbk_{\xi}^{\otimes N}\otimes\Bbbk_\rho)^{ss}=\cX_\dm\]
    for every $N\ge n_0$.
\end{proof}
\begin{Remark}\label{finite_dimension_of_ring_of_invariants}
    Observe that, since $\overline{W}$ is projective, $\oH^0(\overline{W},M^{k})$ has finite dimension. So also the $(G\times \Gm)$-invariants of $\oplus_k \oH^0(\overline{W}, M^k)$ are a finite dimensional vector space.
\end{Remark}
\begin{proof}[Proof of \Cref{prop_ideal_sheaf_unstable_locus_pulls_back}] We will use the notations and assumptions of \Cref{lm:embedding} and \Cref{lm:embedding 2}. In particular, we have:
\begin{enumerate}
    \item a projective variety $\overline{W}$ with an open embedding $W\subseteq \overline{W}$, an ample line bundle $M$ and a map $\overline{p}\colon \overline{W}\to X$, as in \Cref{lm:embedding};
    \item 
     $M=\overline{p}^*L^{\otimes k}\otimes H$, where $H$ is ample over $X$ and admits a $G$-invariant section, which we denote by $h$, given by some multiple the complement of $W\subseteq \overline{W}$; and
    \item two characters $\rho$ and $\xi$ for $G\times \Gm$, so that $\cL=\Bbbk_\rho$ and $\xi$ is as in \Cref{lm:embedding 2}.
\end{enumerate}Moreover, from \Cref{lm:embedding 2} we can replace $N$ with $N'>N$ if necessary.
    Observe then that the sections of $\spec(\bigoplus_k \oH^0(\overline{W},M^d))$ which are $(\Bbbk_{\xi}^{\otimes N}\otimes\Bbbk_\rho)^d$-semi-invariant are the sections of $\oH^0(\overline{W},M^{Nd})$ which are $\rho^d$-semi-invariant. 

 Now, consider the ideal sheaf of the unstable locus $\cX\smallsetminus \cX_\dm$. Recall its construction: given an affine open $U\subseteq X$, the ideal sheaf in $W$ of the preimage of the unstable locus is generated by the $\rho^n$-semi-invariant sections of $\cO_W(p^{-1}(U))$, for every $n\ge 1$. As this ideal is finitely generated, it will be
generated by a sequence of sections $g_{1,U},\ldots,g_{r,U}$.
We can pick $U$ which is the complement of the zero locus of a section $s_U$ of $L^{\otimes k}$, and furthermore we can assume that $p^{-1}(U)$ is schematically dense in $W$. This can be achieved by choosing the section $s_U$ so that its zero locus does not contain the image of any of the associated primes of $W$.

There is $m_0$ such that for every $m\ge m_0$, these
sections will lift to global sections of $\oH^0(\overline{W},M^{\otimes m})$ up to replacing $g_{i,U}$ with $g_{i,U}\otimes \overline{p}^*s_U^{\otimes km}\otimes h^m$.
Moreover, as $s_U$ and $h$ are invariant, the section $g_{i,U}\otimes \overline{p}^*s_U^{\otimes km}\otimes h^m$ is still $\rho^{n_i}$-semi-invariant for some $n_i$,
as it is such when restricted to $p^{-1}(U)$ which is schematically dense in $W$. So we can find sections $\gamma_i$ as above which are
$\rho^{n_i}$-semi-invariant in $\oH^0(\overline{W},M^m)$ for every $m\ge m_0$.

Now, we proceed as follows. First we pick $N$ to be a positive multiple of $m_0$; this is possible from \Cref{lm:embedding 2}. Next, we pick $\gamma_i\in \oH^0(\overline{W},M^{\otimes Nd})$ which is $\rho^d$-semi-invariant; this is again possible up to replacing $g_{i,U}\otimes \overline{p}^*s_U^{\otimes km_0}\otimes h^{m_0}$ with $g_{i,U}\otimes \overline{p}^*s_U^{\otimes km}\otimes h^m$. By how the sections $\gamma_i$ are constructed, the ideal sheaf of the $\Bbbk_{\xi^N\otimes \rho}$-unstable locus of $[\spec(R)/G\times \Gm]$ contains the sections $\gamma_i$, which restrict to generators of the $\cL$-unstable locus of $\cX$.

The other inclusion, namely the fact that the ideal sheaf of the $\cL$-unstable locus contains the restriction of the ideal sheaf of the $\Bbbk_{\xi^N\otimes \rho}$-unstable locus, is clear.
\end{proof}
\subsection{Comparing the quasimaps}
In light of what we proved in the previous subsection, we will make the following.
\begin{Assumption}
We will denote by $n_0$ the positive integer of \Cref{lm:embedding 2}.
    With the assumptions and notations of \Cref{prop_ideal_sheaf_unstable_locus_pulls_back}, we will assume that both the algebra generated by the $\xi^d$-semi-invariants for $d>0$, and one generated by the $(\xi^N\otimes \rho)^d$-semi-invariants for $d>0$, are generated in degree one. Moreover, we assume that the assumptions of \Cref{prop_ideal_sheaf_unstable_locus_pulls_back} point \ref{condition_in_itemize_unstable_locus_pulls_back} hold for every $N\ge n_0$. This is always possible up to replacing $\cL$ and $L$ with some positive multiples of them. 
\end{Assumption}

\begin{Cor}\label[Cor]{cor_length_wrt_curlyX_is_the_same_as_length_wrt_specR} Let $\cC$ be a twisted curve with a quasimap $\cC\to \cX$. Consider the composition \[\cC\xrightarrow{f} \cX\xrightarrow{\iota} [\spec(R)/G\times \Gm]\] with $\iota$ constructed as before. Then for every smooth point $x\in \cC$ which maps not to $\cX_\dm$, we have 
\[
\ell_{\cX,\cL}(x) = \ell_{[\spec(R)/G\times \Gm],\Bbbk_{\xi^N\otimes \rho}}(x)\text{ for every }N\ge n_0.
\]
\end{Cor}
In particular, the quantity $\ell_{[\spec(R)/G\times \Gm],\Bbbk_{\xi^N\otimes \rho}}(x)$ does not depend on $N$ as long as $N\ge n_0$.
\begin{proof}Let $I$ (resp. $J$) the ideal sheaf of the $\cL$-unstable locus (resp. the $\Bbbk_{\xi^{\otimes N}\otimes \rho}$-unstable locus) of $\cX$ over $X$ (resp. of $[\spec(R)/G\times \Gm]$).
    The statement follows from \Cref{prop_ideal_sheaf_unstable_locus_pulls_back} point \ref{condition_in_itemize_unstable_locus_pulls_back}.
\end{proof}
We are finally ready to prove one of the two main results of this section.

\begin{Teo}\label{thm:relative quasimaps are quasimaps} Let $g,n\in \bN$ and let $\beta$ be a class for quasimaps with target $\cX$ and let $\widetilde{\beta}$ be its induced class for the target $[\spec(R)/G\times\Gm]$.
    Then for every $\epsilon>0$ is an integer $N = N(\epsilon,\cX,\cL)$ such that \[\cQ^{\epsilon}_{g,n}([\spec(R)/G\times \Gm],\Bbbk_{\xi^N}\otimes \Bbbk_\rho,\widetilde{\beta}) = \cQ_{g,n}^{\epsilon}(\cX,\cL,\beta).\] 
\end{Teo}
Observe that the right hand side of the equation in \Cref{thm:relative quasimaps are quasimaps} does not depend on $N$. The way to interpret \Cref{thm:relative quasimaps are quasimaps} is as follows: the $\Bbbk_\xi$-semistable locus in $[\spec(R)/G\times\Gm]$ is $\cX$; we need to prove that up to increasing $N$, any quasimap $\cC\to [\spec(R)/G\times\Gm]$ which does not land in $\cX$ will have a basepoint of length greater than $1/\epsilon$, thus forcing $\epsilon$-stable quasimaps to land in $\cX$.
\begin{proof} We begin by replacing $n_0$ with $\max(n_0, \lceil\frac{3}{\epsilon}\rceil)$.
From \Cref{cor_length_wrt_curlyX_is_the_same_as_length_wrt_specR}, for every $N\ge n_0$, any quasimap in $\cQ_{g,n}^{\epsilon}(\cX,\cL,\beta)$,
composed with the open embedding $\cX\to [\spec(R)/G\times\Gm]$, is also a stable quasimap in $\cQ^{\epsilon}_{g,n}([\spec(R)/G],\Bbbk_{\xi^N}\otimes \Bbbk_\rho,\widetilde{\beta})$. We now show the other inclusion for $N\gg n_0$.

The argument is in two steps. First we argue that if an $\epsilon$-stable quasimap with target $[\spec(R)/G\times\Gm]$ factors as \[\cC\to \cX\to [\spec(R)/G\times\Gm],\] then it is in $\cQ_{g,n}^{\epsilon}(\cX,\cL,\beta)$. Then, we argue that up to possibly increasing $N$, any $\epsilon$-stable quasimap $\cC\to [\spec(R)/G\times \Gm]$ factors as $\cC\to \cX\to [\spec(R)/G\times \Gm]$.

Consider an $\epsilon$-stable quasimap \[(\widetilde{\phi}\colon \cC\to [\spec(R)/G\times\Gm])\in \cQ^{\epsilon}_{g,n}([\spec(R)/G\times\Gm],\Bbbk_{\xi^{n_0}}\otimes \Bbbk_\rho,\widetilde{\beta})\] which factors as $\cC\xrightarrow{\phi} \cX\to [\spec(R)/G\times\Gm]$. Since $\Bbbk_{\xi}|_\cX=\pi^*L$ and $\Bbbk_{\rho}|_\cX=\cL$ we have
\[
\widetilde{\phi}^*(\Bbbk_{\xi^{N}}\otimes \Bbbk_\rho)^{\epsilon} = f^*L^{\epsilon N}\otimes \phi^*\cL^{\epsilon}.
\] 
It belongs to $\cQ_{g,n}^{\epsilon}(\cX,\cL,\beta)$ from   \Cref{cor_length_wrt_curlyX_is_the_same_as_length_wrt_specR}, and since $N\ge n_0\ge \lceil\frac{3}{\epsilon}\rceil$.
Thus to show that \[\cQ_{g,n}^{\epsilon}(\cX,\cL,\beta) = \cQ^{\epsilon}_{g,n}([\spec(R)/G\times\Gm],\Bbbk_{\xi^N}\otimes \Bbbk_\rho,\widetilde{\beta})\]
it suffices to show that there is $N=N(\epsilon,\cX,\cL)\gg 0 $ so that any stable quasimap in the moduli space $\cQ^{\epsilon}_{g,n}([\spec(R)/G\times\Gm],\Bbbk_{\xi^N}\otimes \Bbbk_\rho,\widetilde{\beta})$ lands in $\cX$.

\begin{Lemma}\label{lm:semi-invariants}
    For every sufficiently big integer $N_0$, there exists an integer $N$ such that every $(\xi^{N}\otimes \rho)^k$-semi-invariant of $R$ can be written as a sum of monomials of the form $f\cdot f'$ where $f'$ is a $\xi^{N_0k}$-semi-invariant.
\end{Lemma}
\begin{proof}
Consider the bigraded algebra \[A=\oplus_{k\ge 0, \ell\ge 0} R_{\xi^k\otimes \rho^\ell}\] generated, as an $R^{G\times \Gm}$-algebra, by the elements $f\in R$ such that $g\cdot f=\xi^{k}(g)\rho^{\ell}(g)f$ . Observe that $A$ is finitely generated over the ring of invariants, as we can regard it as the invariant subring $R[u,v]^{G\times\Gm}$ for the $G\times\Gm$-action that scales $u$ by $\xi$ and $v$ by $\rho$. 
Denote $f_1,\ldots, f_r$ the homogeneous generators such that $f_i$ has degree $a_i$ in $u$ and $b_i$ in $v$ with $b_i\neq 0$, and $f'_1,\ldots, f'_s$ have degree $a_i$ in $u$ and $0$ in $v$. Set $M=\lceil \mathrm{max}\{\frac{a_i}{b_i}\} \rceil$. 

Let $f$ be a $(\xi^N\otimes \rho)^k$-semi-invariant, i.e. $\deg(f)=(Nk,k)$ for some $k>0$. We can write $f$ as a sum of monomials of the form
\[ c_{I,J} \prod_{h=1}^{r} f_{h}^{i_h} \prod_{h=1}^{s} {f'_{h}}^{j_h}\text{ with }c_{I,J}\in R^{G\times \Gm}. \]
Observe that we must have $\sum_h a_h i_h + \sum a'_h j_h= Nk$ and $\sum b_h i_h=k$, thus
\begin{align*}
    \sum a_hi_h = \sum \frac{a_h}{b_h}\cdot b_hi_h \leq \sum M b_hi_h= Mk
\end{align*}
and in particular $\sum a'_h j_h \geq (N-M)k$. Picking $N=M+N_0$, we deduce that every $(\xi^{N}\otimes \rho)^k$-semi-invariant can be written as a sum of monomials of the form $f\cdot f'$ with $f'$ a $\xi^{N_0k}$-semi-invariant.
\end{proof}
Let $I$ be the ideal sheaf of the $\Bbbk_\xi$-unstable locus, which is the ideal generated by the subring of $\xi^d$-semi-invariants for $d>0$. Recall that the $\Bbbk_\xi$-stable locus coincide with $\cX$. Moreover, we are assuming that the ring of $\xi^d$-semi-invariants is generated in degree one. In particular, any $\xi^m$-semi-invariant can be written as a homogeneous polynomial of degree $m$ in the generators of $I$.

Thanks to \Cref{lm:semi-invariants} we know that for every $N_0$ there exists $N$ such that every $(\xi^{N}\otimes \rho)^k$-semi-invariant is a polynomial formed by monomials containing $\xi^{N_0k}$-semi-invariants. In particular, we can pick $N_0$ large enough so that every $\xi^{N_0}\otimes \rho$-semi-invariant
will be contained in $I^{\lceil \frac{1}{\epsilon}\rceil}$, as it will be a homogeneous polynomial in elements of $I$ appearing with high degree. Then the ideal sheaf of the unstable locus for $R$ with respect to the character $\xi^N\otimes\rho$ will be contained in $I^{\lceil \frac{1}{\epsilon}\rceil}$: any $\epsilon$-stable quasimap must land in $\cX$, otherwise there would be a basepoint with length greater than $\frac{1}{\epsilon}$.\end{proof}

\subsection{$0^+$-quasimaps}\label{sec:qmaps}
Observe that, from the main result of section \Cref{subsection_epsilon_stable}, it is not clear that the values of $\epsilon$ for which the spaces $\cQ_{g,n}^{\epsilon}(\cX,\cL,\beta)$ change, don't accumulate at 0, since the choice of $N$ in \Cref{thm:relative quasimaps are quasimaps} depends on $\epsilon$. The goal of this subsection is then to argue that 0 is not such an accumulation point.
The following definition is the analogue of $0^+$-quasimaps, to ous setting:
\begin{Def}\label{def_stable_qmap}
    Let $(\phi:\sC\to \cX, \Sigma_1,\ldots, \Sigma_n) $ be $n$-marked quasimap of genus $g$. It is \emph{stable} if $\omega_\sC( \Sigma_1+\ldots+\Sigma_n)\otimes \phi^*(p^*L^3\otimes \cL^{\epsilon})$ is ample for every $0<\epsilon$.
\end{Def}

The following lemma shows that one can rephrase the simpler stability condition above, purely in terms of the map $\phi$.
\begin{Lemma}\label{lemma_equivalent_conditions_for_stability}
    Let $(\phi\colon \sC\to \cX, \Sigma_1,\ldots, \Sigma_n) $ be a quasimap of genus $g\neq 1$. The following are equivalent:
    \begin{enumerate}
     \item $\phi$ is a twisted map (see \cite[Definition 2.2]{dli1}),
    \item the line bundle $\omega_\sC( \Sigma_1+\ldots+\Sigma_n)\otimes \phi^*(p^*L^3\otimes \cL^{\epsilon})$ is nef, and if $\sD\subseteq \sC$ is an irreducible component with coarse space $D$ such that $\omega_\sC( \Sigma_1+\ldots+\Sigma_n)\otimes \phi^*(p^*L^3\otimes \cL^{\epsilon})$ has degree zero on $\sD$, then $D\cong \bP^1$ and there are two geometric points $p_1,p_2\in \sD$ such that $\phi(p_1)$ and $\phi(p_2)$ are not isomorphic in $\cX$, and
        \item\label{randomitem} $\phi$ is stable.

        \end{enumerate} 
\end{Lemma}
\begin{proof}
First we prove that (1) and (2) are equivalent. Let $\sD\subseteq \sC$ be as in (2), with attaching nodes $n_1,n_2$ which, as usual, could be gerbes. Let's first assume that $\phi(\sD)\subseteq \cX_\dm$. From \cite{dli1}*{Proposition 2.6} the map $(\phi|_\sD\to \cX, n_1,n_2)$ is twisted if and only if it is a twisted \textit{stable} map in the sense of \cite{AV_compactifying}, if and only $\cD\to \cX_\dm$ is finite, if and only if there are two points of $\cD$ mapping to non-isomorphic points in $\cX_\dm$.
If instead $\phi(\sD)\not \subset \cX_\dm$, then not all points of $\sC$ map to isomorphic points on $\cX$: from the quasimap condition the generic locus of $\sD$ maps to $\cX_\dm$, but $\phi(\sD)\not \subset \cX_\dm$.


We now prove the equivalence between (2) and (3). Let $\sD$ be as in (2). It is immediate to check that $\phi^*p^*L^{3}$ is not ample when restricted to $\cD$, so $\sD\to \sX\to X$ is contracted to a point. Now, recall that $\cL$ descends to an ample line bundle on $\cX_\dm$, and let's first assume that $\phi(\sD)\subseteq \cX_\dm$. Then $\cD\to \cX_\dm$ is finite if and only if there are two points $p_1$ and $p_2$ as above, if and only if $\phi^*\cL$ is ample. 

If instead $\phi(\sD)\not\subseteq \cX_\dm$ we need to show that $\phi^*\cL$ is ample. But since $\phi(\sD)\not\subseteq \cX_\dm$, there is an open $U\subseteq X$ containing $\pi\circ\phi(\sD)$ and a section of $\cL|_{\pi^*\cU}$ which does not vanish identically when pulled back to $\sD$, but has a zero. In particular, $\phi^*\cL$ has a non-zero and vanishing global section, so it has positive degree.
\end{proof}

\begin{Def}
We define the category $\cQ_{g,n}(\cX,\cX_\dm,\beta)$ of stable quasimaps to $\cX$ with respect to $\cX_\dm$ of class $\beta$.
\end{Def}
\begin{Teo}\label{thm:quasimaps are dm}
    The fibered category $\cQ_{g,n}(\cX,\cX_\dm,\beta)$ is a proper Deligne-Mumford stack and coincide with $\cQ_{g,n}^{\epsilon}(\cX,\cX_\dm,\beta)$ for $\epsilon >0$ small enough.
\end{Teo}
\begin{proof}
It follows from the definition of stable and $\epsilon$-stable quasimap, together with \Cref{rmk_pull_back_of_L_is_ample}, that it suffices to prove that there is an upperbound for the lenght of a basepoint for a quasimap in $\cQ_{g,n}(\cX,\cX_\dm,\beta)$. Consider an embedding
\[
\iota\colon \cX\hookrightarrow[\spec(R)/G\times\Gm]
\]
and two characters $\xi$ and $\rho$ as in \Cref{prop_ideal_sheaf_unstable_locus_pulls_back}. Consider also a class $\widetilde{\beta}$ extending $\beta$. Then a quasimap $\phi\colon \sC\to \cX$ of class $\beta$ will induce a quasimap $\iota\circ\phi$ of class $\widetilde{\beta}$. Those are bounded from \cite[\S 2.4.3 (i)]{orbifold_qmap_theory}. Consider then a bounding family $\cC\to B$ with a map $\phi\colon \cC\to [\spec(R)/G\times\Gm]$. Consider $\phi^*I$, the pull-back of the ideal sheaf $I$ of the unstable locus for $[\spec(R)/G\times\Gm]$ with respect to $\Bbbk_{\xi^{N}\otimes \rho}$, and consider
\[\phi^*I\to \cO_{\cC}\to \operatorname{Coker}\to 0.
\]
For every $b\in B$, from the right-exactness of tensor products, such a map will commute with base change, so \[\operatorname{length}(\operatorname{Coker}|_b) = \sum _{x\text{ base point for }\phi_b\colon \sC_b\to [\spec(R)/G\times\Gm]}\ell(x).\] In particular, the closed substack of $\cC$ corresponding to $\operatorname{Coker}$ is finite over $B$, and since $B$ is of finite type, the fibers have bounded length. 
\end{proof}

\section{Applications of \Cref{teo_intro_cX_contains_open_proper_dm}}\label{section_examples_qmaps}
In this section we report two examples when one can use \Cref{teo_intro_cX_contains_open_proper_dm} to compactify the space of maps from curves to certain algebraic stacks. We first show how to use \Cref{teo_intro_cX_contains_open_proper_dm} to construct a compact moduli of
varieties fibered in Calabi--Yau pairs. Then we show that if $\cX$ is a quotient of a separated Deligne--Mumford stack by a torus, and the set of properly stable points is dense in $\cX$, then $\cX$ always contains an open substack which is a proper Deligne--Mumford stack.

\subsection{Moduli of fibered Calabi--Yau pairs}
We begin by recalling that:
\begin{enumerate}
    \item there is an algebraic stack $\cX:=\cD\cP_m^{\CY}$ with a projective good moduli space parametrizing pairs $(S,cD)$ where $\dim(S)\le 2$, $D\subseteq S$ is a divisor such that the singularities of $(S,cD)$ are semi-log-canonical, $K_S+cD\sim_\mathbb{Q}0$, the sheaf $(\omega_S^{\otimes m})^{*}$ is an ample Cartier divisor, and $S$ admits a $\bQ$-Gorenstein smoothing; see \cite{blum2024good}*{Theorems 1.1, 1.2 and 3.6}.
    \item When $c<1$ and $m\gg 0$, there is an open, dense and proper Deligne--Mumford stack $\cX_\dm:=\cD\cP_m^{\operatorname{KSBA}}\subseteq \cX$ parametrizing pairs $(S,cD)$ such that $(S,(c+\epsilon)D)$ is still semi-log-canonical for some $0<\epsilon \ll 1$ \cite{blum2024good}*{Theorems 1.1 and 1.2}.
    \item  For $c<1$ and $m\gg1$, if $p\colon(\sS,c\sD)\to \cX$ is the universal family, then from cohomology and base change $p_*\cO_\cS(d K_{\sS/\cX} + d(c+\epsilon)\sD)$ is a vector bundle, and $\det\left(p_*\cO_\sS(d K_{\cS/\cX} + d(c+\epsilon)\sD)\right)$ restricted to $\cX_\dm$ descends to an ample line bundle on $X_\dm$ \cite{KP17}.
    \item $\cX$ is a global quotient stack \cite{ascher2023moduli}*{proof of Theorem 3.7}.
    \item For $m\gg 1$, the algebraic stack $\cX$ admits a projective good moduli space $q:\cX\to X$ \cite{blum2024good}*{\S 6.4}.
    \item The map $q$ identifies two pairs $(S_1,cD_1)$ and $(S_2,cD_2)$ if and only if they are $S$-equivalent \cite{blum2024good}. 
\end{enumerate}
{We aim at proving \Cref{teo_intro_kodairadim1}, whose content we recall below for the convenience of the reader.
\begin{Teo}
    Set $\cX:=\cD\cP^{\CY}_m$ and $\cX_\dm=\cD\cP^{\operatorname{KSBA}}_m$, and assume that $\cX_\dm$ has complement of codimension at least 2. Then the assumptions of \Cref{teo_intro_cX_contains_open_proper_dm} apply for the inclusion $\cX_\dm\subseteq \cX$. In particular:
    \begin{enumerate}
        \item the stack $\cQ_g(\cX,\cX_\dm,\beta)$ compactifies the space of maps $\pi\colon(Y,cD)\to C$ with fibers in $\cX$ such that:
        \begin{itemize}
            \item[(Sing)] the curve $C$ is smooth and the generic fiber of $\pi$ has klt singularities,
            \item[(Stab)] either $\omega_C$ is ample, or not all the fibers of $\pi$ are $S$-equivalent,
            \item[(Num)] the family $\pi$ comes from a map $C\to \cX$ of class $\beta$ from a curve of genus $g$.
        \end{itemize}
        \item the boundary of $\cQ_g(\cX,\cX_\dm,\beta)$ parametrizes families $\pi:(\sY,c\sD)\to \sC$ of pairs in $\cX$ with fibers in $\cX$, fibered over a twisted curve $\sC$, such that:
        \begin{itemize}
            \item[(Sing)] the set $\Delta:=\{p\in \sC:(\sY_p,(c+\epsilon)\sD_p)$ does \underline{not} have semi-log-canonical singularities for any $0<\epsilon \ll 1\}$ is a finite union of smooth points $\sC$,
            \item[(Stab)] if $\sR\subseteq \sC$ is an irreducible component such that $\deg(\omega_\sC |_\sR)< 0$, then not all the fibers of $\pi|_\sR\colon(\sY|_\sR,c\sD|_\sR)\to \sR$ are $S$-equivalent; whereas if $\deg(\omega_\sC |_\sR)=0$, then not all fibers of $\pi|_\sR$ are isomorphic, and
            \item[(Num)] the family $\pi$ comes from a map $\sC\to \cX$ of class $\beta$ from a twisted curve of genus $g$.
        \end{itemize}
    \end{enumerate}
\end{Teo}}
\begin{Remark}We are not aware of an example where $\cX_\dm$ does not have complement of codimension at least 2.    
\end{Remark}
\begin{Remark}\label{remark_stefano_some_fibers_are_not_slc}
    As stated, the fibrations $\pi\colon (Y,cD)\to C$ in part (1) of \Cref{teo_intro_kodairadim1}
    must have \textit{all} fibers that are in $\cD\cP^{\CY}_m$ and the generic fiber in $\cD\cP^{\operatorname{KSBA}}_m$.
    However, if one has a morphism $\pi':(Y,D)\to C$ over a smooth curve $C$ such that $D$ does not have any vertical component and the generic fiber of $\pi$ is parametrized by $\cD\cP^{\operatorname{KSBA}}_m$, one can proceed as follows. Let $U\subseteq C$ be the open locus where the fibers of $\pi$ are in $\cD\cP^{\CY}_m$.
    One can perform a sequence of root stacks $\cC\to C$ so that the generic morphism $ U\to \cD\cP^{\CY}_m$
    extends to a morphism $\cC\to \cD\cP^{\CY}_m$ (by applying \cite{bresciani2022arithmetic} to the generic morphism to $\cD\cP^{\operatorname{KSBA}}_m$). If we denote by $\Sigma_1,\ldots,\Sigma_n$ be the stacky points on $\cC$, the map $\phi\colon(\cC,\Sigma_1,\ldots\Sigma_n)\to \cD\cP_m^{\CY}$ is \textit{pointed} stable, and so one can compactify the space of such maps using \Cref{teo_intro_cX_contains_open_proper_dm}.

    The map $\phi$ induces a family $(\cX,c\Delta_{\cX})\to \cC$ whose coarse moduli space $(X,c\Delta)\to C$ is obtained from $(Y,cD)$ as follows. For each point $x\in C$ such that $(Y_x,cD_x)$ is not slc, consider $U$ a neighbourhood of $x$ such that for any $z\in U\smallsetminus\{x\}$ the fiber $(Y_z,(c+\epsilon)D_z)$ is slc.
    Consider $W_U\to Y_U$ a log-resolution of the pair $(Y_U,(1+\epsilon)D_U + F)$ where $F$ is the fiber over $x$ with the reduced structure. If we denote by $cD_W$ the proper transform of $cD_U$ and with $E$ the exceptional divisors of $W_U\to Y_U$ with the reduced structure, the pair $(X_U,c\Delta_U)$ is a weak minimal model for $(W_U,cD_W + E)$ over $U$ and a canonical model for $(W,(c+\epsilon)D_W + E)$ over $U$.  
\end{Remark}
\begin{proof}[Proof of \Cref{teo_intro_kodairadim1}]{To prove that $\cX_\dm$ is the semistable locus relative to $X$ for a line bundle on $\cX$ we can use \Cref{cor_when_X_dm_is_ss_locus_for_a_lb} with the line bundle $\det\left(p_*\cO_\sS(d K_{\cS/\cX} + d(c+\epsilon)\sD)\right)$. The only other statement one needs to check is that if $f\colon\sC\to \sX$ is stable, then it satisfies condition (2) part (Stab) in \Cref{teo_intro_kodairadim1}. This follows from \Cref{lemma_equivalent_conditions_for_stability}: since the generic points of $\sC$ map to $\sX_\dm$, if $f$ is stable then any irreducible component $\sR\subseteq \sC$ such that $\deg((\omega_\sC)|_\sR)<0$ must map quasi-finitely to $X$, and so from point (6) above there must be two fibers of $\pi|_\cR$ not $S$-equivalent. Assume $\sR$ is such that $\deg((\omega_\sC)|_\sR)=0$. Since $f$ is a quasimap, the generic point of $\sR$ maps to a point $x\in\cX_\dm$, and $f$ is stable unless all the points of $\sR$ map to $x$. }    
\end{proof}
\subsection{Toric quotients}\label{section_toric_quotients}
We prove that if $\cX$ is a quotient of a Deligne--Mumford stack by a split torus, and there is an open dense $U\subseteq X$ such that $\cX\times_XU$ is a Deligne--Mumford stack, then $\cX$ always contains an open substack which is a proper Deligne--Mumford stack. So for such stacks, the enlargement $\cX\subseteq \widetilde{\cX}$ of \Cref{thm_intro_simplified} is not needed{: if $\cX$ is a global quotient, \Cref{teo_intro_cX_contains_open_proper_dm} applies}.
\begin{Teo}\label{teo_toric_stack_has_open_dense_dm}
    Let $\cX$ be an algebraic stack with a good moduli space $\cX\to X$ which is separated. Assume that:
    \begin{itemize}
        \item there is a morphism $\cX\to \cB \Gm^n$,
    representable in Deligne--Mumford stacks, and
    \item there is a dense open $U\subseteq  X$ such that $\cX\times_XU$ is Deligne--Mumford.
    \end{itemize} Then there are finitely many hyperplanes $H_1,...,H_m$ in the space $\mathbf{X}(\Gm^n)_\bQ$ of rational characters of $\Gm^n$ such that, if $\mu \notin \bigcup H_i$, then $\cX(\Bbbk_\mu)^{ss}_X$ is a separated Deligne--Mumford stack, with coarse moduli space projective over $X$ and such that $\cX\times_XU\subseteq \cX(\Bbbk_\mu)^{ss}_X$.
\end{Teo}

\begin{proof}
We will denote by $G$ a group which is a central extension of $\Gm^r$ by a finite group $F$. 

\textbf{Affine case.} We will assume that $\cX=[\spec(A)/G]$, and that $\spec(A)$ has a $G$-fixed point. In particular, our assumptions guarantee that there is a morphism $\psi:G\to \Gm^n$ with finite kernel.  

    \textbf{Step 1.} We can assume $G=\Gm^r\times F$.
    
    From \cite{brion2015extensions}, there is a finite subgroup $F<G$ and a surjective morphism $\Gm^r \rtimes F \to G$ with finite kernel. As $\Gm^r$ is contained in the center of $G$, the product is direct. In particular, there is a morphism $[\spec(A)/\Gm^r\times F]\to [\spec(A)/G]$ which is separated and a gerbe, so if $\mu$ is a character of $\Gm^n$, then $[\spec(A)/\Gm^r\times F](\Bbbk_\mu)^{ss}$ satisfies the condition of the theorem if and only if $[\spec(A)/G](\Bbbk_\mu)^{ss}$ does. In particular, we can assume that $G=\Gm^r\times F$.

    \textbf{Step 2.} We can assume that $F=\{1\}$.

    Indeed, up to post-composing $\cX\to \cB\Gm^n$ with the map $\cB \Gm^n\to \cB\Gm^n$ induced by $(t_1,\ldots ,t_n)\mapsto (t_1^N,\ldots ,t_n^N)$ for $N\gg 1$, we can assume that the morphism $\psi:\Gm^r\times F\to \Gm^n$ is trivial if restricted to $F$.  In particular, we can factor $[\spec(A)/\Gm^r\times F]\to \cB\Gm^n\times \spec(A^{\Gm\times F})$ as
    \[[\spec(A)/\Gm^r\times F]\to [\spec(A^F)/\Gm^r]\to \cB\Gm^n\times \spec(A^{\Gm\times F})\]
    and again, if we can find such a $\mu$ for $[\spec(A^F)/\Gm^r]$, the same $\mu$ will also work for $[\spec(A)/\Gm^r\times F]$. 

    \textbf{Step 3.} We can assume that $\spec(A)=\bA^N$.

    Indeed, there is an equivariant closed embedding $i:\spec(A)\hookrightarrow \bA^N$. From the affine Hilbert-Mumford criterion, since $i$ is a closed embedding, we have that $i^{-1}(\bA^N(\Bbbk_\mu)^{ss}) =\spec(A)(\Bbbk_\mu)^{ss}$. Since $\Gm^r$ is reductive, there is a closed embedding $\spec(A^{\Gm^r})\hookrightarrow \bA^N/\!\!/\Gm^r$ where we denoted by $\bA^N/\!\!/\Gm^r$ the good moduli space of $[\bA^N/\Gm^r]$. So if we find such a $\mu$ for $[\bA^N/\Gm^r]$, the same $\mu$ will also work for $[\spec(A)/\Gm^r]$.

    \textbf{The $\bA^N$ case.} Up to performing a change of coordinates, we can assume that the action is diagonal, given by $N$ characters $\chi_1,\ldots,\chi_N$.
     
    Given a character $\widetilde{\mu}$, recall that:
    \begin{itemize}
        \item from \cite{Hos}*{Proposition 2.5} a point $p\in \bA^N$ is $\widetilde{\mu}$-stable
    if, for every one-parameter subgroup $\lambda:\Gm\to \Gm^r$ such that $\lim_{t\to 0}\lambda(t)p$ exists, we have that $<\lambda, \widetilde{\mu}>\le 0$, and
    \item given $\lambda\in \mathbf{X}(\Gm^r)^*$, if a point $p$ is such that $\lim_{t\to 0}\lambda(t)p$ exists and if $<\lambda,\chi_j><0$, then the $j$-th coordinate of $p$ is 0.
    \end{itemize} 
    
    Since the generic stabilizer is finite, the morphism $\Xi:\Gm^r\to \Gm^N$ that sends $g\mapsto (\chi_1(g),\ldots, \chi_N(g))$ has finite kernel. In particular, the point $(1,\ldots,1)$ is a stable point, as there is no one-parameter subgroup $\lambda$ such that $\lim_{t\to 0}\lambda(t)p$ exists. Moreover,
    if we choose $\widetilde{\mu}$ such that the hyperplane in
    $\mathbf{X}(\Gm^r)^*$
    given by $<\cdot, \widetilde{\mu}>=0$ is different from $<\cdot, \chi_i>=0$ for every $i$, then the stable locus and the semistable locus agree. Indeed, if $p$ is a point such that,
    for a given $\lambda$, we have that $<\lambda,\widetilde{\mu}>=0$ and $\lim_{t\to 0}\lambda(t)p$ exists, then $\lambda$ is away from the hyperplanes $<\cdot, \chi_i> = 0$. In particular, we can slightly perturb $\lambda$ in a way such that $\lim_{t\to 0}\lambda(t)p$ still exists, but $<\lambda,\widetilde{\mu}> \gg 0$, so $p$ was not semistable.
    
    Finally, since $\psi:G\to\Gm^n$ has finite kernel, the morphism
    $\Psi:\mathbf{X}(\Gm^n)_{\mathbb{Q}}\to \mathbf{X}(\Gm^r)_{\mathbb{Q}}$ induced by $\psi:G\to \Gm^n$ is surjective.
    So as long as we take $\mu$ away from $\bigcup\Psi^{-1}(\chi_i)$,
    the character $\widetilde{\mu}=\Psi(\mu)$ will work. 

    \textbf{General case.} Since $\cX\to \cB\Gm^n$ is representable in Deligne--Mumford stacks, from \cite{di2024effective}*{Corollary 4.7 \& Theorem 2.19} the stabilizers $G$ of the geometric points are central extensions of a split torus $\Gm^r$ by a finite group $F$.
    In particular, for every $p\in X$ there is an \'etale neighbourhood of $p$ of the form $\spec(A^G)\to X$ such that the following diagram is cartesian:
    \[
    \xymatrix{[\spec(A)/G]\ar[r] \ar[d] & \cX\ar[d] \\ \spec(A^G)\ar[r] & X}
    \]
    and there is a homomorphism $\psi:G\to \Gm^n$ with finite kernel. So for each neighbourhood as above we can find a finite set of hyperplanes of $\mathbf{X}(\Gm^n)_{\mathbb{Q}}$ as in the statement of the theorem.
    Now the desired result follows since $X$ is quasicompact, so we can cover $X$ with finitely many \'etale charts as above. 
\end{proof}
\begin{example}
We give an application of \Cref{teo_toric_stack_has_open_dense_dm} to study stable maps to a specific quotient. Consider the action of $\Gm$ on $\bA^2$ with weights 1 and $-1$. The good moduli space $[\bA^2/\Gm]$ is $\bA^1$, and the quotient map $[\bA^2/\Gm]\to \bA^1$ is an isomorphism away from a point (which we set to be 0). We can consider the stack $\sP^1$ obtained by replacing the chart in $\bP^1$ given by $x_0\neq 0$ with $[\bA^2/\Gm]$.
More specifically, this is the gluing of $\bA^1_t$ and $[\bA^2_{x,y}/\Gm]$ along $\bA^1\smallsetminus\{0\}$ and $[\bA^2_{x,y}\smallsetminus\{xy=0\}/\Gm]$ via the map that sends $t\mapsto \frac{1}{xy}$. The stack $\sP^1$ is an algebraic stack with $\bP^1$ with a good moduli space. It is smooth and it has a Zariski open cover where it is a quotient of a separated scheme by a torus.
Then from the main result of \cite{di2024effective}, it is globally a quotient of a separated Deligne--Mumford stack by a torus. Indeed in \textit{loc. cit.} we give a criterion for when a smooth algebraic stack with a good moduli space admits such a presentation, and the criterion is Zariski local on the good moduli space. Therefore \Cref{teo_toric_stack_has_open_dense_dm} applies, and one can check explicitly that $\sP^1$ contains $\bP^1$ as an open substack, since the quotient map $[\bA^2/\Gm]\to \bA^1$ admits a section given by $[\bA^2\smallsetminus\{x=0\}/\Gm]\to \bA^1$. This extends to a section $\bP^1\to \sP^1$. 

Consider the map $\bA^1_t\to \bA^2_{x,y}$, $x\mapsto t-a$ and $y\mapsto t-b$ for $a,b\in k$. This gives a map $\bA^1\to [\bA^2/\Gm]$ which, if composed to $[\bA^2/\Gm]\to \bA^1\hookrightarrow\bP^1$, compactifies to a map $\phi:\bP^1\to\sP^1$ whose composition with the good quotient of $\sP^1$ gives $\bP^1\to \sP^1\to \bP^1$ of degree 2, and not ramified at 0 (the point on $\bP^1$ corresponding to a polystable and not stable point of $\sP^1$).  If we take as open proper Deligne--Mumford of $\sP^1$ the locus in $[\bA^2/\Gm]$ where $x\neq0$, and we denote by $C$ the domain of $\phi$, then we have a unique point in $C$ which does not map to $\sP^1_\dm$. If one takes the closure of the locus of maps $\bP^1\to \sP^1$ constructed as above, it is straightforward to check that the nodal curves that one finds on the boundary are either $D$ the nodal union of two $\bP^1$s, with the composition $D\to \sP^1\to \bP^1$ being unramified over $0$, or the nodal union $D$ of three $\bP^1$s, with the two $\bP^1$s which have a single node mapping finitely to the good quotient of $\sP^1$, and the $\bP^1$ with two nodes mapping to $0$. The map $\phi:D\to \sP^1$ will still be a quasimap, so $\phi^{-1}(\sP^1\smallsetminus\sP^1_\dm)$ will still be a unique smooth point.
\end{example}

\section{Extended weighted blow-ups}\label{section_extended_blowups}
In this section we first introduce extended weighted blow-ups, which are a mild generalization of weighted blow-ups, and are more amenable for compactifying moduli spaces of maps to algebraic stacks. Then we prove \Cref{thm_you_can_embed_a_Stack_in_one_with_open_dm_by_taking_deformations_to_nc}, which allows to embed any algebraic stack $\cX$ with a good moduli space $\cX\to X$ and a properly stable point in an algebraic stack $\widetilde{\cX}$ with the same good moduli space $\widetilde{\cX}\to X$, and such that $\widetilde{\cX}$  satisfies \Cref{assumptions_epsilon_qmaps} if $\cX=[W/G]$ for $G$ reductive and $X$ is projective.
In particular, we can combine \Cref{thm_you_can_embed_a_Stack_in_one_with_open_dm_by_taking_deformations_to_nc} with the results of \Cref{sec:qmaps} to get the following
\begin{Cor}
    Let $\cX$ be an algebraic stack with a good moduli space $p:\cX\to X$, and with an open dense $U\subseteq X$ such that $\cX\times_XU$ is Deligne--Mumford. Assume that $\cX=[W/G]$, where $G$ is a reductive group. Then there is a schematically dense open embedding $\cX\subseteq \widetilde{\cX}$ as in \Cref{thm_you_can_embed_a_Stack_in_one_with_open_dm_by_taking_deformations_to_nc}, and there is an algebraic stack $\cQ_{g,n}(\widetilde{\cX}, \widetilde{\cX}(\cL_\dm)^{ss}_X,\beta)$ parametrizing stable quasimaps $\sC\to \widetilde{\cX}$ of genus $g$ and class $\beta$. 
    \end{Cor}
\begin{Remark}
    Recall locus $\cS:=\{x\in X:p^{-1}(x)$ is a Deligne--Mumford gerbe$\}$ is open in $X$ from \cite{ER21}*{Proposition 2.6}. So if $X$ is irreducible, if there is a point in $\cS$ then $\cS$ is a dense open subscheme of $X$.
\end{Remark}
\subsection{Extended weighted blow-ups}
In this subsection we define \textit{extended weighted blow-ups}.
On a first approximation, an extended blow-up of an algebraic stack $\cX$ along an ideal $\cI$ is the $\Gm$-quotient of the deformation of $\cX$ to the weighted normal cone determined by $\cI$.
More specifically, if we denote by $\cZ\subset\cX$ the closed substack defined by $\cI$, the deformation to the weighted normal cone of $\cX$ is (a subset of) the weighted blow-up of $\cX\times\bA^1$ along $\cZ\times\{0\}$, where the parameter $T$ of $\bA^1$ is defined of degree $-1$. The induced $\Gm$-action on $\cX\times\bA^1$ extends then to
the deformation of the normal cone, and the extended weighted blow-up is the $\Gm$-quotient of the latter.

More formally, we begin with the following

\begin{Def}[\cite{QR}*{Definition 3.1.1}]\label{def_weighted_embedding}
   Let $\cX$ be an Artin stack. A weighted embedding $\cY_{\bullet} \hookrightarrow \cX$ is defined by a sequence of closed embeddings $\{ \cY_n=V(\cI_n) \hookrightarrow \cX \}_{n \ge 0}$ such that: \begin{itemize}
        \item $\cI_0 \supset \cI_1 \supset \dots \supset \cI_n \supset \dots$
        \item $\cI_n \cI_m \subset \cI_{n+m}$
        \item Locally in the smooth topology on $X$, there exists a sufficiently large positive integer $d$ such that for all integers $n \ge 1$, \[\cI_n=\left(\cI_1^{l_1}\cI_2^{l_2} \cdots \cI_d^{l_d} \; : \; l_i \in \bN, \ \sum_{i=1}^d il_i=n\right)\] in which case, we say $\cI_\bullet$ is generated in degrees $\le d$. 
    \end{itemize} 
    The last condition can be understood as a needed condition for $\bigoplus_n \cI_n$ to be of finite type.
    Furthermore, we set $\cI_n=\cO_{\cX}$ for $n\leq 0$ and we call the sequence of ideals $\{\cI_n\}_{n\in \bZ}$ a \textit{weighted ideal sequence}. 
\end{Def} 

Consider then the graded $\cO_\cX$-algebra $\cA:=\oplus_{n\in\bZ} \cI_n$. We will use an auxiliary variable $T$ to denote the degree, so in particular we will write an element of $\cA$ as $\sum_{j=-k}^ka_jT^j$, where $a_j$ is an element of $\cI_j$ in degree $j$. Observe that:
\begin{enumerate}
    \item there is a morphism $\pi^\#:\cO_\cX\to \cA$, which sends $\cO_\cX$ to the degree zero component;
    \item the previous morphism has a left inverse $i^\#:\cA \to \cO_\cX$ defined by $\sum_{j=-k}^ka_jT^j\mapsto \sum_{j=-k}^ka_j$.
\end{enumerate}
Moreover, we have the following.
\begin{Lemma}\label{lm:action}
    There is a strict $\Gm$-action on $\spec_{\cO_\cX}(\cA)$ over $\cX$ in the sense of \cite{Rom05}*{Definition 1.3}. In particular, there is a quotient $[\spec_{\cO_\cX}(\cA)/\Gm] \to \cX$ whose formation commutes with base change.
\end{Lemma}
\begin{proof}
    The grading of $\cA$ induces a co-action
    $\cA \longrightarrow \cA[u^{\pm 1}],\text{ } aT^j\longmapsto aT^ju^j.$
    We can then define an action $\mu\colon\Gm \times \spec_{\cO_\cX}(\cA) \to \spec_{\cO_\cX}(\cA)$ as follows: \begin{center} the object $(\lambda \in \cO_S^*(S), \xi\colon S\to\cX, f\colon \xi^*\cA \to \cO_S)$ is mapped to $(\xi, \xi^*\cA \to \xi^*\cA[u^{\pm 1}] \overset{{u=\lambda}}{\longrightarrow} \xi^*\cA \overset{f}{\to} \cO_S).$\end{center}
    Moreover, as the automorphism group of $(\lambda,\xi,f)$ coincides with the one of $(\xi,f)$ and $\xi$, the morphism between automorphism groups is just the identity; in this way, we get a morphism of algebraic stacks.

    To check that the morphism above defines a strict action, we need to check that the identity (regarded as a $2$-morphism) makes the two following diagrams of algebraic stacks $2$-commute:
    \[
    \begin{tikzcd}
        \Gm ^{\times 2}\times \spec_{\cO_\cX}(\cA) \ar[r, "m\times \id"] \ar[d, "\id\times\mu"] & \Gm \times \spec_{\cO_\cX}(\cA) \ar[d, "\mu"] \\
        \Gm \times \spec_{\cO_\cX}(\cA) \ar[r, "\mu"] & \spec_{\cO_\cX}(\cA)
    \end{tikzcd} 
    \begin{tikzcd}
        \Gm \times \spec_{\cO_\cX}(\cA) \ar[r, "\mu"] & \spec_{\cO_\cX}(\cA) \\
        \spec_{\cO_\cX}(\cA). \ar[u, "1\times\id"] \ar[ur, "\id"]
    \end{tikzcd}
    \]
    This is pretty straightforward, and it follows from the fact that the action morphism is a morphism over $\cX$, where the $2$-commutativity is strict (meaning that the $2$-morphism making the diagram $2$-commute is actually the identity); we omit the details. Then from \cite{Rom05}*{Theorem 4.1} it follows that there exists an algebraic stack $[\spec_{\cO_\cX}(\cA)/\Gm]\to \cX$ whose formation commute with base change.
\end{proof}
Given \Cref{lm:action}, we observe also the following:
\begin{enumerate}
    \setcounter{enumi}{2}
    \item if we denote by $\pi$ the morphism induced by $\pi^\#$ on sections, then $\pi$ is $\Gm$-equivariant and the induced morphism $[\spec_{\cO_\cX}(\cA)/\Gm]\to X$ is a good moduli space as the $\Gm$-invariant functions are those on degree zero;
    \item if we denote by $i$ the morphism $\cX\to [\spec_{\cO_\cX}(\cA)/\Gm]$ induced by $i^\#$, this is an open embedding \cite{QR}*{Proposition 4.3.4};
    \item there is a morphism $ [\spec_{\cO_\cX}(\cA)/\Gm]\to \cX\times[\bA^1/\Gm]$, given by the inclusion $\bigoplus_{n\le 0}\cO_\cX[T^{-1}]\to \cA$, which induces an isomorphism of $i(X)$ with the preimage of $X\times \{1\}$ and of its complement with the complement of $\{T^{-1}\}=0$; in particular, we have that $[\spec_{\cO_\cX}(\cA)/\Gm]\smallsetminus i(X)$ is Cartier.
    \item from \cite{QR}*{Remark 3.2.5}, we have that $[\spec_{\cO_\cX}(\cA)\smallsetminus V(\cA^+)/\Gm]$ is the weighted blow-up of $\cX$ along $\cI_\bullet$.
\end{enumerate}
\begin{Def}\label{def_ext_w_blowup} Given a weighted ideal sequence $\{\cI_n\}$ as above, we define \[
\EB_{\cI_\bullet}\cX:=[\spec_{\cO_\cX}(\cA)/\Gm]
\]
     the \textit{extended weighted blow-up} of $\cX$ along $\cI_\bullet$. When $\cI_n = \cI_1^n$ for $n\ge 1$, we call $\EB_{\cI_\bullet}\cX$ simply an \textit{extended blow-up} and we denote it by $\EB_{\cI}\cX$.
\end{Def}
\begin{Remark}
    An extended weighted blow-up $\pi:\EB_{\cI^\bullet}\cX\to \cX$ is such that $\pi_*\cO_{\EB_{\cI^\bullet}\cX}=\cO_{\cX}$. Indeed, one can check this smooth locally over $\cX$ so we can assume that $\cX$ is a scheme, in which case it follows from point (3) above.
\end{Remark}

\begin{example}[Extended blow-up of $0\in \bA^2$]\label{example_ext_blowup_of_origin_in_A2} In this example we work out explicitly the extended blow-up of $0\in \bA^2$.
    The deformation to the normal cone of the origin $V(x,y)\subseteq \bA^2$ can be identified with $\spec(k[x,y,T,X,Y]/(XT-x, YT-y)$ where:
    \begin{itemize}
        \item $X,Y$ are the generators of the maximal ideal in degree 1,
        \item $T$ has degree $-1$,
        \item $x,y$ have degree 0,
        \item the map $X\to \spec(I^\bullet)$ is given by sending $T\mapsto 1$, and 
        \item the map $\spec(I^\bullet)\to X$ is given by the inclusion of the degree 0 component.
    \end{itemize}
    Therefore we can identify it with $\bA^3=\spec(k[X,Y,T])$, which has a $\Gm$-action with weights $(1,1,-1)$.
    The map $\bA^3\to \bA^2$ given by $(a,b,u)\mapsto (au,bu)$
    induces the good moduli space morphism \[\pi:[\bA^3/\Gm]\to \bA^2\] which is the extended blow-up of $0\in \bA^2$. The section $i$ of point (4) is given by the locus where $T=1$ in $[\bA^3/\Gm]$. As for point (6), we can identify $[\{(a,b,t):(a,b)\neq (0,0)\}/\Gm]\hookrightarrow[\bA^3/\Gm]$ with the blow-up of the origin in $\bA^2$ as in \cite{Inc}*{\S 2.2}.
    
    Similarly, when a group $G$ acts on $\bA^2$, we can extend its action to $\bA^3$,
    by acting trivially on the $T$ component. So as before one has two morphism $[\bA^2/G]\to[\bA^3/G\times \Gm] $ and $[\bA^3/G\times \Gm]\to [\bA^2/G]$ whose composition is the identity, where the first is the complement of a Cartier divisor, and $[\bA^3/G\times \Gm]$ contains as an open the blow-up of $\cB G$ in $[\bA^2/G]$. 
\end{example}
\begin{Def}\label{def_gen_blowup_comes_with_a_cartier_divisor} Given a stack $\cX$ and a weighted ideal sequence $\cI^\bullet$, there is an effective Cartier divisor $\cE\subseteq \EB_{\cI^\bullet}X$ which restrict to the exceptional divisor of the blow-up, and which in the examples given above it agrees with the vanishing locus of $T^{-1}$. We call such an effective Cartier divisor the \textit{exceptional divisor}. Observe that the exceptional divisor of an extended weighted blow-up $\EB_{\cI^\bullet}X$ is the scheme-theoretic closure of the (usual) exceptional divisor of the weighted blow-up contained (as an open) in $\EB_{\cI^\bullet}X$.
\end{Def}
\begin{Remark}
    One can check that if $\cZ\subseteq \cX$ is a smooth closed substack of $\cX$ with ideal $\cI$, and $\cI^\bullet = \{\cI^n\}$ then the extended blow-up $\EB_{\cI^\bullet}\cX$ is smooth.
\end{Remark}
\subsection{Proof of \Cref{thm_you_can_embed_a_Stack_in_one_with_open_dm_by_taking_deformations_to_nc}}
We begin with the following.
\begin{Lemma}\label{lemma_ext_blowup_is_relative_ss_locus}
    Let $\cX$ be an algebraic stack and $\cI_\bullet:=\{\cI_n\}_{n\in \bZ}$ a weighted ideal sequence of $\cO_\cX$. Let $\theta:\Gm\to \Gm$ be the character $t\mapsto t^{-1}$.
    Then $(\EB_{\cI_\bullet} \cX)(\Bbbk_{\theta} )^{ss}_{\cX}$
    is the weighted blow-up of $\cX$ along $\cI_\bullet$. Moreover, if $\pi:\cX\to X$ is a good moduli space and $\cI_n=\cI_1^n$, then $(\EB_{\cI_\bullet} \cX)(\Bbbk_{\theta} )^{ss}_{X}\cong \Bl^\pi_\cI\cX$, where the latter is a saturated blow-up (see \cite{ER21}*{Definitions 3.1 and 3.2}).
\end{Lemma}
\begin{proof}If we denote by $p:\EB_{\cI_\bullet}\cX\to \cX$ the projection, observe that $
    \pi_*\Bbbk_{\theta}^{\otimes n} = \cI^n.$
    This can be checked smooth locally over $\cX$, where it suffices to observe that the sections of $\spec_{\cO_\cX}(\bigoplus_{n\in \bZ}\cI^n)$ which are $\Bbbk_{\theta}^{\otimes n}$-semiinvariant are homogeneous of degree $n$.
    The first part now follows since $(\EB_{\cI_\bullet}\cX)(\Bbbk_{\theta}  )^{ss}_{\cX} = [\spec(\bigoplus_{n\in \bZ} \cI_n)\smallsetminus V(<\bigoplus_{n>0}\cI_n>)/\Gm]$ and from \cite{QR}*{Remark 3.2.4 and \S 1.1} the latter is the desired weighted blow-up.

     For the moreover part, observe that $(\EB_{\cI} \cX)(\Bbbk_{\theta} )^{ss}_{X}$ is the locus of $\EB_{\cI_\bullet} \cX$ given by the complement of \[V((\pi\circ p)^{-1}(\pi\circ p)_*(\bigoplus_{n\in \bZ} \Bbbk_{\theta}^{\otimes n})) = V((\pi\circ p)^{-1}\pi_*(\bigoplus_{n\in \bZ} \cI^n))=V(p^{-1}\pi^{-1}\pi_*(\bigoplus_{n\in \bZ} \cI^n))=:\cZ.\]
Since $\pi^{-1}\pi^*\cI_n\subseteq \cI^n$, we have that $V(\bigoplus_{n\ge 1}\cI_n)\subseteq \cZ$ so the latter is the locus in the weighted blow-up of $\cI_\bullet$ given by the complement of $V(\pi^{-1}\pi_*(\bigoplus_{n\in \bZ} \cI_n))$, namely $\Bl^\pi_\cI\cX$.
\end{proof}
\begin{proof}[Proof of \Cref{thm_you_can_embed_a_Stack_in_one_with_open_dm_by_taking_deformations_to_nc}] By the main theorem in \cite{ER21} there is a sequence of saturated blow-ups $\cX_n\to \cX_{n-1}\to \ldots \to \cX_1\to \cX$ such that $\cX_n$ is Deligne--Mumford. As taking an extended blow-up commutes with open embeddings, it suffices to prove by induction on $n$ that there is an algebraic stack $\widetilde{\cX}_j$ such that:
\begin{enumerate}
    \item there is an open embedding $i_j:\cX\to\widetilde{\cX}_j$ with a projection $\pi_j:\widetilde{\cX}_j\to \cX$ inducing isomorphisms on good moduli spaces and such that $\pi_j\circ i_j$ is isomorphic to the identity,
    \item if $U\subseteq X$ consists of the open subset of properly stable points in $X$, then $\pi_j$ is an isomorphism on the preimage of $U$,
    \item there is a morphism $\xi_j:\widetilde{\cX}_j\to \cB\Gm^{\oplus j}$ and a character $\theta_j:\Gm^{\oplus j}\to \Gm$ such that if we denote by $\cL_j:=\xi_j^*\Bbbk_{\theta_j}$, the semistable locus for $\cL_j$ over $X$ is isomorphic to $\cX_j$.
\end{enumerate}
For the case $n=1$ we can take the extended blow-up corresponding to the saturated blow-up $\cX_1\to \cX$, and apply \Cref{lemma_ext_blowup_is_relative_ss_locus}: all the points above follow from the properties of extended blow-ups and from \Cref{lemma_ext_blowup_is_relative_ss_locus}.

Assuming that the desired result holds for $j$, we show it holds for $j+1$. Let $\cZ_j\subseteq \cX_j$ the closed locus we would blow-up to obtain $\cX_{j+1}$, and let $\widetilde{\cZ}_j$ be its closure in $\widetilde{\cX}_j$. As $\cX_j$ is open in $\widetilde{\cX}_j$, the extended blow-up $\tau:\EB_{\widetilde{\cZ}_j}\widetilde{\cX}_j\to \widetilde{\cX}_j$ restricts to $\EB_{\cZ_j}\cX_j$. From \Cref{lemma_ext_blowup_is_relative_ss_locus} there is a morphism $\EB_{\cZ_j}\cX_j\to \cB \Gm$ and a character $\theta_{j+1}$ such that the semistable locus over the good moduli space $X_j$ of $\cX_j$ for $\Bbbk_{\theta_{j+1}}$ is $\cX_{j+1}$ . The morphism $\EB_{\cZ_j}\cX_j\to \cB \Gm$ comes from the $\Gm$-quotient presentation of the definition of extended blow-up, so it extends to a morphism $\EB_{\widetilde{\cZ}_j}\widetilde{\cX}_j\to \cB \Gm$. There is a morphism $\iota:\widetilde{\cX}_j\to \EB_{\widetilde{\cZ}_j}\widetilde{\cX}_j$ from the properties of extended blow-ups, and it is straightforward to check that points (1) and (2) follow by taking $i_{j+1}:=\iota\circ i_i$ and $\pi_{j+1} = \pi_j\circ \tau$. Point (3) instead follows from \Cref{prop_stable_locus_of_pi*L^n_otimes_G}, we explain how. By induction we have a line bundle $\cL$ on $\widetilde{\cX}_j$ whose semistable locus over $X$ is $\cX_j\subseteq \widetilde{\cX}_j$, the $j$-th step of the algorithm of Edidin and Rydh.
We know that the $(j+1)$-th
step $\cX_{j+1}\to \cX_j$ is a saturated blow-up, which by \Cref{lemma_ext_blowup_is_relative_ss_locus} is the relative semistable locus, over the good moduli space $X_j$ of $\cX_j$, for a line bundle $\cM$ of the extended blow-up $\EB_{\cZ_{j+1}}\cX_{j+1}$. Then \Cref{prop_stable_locus_of_pi*L^n_otimes_G} allows us to take both semistable loci simultaneously, so $\cX_{j+1}$ will be the $\pi^*\cL^{\otimes m}\otimes \cM$-semistable locus of $\widetilde{\cX}_{j+1}$ over $X$, which proves (3). 

Finally, as the weighted blow-up $\EB_\cI \cX$ admits a representable (in fact, affine) morphism $\EB_\cI\to \cX\times \cB \Gm$, if $\cX$ is a global quotient then also $\widetilde{\cX}_i$ is a global quotient for every $i$.
\end{proof}

\section{An application of \Cref{thm_intro_simplified} and  \Cref{thm_you_can_embed_a_Stack_in_one_with_open_dm_by_taking_deformations_to_nc}}\label{section_examples_extended}
\Cref{thm_intro_simplified} now follows from the results in \Cref{sec:qmaps} and \Cref{section_extended_blowups}: from \Cref{section_extended_blowups} given an algebraic stack $\cX$ with a dense properly stable locus, one compactify the space of maps to $\cX$ as follows. Either the assumptions of \Cref{teo_intro_cX_contains_open_proper_dm} apply, or one can enlarge $\cX\subseteq \widetilde{\cX}$ using \Cref{thm_you_can_embed_a_Stack_in_one_with_open_dm_by_taking_deformations_to_nc}, so that the assumptions of \Cref{teo_intro_cX_contains_open_proper_dm} apply to $\widetilde{\cX}$. One can use the inclusion $\cX\subseteq \widetilde{\cX}$ and \Cref{teo_intro_cX_contains_open_proper_dm} to compactify the space of maps to $\cX$, and the projection morphism $\widetilde{\cX}\to \cX$ to give a modular interpretation to maps to $\widetilde{\cX}$ in terms of maps to $\cX$.

In this section we give an example of how to use the \Cref{thm_you_can_embed_a_Stack_in_one_with_open_dm_by_taking_deformations_to_nc} and \Cref{teo_intro_cX_contains_open_proper_dm} to compactify the space of maps to the GIT moduli space of plane cubics, which we denote by $\sP^1$. Recall that this GIT moduli space is defined as $\sP^1:= [\bP(\oH^0(\cO_{\bP^2}(3)))^{ss}/\PGL_3]$, where $\bP(\oH^0(\cO_{\bP^2}(3)))^{ss}$ is the open in $\bP(\oH^0(\cO_{\bP^2}(3)))$ parametrizing plane cubics with at most nodal singularities. This algebraic stack both does not have an open substack which is proper and Deligne--Mumford, and it does not have (as far as we know) a modular enlargement containing a proper and Deligne--Mumford open substack.
We show that also in this case, from the properties of extended blow-ups, one can still carry a detailed analysis of the objects on the boundary of $\cQ_g(\widetilde{\sP}^1,\widetilde{\sP}^1_\dm,\beta)$ for an extended blow-up $\widetilde{\sP}^1\to \sP^1$, by post-composing a map $\cC\to\widetilde{\sP}^1$ with $\cC\to\widetilde{\sP}^1\to \sP^1$. In \Cref{prop_conditions_induced_by_phi_tilda_on_phi_and_f} we list some of the properties of maps $C\to \widetilde{\sP}^1$ appearing as degenerations in the quasimaps moduli space of maps $\bP^1\to \bP(\oH^0(\cO_{\bP^2}(3)))^{ss}/\PGL_3]$ coming from pencils of cubics.

\subsection{GIT moduli space of plane cubics}We begin by reporting a few classical results, which one might find in \cites{MFK, ascher2023moduli}
\begin{enumerate}
    \item each smooth cubic $C$ is GIT-stable,
    \item each nodal cubic can be degenerated, along $[\bA^1/\Gm]$, to $x_0x_1x_2$
    \item each smooth plane cubic can be written as $x_0^3+ x_1^3 + x_2^3 + \lambda x_0x_1x_2$, \cites{nakamura2004planar, hesse1897ludwig}, so each plane cubic contains $\bmu_3\rtimes S_3$ as automorphism group,
    \item the stabilizer of $x_0x_1x_2$ is $\Gm^3/{\Gm}_{, \Delta}\rtimes S_3$ where ${\Gm}_{, \Delta}$ is the diagonal subgroup of $\Gm^3$ and the $S_3$ action permutes the coordinates $x_0, x_1, x_2$. This is straightforward once we realize that any automorphism of $\PGL_3$ which fixes $x_0x_1x_2=0$ fixes the three nodes $[1,0,0]$, $[0,1,0]$ and $[0,0,1]$,
    \item the GIT quotient $\bP(\oH^0(\cO_{\bP^2}(3)))/\!\!/\PGL_3$ is isomorphic to $\bP^1$.
\end{enumerate}
In particular, we introduce the following
\begin{Notation}
    We denote by $G:=\Gm^3/{\Gm}_{, \Delta}\rtimes S_3$, the stabilizer of $x_0x_1x_2$, and by $\cG:=\Gm^3\rtimes S_3$ the group which surjects to $G$ in the obvious way. We further denote $\sP^1:= [\bP(\oH^0(\cO_{\bP^2}(3)))^{ss}/\PGL_3]$, and if $\cB G\hookrightarrow \sP^1$ is the inclusion of the point corresponding to $x_0x_1x_2$, we assume that it corresponds to the point at infinity in $\bP^1$.
\end{Notation}

\begin{Lemma}\label{lemma_luna_slice_curlyP1}
    There is a Luna slice for $\sP^1$ at $x_0x_1x_2$ isomorphic to $[\bA^3/G]$ where the action of $G$ on $\bA^3$ comes from an action of $\cG$ defined as follows:
    \begin{align*}&((\lambda,\mu,\nu),\sigma)\in \Gm^3\rtimes S_3\text{ acts as }\\((\lambda,\mu,\nu),\sigma)*&(x_1,x_2,x_3) = \left(\frac{\lambda^2}{\mu\nu}x_{\sigma^{-1}(1)}, \frac{\mu^2}{\lambda\nu}x_{\sigma^{-1}(2)}, \frac{\nu^2}{\lambda\mu}x_{\sigma^{-1}(3)}\right).
    \end{align*}
\end{Lemma}
\begin{proof} We denote by $x$ the point $x_0x_1x_2$.
It suffices to study the action of $\cG$ on $T_x$, the tangent space of $x\in \bP(\oH^0(\cO_{\bP^2}(3)))$, and the image of the map on tangent spaces $T_{\Id}\PGL_3\to T_x$ given by the inclusion of the orbit of $x$, and 
where $\Id$ is the identity in $\PGL_3$. Recall that we can identify the tangent space of the identity of $\PGL_3$ as the matrices as follows
\[
\begin{bmatrix}1+\epsilon_1&\epsilon_4&\epsilon_7\\\epsilon_2&1+\epsilon_5&\epsilon_8\\\epsilon_3&\epsilon_6&1+\epsilon_9
\end{bmatrix}\text{ such that }\epsilon_1+\epsilon_5+\epsilon_9=0\text{ and }\epsilon_i\epsilon_j=0\text{ for every }i,j.\]
Then the image of $T_{\Id}\PGL_3\to T_xX$ can be identified with the following forms \begin{align*}\{((1+\epsilon_1)x_0+\epsilon_2x_1+\epsilon_3x_2)&(\epsilon_4x_0+(1+\epsilon_5)x_1+\epsilon_6x_2)(\epsilon_7x_0+\epsilon_8x_1+(1+\epsilon_9)x_2):\\&\epsilon_i\epsilon_j=0, \epsilon_1+ \epsilon_5+\epsilon_9=0\}.\end{align*}
The above expression simplifies to
\[x_0x_1x_2 + \epsilon_7x_0^2x_1 +\epsilon_4x_0^2x_2+\epsilon_8x_0x_1^2 +\epsilon_2x_1^2x_2 + \epsilon_6x_0x_2^2 + \epsilon_3x_1x_2^2.\]
In particular, in the affine chart of $\bP(\oH^0(\cO_{\bP^2}(3)))$ given by $x_0x_1x_2\neq 0$, we can take as coordinates for $T_xX$ the following ones:
\[\frac{x_0^3}{x_0x_1x_2}, \frac{x_1^3}{x_0x_1x_2}, \frac{x_2^3}{x_0x_1x_2}, \frac{x_0^2x_1}{x_0x_1x_2},
\frac{x_0^2x_2}{x_0x_1x_2},
\frac{x_1^2x_0}{x_0x_1x_2},
\frac{x_1^2x_2}{x_0x_1x_2},
\frac{x_2^2x_0}{x_0x_1x_2},
\frac{x_2^2x_1}{x_0x_1x_2}.
\] Then the image of $T_{\Id}\PGL_3\to T_xX$ can be identified with the linear subspace generated by \[
<\frac{x_0^2x_1}{x_0x_1x_2},
\frac{x_0^2x_2}{x_0x_1x_2},
\frac{x_1^2x_0}{x_0x_1x_2},
\frac{x_1^2x_2}{x_0x_1x_2},
\frac{x_2^2x_0}{x_0x_1x_2},
\frac{x_2^2x_1}{x_0x_1x_2}>.
\] Since the linear subspace $<\frac{x_0^3}{x_0x_1x_2}, \frac{x_1^3}{x_0x_1x_2}, \frac{x_2^3}{x_0x_1x_2}>$ is $G$-invariant, we can take it as our $W$. It is straightforward to check now that the action of $G$ on $W$ is the desired one.  
\end{proof}
The advantage of our local description is that we can explicitly determine how many extended blow-ups are needed to apply \Cref{thm_you_can_embed_a_Stack_in_one_with_open_dm_by_taking_deformations_to_nc}. In particular:
\begin{Lemma}
Consider $\bA^3$ with the action of $G$ as in \Cref{lemma_luna_slice_curlyP1}.
The extended blow-up of the origin $\EB_0\bA^3$ contains a Deligne--Mumford stack $(\EB_0\bA^3)_\dm$ which is proper over the good moduli space of $\EB_0\bA^3\to \bA^1$, and the stabilizer of the unique closed point of $(\EB_0\bA^3)_\dm$ over $0\in \bA^1$ is $\bmu_3^{\otimes 3}/\bmu_{3,\Delta}\rtimes S_3$.
\end{Lemma}
\begin{proof}Recall that we can identify $\EB_0\bA^3$ with $[\bA^4/G\times \Gm]$ where the action of $G$ on the first three components is the original one and trivial on the last component. Similarly, the action of $\Gm$ on the first three components is with weight $1$ and on the last component is with weight $-1$. It is now straightforward to check that the point $(1,1,1,0)$ has $\bmu_3^{\otimes 3}/\bmu_{3,\Delta}\rtimes S_3$ as stabilizers. 
\end{proof}
Since the formation of extended blow-ups commutes with \'etale base change, we have
\begin{Cor}
    The extended blow-up of $\cB G \subseteq \sP^1$, which we denote by $\widetilde{\sP}^1$, contains a
    proper Deligne--Mumford stack, which we denote by $\widetilde{\sP}^1_\dm$. The coarse moduli space of $\widetilde{\sP}^1_\dm$ is $\bP^1$. 
\end{Cor}
\subsection{Space of maps to $\sP^1$}
We now focus on the problem of compactifying maps to $\sP^1$, in a specific example. Consider a moduli space $\cU$ of pencils of plane cubics, i.e. of maps $\bP^1\to \bP(\oH^0(\cO_{\bP^2}(3)))$, up to $\PGL_3$ (namely, up to change of coordinates) and such that each map $\bP^1\to \bP(\oH^0(\cO_{\bP^2}(3)))$ intersects the discriminant locus transversally (this is a $\PGL_3$-invariant condition). In particular:
\begin{enumerate}
    \item there is a fibration $\cC\to \bP^1$ where each fiber is a plane cubic, the generic fiber is smooth and the singular fibers have a single node,
    \item there are exactly 12 nodal singular fibers, as the discriminant is a hypersurface of degree 12 in $\bP(\oH^0(\cO_{\bP^2}(3)))$,
    \item if we denote by $\phi:\bP^1\to \sP^1$ the induced morphism, then $\phi^{-1}\cB G=\emptyset$ as the fibers of $\cC\to\bP^1$ have a single node,
    \item if we denote by $f:\bP^1\to \sP^1\to\bP^1$ the composition of $\phi$ with the good moduli space map, then $f$ has degree 12 and it is unramified over $\infty$ (the point corresponding to $x_0x_1x_2$).
\end{enumerate}
Then there is a map $\cU\to \cQ_{0}(\widetilde{\sP}^1,\widetilde{\sP}^1_\dm)$. Assume it is a monomorphism. 
\begin{Notation}\label{notation_sections}
    We denote by $\overline{\cU}$ the closure of $\cU$ inside $\cQ_{0}(\widetilde{\sP}^1,\widetilde{\sP}^1_\dm)$. Given a map $\psi:C\to\sP^1\to \widetilde{\sP}^1$ corresponding to a point in $\cU$, there are 12 points in $C$ such that $\psi(p)\not \in \widetilde{\sP}^1_\dm$. We denote by \[\widetilde{\phi}:(C\xrightarrow{\pi} B;x_1,\ldots,x_{12})\to \sP^1\]
    the objects $\widetilde{\phi}:C\to\sP^1$ of $\overline{\cU}(B)$, such that there are 12 \textit{sections} $x_1,\ldots,x_{12}$ of $\pi$ that satisfy \[\widetilde{\phi}^{-1}(\widetilde{\sP}^1\smallsetminus\widetilde{\sP}^1_\dm) = \{x_1,\ldots,x_{12}\}.\] 
\end{Notation}
The following remark follows from the definition of stable quasimap (\Cref{def_stable_qmap}), and since the coarse moduli space of $\widetilde{\sP}^1_\dm$ is $\bP^1$.
\begin{Remark}
    Let $\widetilde{\phi}:\cC\to \widetilde{\sP}^1$ a morphism in $\overline{\cU}$ inducing $f:C\to \bP^1$ on good moduli spaces. Then if $\{x_1,\ldots,x_m\}=\widetilde{\phi}^{-1}(\widetilde{\sP}^1\smallsetminus\widetilde{\sP}^1_\dm)$, for every $0<\epsilon $ the line bundle $\omega_C(\epsilon\sum x_i)\otimes f^*\cO_{\bP^1}(3)$ is ample on $C$.
\end{Remark}

\begin{Lemma}\label{lemma_irred_cp_which_map_to_infty_map_to_exc_locus}
     Let $R$ be a DVR with generic point $\eta$ and closed point $p$ and consider the pair $(C\to \spec(R);x_1,...,x_{12})\xrightarrow{\widetilde{\phi}}\widetilde{\sP}^1$ corresponding to a morphism $\spec(R)\to \overline{\cU}$ witn $12$ sections as in \Cref{notation_sections}. Assume that $\eta\to\overline{\cU}$ maps into $\cU$, and let $D\subseteq C_p$ be an irreducible component of the central fiber of $C\to \spec(R)$. Then:
     \begin{enumerate}
         \item if $f:D\to C\to \bP^1$ is finite and $x_i\in D$ is a marked point, $\widetilde{\phi}(x_i)$ is not contained in the exceptional divisor of $\widetilde{\sP}^1$,
         \item if $D\to C\to \bP^1$ maps to a point $z$, then $D$ is contained in the exceptional locus of $\widetilde{\sP}^1$ if and only if $z=\infty$, and
         \item if $f:D\to C\to \bP^1$ is finite $x\in D$ is a smooth point of $D$, then the multiplicity of $f^{-1}(\infty)$ at $x$ is $\#\{i:x_i=x\}$.         
     \end{enumerate}
\end{Lemma}
\begin{proof}
    Let $\cO(-1)$ be the line bundle on $\widetilde{\sP}^1$ given by the exceptional divisor, let $\cL=\widetilde{\phi}^*\cO(-1)$ and let $t$ be the section of $\cL$ given by the pull-back of the exceptional locus.
    By construction of the morphism $\eta\to \cU$, the section $t$ does not vanish along the generic fiber of $C$. Then $V(t)$ is, set theoretically, a union of irreducible components of the special fiber $C_p\subseteq C$ that don't map finitely to $\bP^1$, so point (1) follows since the markings $x_i$ are smooth on the special fiber.
    
    For (2), since the exceptional locus is contained in $\pi^{-1}\cB G$ from \Cref{def_gen_blowup_comes_with_a_cartier_divisor}, if $D\subseteq C_p$ maps to a point different from $\infty$ via $f$, then $D$ is not contained in $V(t)$. We are left with showing that if $f(D)=\infty$ then $t$ vanishes on $D$. By assumption, the generic point of $\widetilde{\phi}(D)$ maps to the Deligne--Mumford locus in $\widetilde{\sP}^1$, and there is a unique point of $\widetilde{\sP}^1_\dm$ over $\infty$. Such a point is in the exceptional locus, so $V(t)$ vanishes over it.

    Finally for (3) observe that the Cartier divisor $f^{-1}(\infty)$ is of the form $\cO_C(a_1x_1+\ldots +a_{12}x_{12})\otimes \cO_C(\Delta)$, where $\Delta$ is effective and supported on the irreducible components of the special fiber mapping to $\infty$ and $a_i\ge 1$. So it suffices to prove that $a_i=1$ for every $i$, which can be checked on the generic fiber. This follows by assumption, as the maps in $\cU$ intersect the discriminant locus transversally.
\end{proof}
We will need the next auxiliary lemma
\begin{Lemma}\label{lemma_when_a_map_from_a_dvr_to_An_lifts_to_the_blow_up_without_the_proper_transform_of_the_axis}
    Let $R$ be a DVR with generic point $\eta$ and closed point $p$, let $\Gm^r$ be a torus acting on $\bA^n$ diagonally, let $\pi:B_0[\bA^n/\Gm^r]\to [\bA^n/\Gm^r]$ be the blow-up of $[0/\Gm^r]\hookrightarrow [\bA^n/\Gm^r]$.  Let $\phi:\spec(R)\to [\bA^n/\Gm^r]$ be a morphism sending the generic point $\eta\in\spec(R)$ to a point in $[\bA^n/\Gm^r]$ which does not lie in $[\{x_i=0\}/\Gm^r]\subseteq [\bA^n/\Gm^r]$ for any $i$, where $x_i$ are the coordinates on $\bA^n$, and assume that $\phi(p)$ maps instead to $[0/\Gm^r]$. Then $\phi$ lifts uniquely to a morphism $\widetilde{\phi}:\spec(R)\to B_0[\bA^n/\Gm^r]$. Moreover, if we denote by $x_1,...,x_n$ the coordinates in $\bA^n$, we have that $\widetilde{\phi}(p)$ lies along the closed substack of $B_0[\bA^n/\Gm^r]$ given by the proper transform of the axis in $[\bA^n/\Gm^r]$ if and only if the there are $i,j$ such that $\phi^*x_i$ and $\phi^*x_j$ vanish at $p$ with different orders.
\end{Lemma}
\begin{proof}
    The existence of $\widetilde{\phi}$ follows from the valuative criterion for properness since $B_0[\bA^n/\Gm^r]\to [\bA^n/\Gm^r]$ is representable and proper, and an isomorphism away from $[0/\Gm^r]$. For the moreover part, one can identify the extended blow-up $\EB_0[\bA^n/\Gm^r]\to [\bA^n/\Gm^r]$ with $[\bA^{n+1}/\Gm^{r+1}]$ as in \Cref{example_ext_blowup_of_origin_in_A2} for the case $n=2$, and if $\bA^{n+1}$ has coordinates $X_1,...,X_n,T$, the map $[\bA^{n+1}/\Gm^{r+1}]\to [\bA^{n}/\Gm^{r}]$ sends $x_i\mapsto X_iT$. A morphism $\spec(R)\to [\bA^n/\Gm^r]$ corresponds to $r$ line bundles $\cL_1,\ldots,\cL_r$, together with a section $\phi^*x_i\in \oH^0(\spec(R),\cG_i)$ for $1\le i\le n$, where each $\cG_i$ is a linear combination of $\cL_1,\ldots,\cL_r$.
    If we denote by $m_i$ is the order of vanishing of $\phi^*x_i$ for every $i$ and $m=\min\{m_i\}$, the sections $\phi^*x_i$ lift uniquely to sections of $\widetilde{x}_i\in \oH^0(\spec(R),\cG_i\otimes \cO_{\spec(R)}(-mp))$. One can check that the lift $\widetilde{\phi}$ comes from the line bundles $\{\cG_i\otimes \cO_{\spec(R)}(-mp)\}_i\cup\{\cO_{\spec(R)}(mp)\}$
    and the sections being the sections $\widetilde{x}_i$ and the one coming from dualizing the inclusion of the ideal sheaf $\cO_{\spec(R)}(-mp)\to \cO_{\spec(R)}$. Now the desired statement follows observing that $m<m_j$ is equivalent to the vanishing of $\widetilde{x}_j$ at $p$.
\end{proof}
\begin{Prop}\label{prop_conditions_induced_by_phi_tilda_on_phi_and_f}
    Let $(C\to B;x_1,...,x_{12})\xrightarrow{\widetilde{\phi}}\widetilde{\sP}^1$ be a $B$-point in $\overline{\cU}$ as in \Cref{notation_sections}, let $\phi:C\to \sP^1$ be the composition of $\widetilde{\phi}$ and the extended weighted blow-up $\widetilde{\sP}^1\to \sP^1$, and let $f:C\to \bP^1$ the composition to the good moduli space map $\sP^1\to \bP^1$. Then:
    \begin{enumerate}
        \item there is a line bundle $\cL$ on $C$ with a section $t\in \oH^0(C,\cL)$ such that, for every $b\in B$ we have an equality in sets $V(t)_b = \{D\subseteq C_b$ irreducible component such that $f_b(D)=\infty\}$,
        \item if $D\subseteq C_b$ is an irreducible component such that $f_b(D)=\infty$, then the map $D\to \sP^1$ factors via $D\to \cB G\to \sP^1$,
        \item $f^*\cO_{\bP^1}(1)\cong \cO_\cC(x_0+ \ldots+ x_{12})\otimes \cL^{\otimes 3}$ and $f^{-1}(\infty)=V(x_0\cdot \ldots\cdot x_{12}\cdot t^3)$, where the latter is a section of $\cO_C(x_0+ \ldots+ x_{12})\otimes \cL^{\otimes 3}$,
        \item let $D\subseteq C_b$ be an irreducible component such that $f(D)$ is not a point and let $n\in D$ be a node which maps to $\infty$. Let $V\to D$ be an \'etale cover of $n$ which is a scheme, let $v\in V$ be a point mapping to $n$, and let $\cC\to V$ be the resulting family of plane cubics. Then the singularities of $\cC$ around each node of $\cC_v$ are isomorphic.
    \end{enumerate}
\end{Prop}
\begin{proof}
    If we denote by $\cL=\widetilde{\phi}^*\cO(1)$ and $t$ the pull-back of a section that vanishes along exceptional locus, then point (1) follows from \Cref{lemma_irred_cp_which_map_to_infty_map_to_exc_locus}. Point (2) instead follows from \Cref{def_gen_blowup_comes_with_a_cartier_divisor}. Point (3) is true away from the vanishing locus of $t$ from \Cref{lemma_irred_cp_which_map_to_infty_map_to_exc_locus}, so it suffices to prove it in an \'etale neighbourhood of $t=0$. Similarly, point (4) can be checked on an \'etale neighbourhood of $t=0$, so we choose a neighbourhood where we will check both, as follows. We already computed a Luna slice around the polystable point of $\sP^1$ in \Cref{lemma_luna_slice_curlyP1},  so from \Cref{subsection_luna_slice} there are diagrams as follows  \[
    \xymatrix{[\widetilde{W}/G\times \Gm]\ar[d] \ar[r]^\gamma &[\bA^4/G\times \Gm]\ar[d]\\
    [W/G]\ar[d]\ar[r]^\beta  &[\bA^3/G]\ar[d]\\
    W/\!\!/G\ar[r]^\alpha& \bA^1}\text{ }\text{ }\text{ }\text{ }\text{ }\text{ }\xymatrix{[\widetilde{W}/G\times \Gm]\ar[d] \ar[r] & \widetilde{\sP}^1\ar[d]\\
    [W/G]\ar[d]^\pi\ar[r]^\iota  &\sP^1\ar[d]\\
    W/\!\!/G\ar[r]^i& \bP^1}
    \]
    with $i$ and $\alpha$ \'etale, with all squares cartesian, and with a point $w\in W/\!\!/G$ mapping to $0\in \bA^1$ and $\infty \in \bP^1$. Recall that $[\bA^4/G\times \Gm]\to [\bA^3/G]$ is the extended blow-up of the origin in $[\bA^3/\Gm]$.
    We can then pull-back the morphism $\phi:C\to \sP^1$ along $\iota$, this will give us an \'etale neighbourhood where we will check conditions (3) and (4). We denote by $U:=C\times_{\sP^1}[W/G]$ and by \[
        \phi_U:U\to [W/G]\to [\bA^3/G]\text{ and }\widetilde{\phi}_U:U\to [\widetilde{W}/G\times\Gm]\to [\bA^4/G\times \Gm]\] the composition of the two second projections with the maps $\beta$ and $\gamma$.

    To check condition (3), observe that if $x_1,x_2,x_3$ are the coordinates of $\bA^3$, then the ring of $G$-invariants consists of $\spec(k[x_1x_2x_3])$. Moreover if $X_1,X_2,X_3,T$ are the coordinates on $\bA^4$, the blowdown morphism $[\bA^4/G\times\Gm]\to [\bA^3/G]$ is induced by $x_i\mapsto X_iT$. Then the composition sends $x_1x_2x_3\to X_1X_2X_3T^3$, so if we denote by $\widetilde{\pi}:[\bA^4/G\times \Gm]\to \bA^1$ is the composition of the extended blow-down $[\bA^4/G\times \Gm]\to [\bA^3/G]$ and the good moduli space map $[\bA^3/G]\to \bA^1$, then \begin{equation}\label{equation}\widetilde{\pi}^{-1}(0)=[V(X_1X_2X_3T^3)/G\times \Gm]. \end{equation} In particular, as $V(T)$ is the exceptional divisor, the Cartier divisor $V(t)\subseteq U$ agrees with $V(\widetilde{\phi}_U^*T)$. Similarly, from \Cref{lemma_irred_cp_which_map_to_infty_map_to_exc_locus}, the Cartier divisor $V(\widetilde{\phi}_U^*(X_1X_2X_3))$ agrees with $V(p_1\cdot\ldots\cdot p_{12})$, where we denoted by $p_i$ the pull-backs of the markings on $C$ to $U$. Now point (3) follows from \Cref{equation}.

    To prove (4), we can take a further \'etale cover of $[\bA^3/G]$ as follows. Recall that $G=(\Gm^3/\Gmdelta)\rtimes S_3$, so one can take the $S_3$-torsor $[\bA^3/(\Gm^3/\Gmdelta)]\to [\bA^3/G]$. This gives a family of maps $\sC\to [\bA^3/(\Gm^3/\Gmdelta)]$, and over the origin one can check from the explicit description of the Luna slice in \Cref{lemma_luna_slice_curlyP1} that the directions $\{x_i=0\}\subseteq \bA^3$ correspond to degenerations of $V(x_0x_1x_2)\subseteq \bP^2$ along which a node does not smooth. Extended blow-ups commutes with flat base changes, so $[\bA^4/G\times \Gm]\times_{[\bA^3/\Gm]}[\bA^3/(\Gm^3/\Gmdelta)]\cong[\bA^4/(\Gm^3/\Gmdelta)\times \Gm]$. 
    Since $\widetilde{\phi}(n)$ maps to $\widetilde{\sP}^1_\dm$, $\widetilde{\phi}_U(n)$ does not map in $[V(X_1X_2X_3)/G\times \Gm]\subseteq [\bA^4/G\times \Gm]$. In other terms, $\widetilde{\phi}_U(n)$ maps to the unique point of $[\bA^4/G\times \Gm]$ which is in $\widetilde{\sP}^1_\dm\times_{\widetilde{\sP}^1}[\bA^4/G\times \Gm]$ and is over $\infty$, namely the point $(1,1,1,0)\subseteq [\bA^4/G\times \Gm]$. Then from \Cref{lemma_when_a_map_from_a_dvr_to_An_lifts_to_the_blow_up_without_the_proper_transform_of_the_axis} the order of vanishing of $\widetilde{\phi}^*_Ux_i$ does not depend on $i$. The conclusion now follows from \Cref{lemma_when_a_map_from_a_dvr_to_An_lifts_to_the_blow_up_without_the_proper_transform_of_the_axis}.
    \end{proof}
\appendix
\section{Extended blow-ups, rational pointed curves and Hassett's conjecture}\label{appendix}
In this appendix, we give a criterion for when a stack can be obtained as an extended weighted blow-up, and we use it in order to prove \Cref{thm_intro_hassett}.
\subsection{A criterion for being an extended blow-up}\label{subsection_criterion_for_ext_wblowups}
In this subsection we give a few criteria under which a morphism of algebraic stacks $\widetilde{\cX}\to \cX$ is an extended weighted blow-up; this will come in handy later.
\begin{Prop}\label{prop:criterion}
    Let $\pi\colon\widetilde{\cX}\to\cX$ be a morphism of algebraic stacks such that:
    \begin{enumerate}
        \item the stacks $\widetilde{\cX}$ and $\cX$ have the same good moduli space, i.e. there is a good moduli space $\rho:\cX\to X$ and the composition $\rho\circ\pi$ also defines a good moduli space for $\widetilde{\cX}$;
        \item the canonical map $\cO_\cX\to \pi_*\cO_{\widetilde{\cX}}$ is an isomorphism;
        \item there is a morphism $\widetilde{\cX}\to [\bA^1/\Gm]$ such that $\widetilde{\cX}\to \cX\times[\bA^1/\Gm]$ is representable;
        \item if $\widetilde{\cE}\subseteq\widetilde{\cX}$ is the Cartier divisor induced by $\widetilde{\cX}\to\cX\times [\bA^1/\Gm]\to [\bA^1/\Gm]$, we have $(\pi|_{\widetilde{\cE}})_*(\cN_{\widetilde{\cE}}^{\otimes d})=0$ for $d>0$, where $\cN_{\widetilde{\cE}}$ is the normal bundle of $\widetilde{\cE}$.
    \end{enumerate}
    Then $\pi\colon\widetilde{\cX}\to \cX$ is an extended weighted blow-up of $\cX$ along the filtration
    \[ \cO_\cX\supset\pi_*\cO(-\widetilde{\cE}) \supset \pi_*\cO(-2\widetilde{\cE}) \supset \ldots \supset \pi_*\cO(-n\widetilde{\cE}) \supset \ldots\]
\end{Prop}
\begin{proof}
    First, observe that (1) implies that $\widetilde{\cX}\to \cX$ is cohomologically affine \cite{Alp}*{Proposition 3.13}. 
    Let $\cY$ be the total space of the $\Gm$-torsor $g\colon\cY\to\widetilde{\cX}$ associated to $\widetilde{\cE}$. As $\cY\to\widetilde{\cX}$ is affine and $\widetilde{\cX}\to\cX$ is cohomologically affine, also $\cY\to\cX$ is cohomologically affine and representable from (3), hence affine \cite{Alp}*{Proposition 3.3}.

    Write $\cY=\spec_{\cO_\cX}(\cA)$. Observe that, as $g\colon\cY\to\widetilde{\cX}$ is the $\Gm$-torsor associated to $\widetilde{\cE}$, we have that $g_*\cO_{\cY}=\oplus_{n\in\bZ} \cO(-n\widetilde{\cE})$. We deduce that $\cA=\oplus_{n\in\bZ} \pi_*\cO(-n\widetilde{\cE})$.    We prove that $\pi_*\cO(n\widetilde{\cE})=\cO_{\cX}$ for $n\geq 0$.
    
    Observe that by (2) the map $\cO_\cX \to \pi_*\cO_{\widetilde{\cX}}$ is an isomorphism, and consider the exact sequence
    \[ 0 \to \cO_{\widetilde{\cX}}(n\widetilde{\cE}) \to \cO_{\widetilde{\cX}}((n+1)\widetilde{\cE}) \to j_* \cN_{\widetilde{\cE}}^{\otimes (n+1)} \to 0. \]
    Pushing it forward along $\pi$, we obtain an isomorphism $\pi_*\cO_{\widetilde{\cX}}(n\widetilde{\cE})\simeq \pi_*\cO_{\widetilde{\cX}}((n+1)\widetilde{\cE})$, because the third term of the sequence vanishes by (4). We can therefore conclude by induction that $\pi_*\cO(n\widetilde{\cE})\simeq\cO_{\cX}$ for $n\geq 0$.

    Set then $\cI_n:=\pi_*\cO_{\widetilde{\cX}}(-n\widetilde{\cE})$ for $n\in\bZ$. Then by pushing down the filtration
    \[ \cO_{\widetilde{\cX}}\supset \cO_{\widetilde{\cX}}(-\widetilde{\cE}) \supset \ldots \supset \cO_{\widetilde{\cX}}(-n\widetilde{\cE}) \supset \ldots \]
    we get a filtration 
    \[ \cO_{\cX}\supset \cI_1 \supset \cI_2 \supset \ldots \supset \cI_n \supset \ldots. \]
    One can check that it satisfies the conditions of \Cref{def_weighted_embedding},
    hence $\cA=\oplus_{n\in\bZ} \cI_n$ and $\widetilde{\cX}=\EB_{\cI_\bullet}\cX$.
\end{proof}
\begin{Lemma}\label{lm:iso of O}
    Let $\pi:\widetilde{\cX}\to \cX$ be a morphism of algebraic stacks,
    assume that there is a schematically dense open immersion $i:\cX\hookrightarrow \widetilde{\cX}$ satisfying $\pi\circ i = \id_{\cX}$. Then $\pi_*\cO_{\widetilde{\cX}}=\cO_\cX$.
\end{Lemma} 
For example, if $i$ is topologically dense and $\widetilde{\cX}$ is reduced then $i$ is schematically dense.
\begin{proof}
    The morphism $\cO_{\widetilde{\cX}}\to i_*\cO_{\cX}$ is injective, since $i$ is schematically dense. So also $\pi_*\cO_{\widetilde{\cX}}\to \pi_*i_*\cO_{\cX}$ is injective. Since the composition $\cO_{\cX} \longrightarrow \pi_*\cO_{\widetilde{\cX}} \longrightarrow \pi_*i_*\cO_{\widetilde{\cX}}=\cO_{\cX}$
    is surjective, the morphism $\pi_*\cO_{\widetilde{\cX}}\to \pi_*i_*\cO_{\cX}$ is also surjective. This implies that $\cO_\cX \to \pi_*\cO_{\widetilde{\cX}}$ is an isomorphism as claimed.
\end{proof}
\begin{Cor}\label{cor:criterion 2}
    In the setting of \Cref{prop:criterion}, suppose that (1) and (3) hold, together with 
    \begin{itemize}
        \item[(5)] the stack $\cX$ is normal, $\widetilde{\cX}$ is integral, and 
       there is a schematically dense open immersion $i:\cX\hookrightarrow \widetilde{\cX}$ satisfying $\pi\circ i = \id_{\cX}$ such that $\pi(\widetilde{\cE})\subset \cX$ has codimension $\geq 2$.
    \end{itemize}
    Then the same conclusion of \Cref{prop:criterion} holds true.
\end{Cor}
\begin{proof}
    Set $\cZ:=\pi(\widetilde{\cE})$ and observe that by (5) and \Cref{lm:iso of O} there is an isomorphism $j:\cO_{\cX} \overset{\sim}{\to} j_*\cO_{\cX\smallsetminus\cZ}$. Consider the morphism $\alpha\colon\cO_{\cX}\simeq \pi_*\cO_{\widetilde{\cX}} \to \pi_*\cO(n\widetilde{\cE})$ for $n>0$. This is injective since $\cO_{\widetilde{\cX}}\to \cO_{\widetilde{\cX}}(n\widetilde{\cE})$ is injective, as $\widetilde{\cX}$ is integral. We now show that it is surjective.
    As $\widetilde{\cX}$ is integral, we have an injection $ \cO_{\widetilde{\cX}}\to \widetilde{j}_*\cO_{\widetilde{\cX}\smallsetminus \widetilde{\cE}}$ given by the restriction to the complement of $\widetilde{\cE}$. So:
    \[0\longrightarrow\cO(n\widetilde{\cE}) \longrightarrow \cO(n\widetilde{\cE})\otimes \widetilde{j}_*\cO_{\widetilde{\cX}\smallsetminus\widetilde{\cE}} \simeq \widetilde{j}_*\widetilde{j}^*\cO(n\widetilde{\cE}) \simeq \widetilde{j}_*\cO_{\widetilde{\cX}\smallsetminus \widetilde{\cE}}.\] 
    By pushing forward along $\pi$ we get an injection 
    \[0 \longrightarrow \pi_*\cO(n\widetilde{\cE}) \longrightarrow \pi_*\widetilde{j}_*\cO_{\widetilde{\cX}\smallsetminus\widetilde{\cE}} \simeq j_*\cO_{\cX\smallsetminus\cZ} \simeq \cO_{\cX}.\]
    This provides a section to $\alpha$, hence $\alpha$ is surjective. Now the same proof of \Cref{prop:criterion} goes through.
\end{proof}
\subsection{Rational pointed curves}\label{subsection_pts_on_p1}

In this subsection we study a \textit{modular} enlargement $\cC^\git_{2n} \subseteq \widetilde{\cC}^{\CY}_{2n}$ of the GIT moduli stack of genus 0 curves with a divisor of degree $2n$, which we denote by $\cC^\git_{2n}$, and which we show is an extended weighted blow-up.  We will use results from \Cref{section_extended_blowups} to prove \Cref{thm_intro_hassett}. In \Cref{subsubsection_ordered_pts} we study the problem for the case of $2n$ \textit{ordered} points on $\bP^1$.
\subsubsection{Unordered points on $\bP^1$}
For $n\geq 1$, let $\cC_{2n}^{\git}$ be the stack of smooth rational curves endowed with a divisor $D=p_1+\ldots+p_{2n}$ such that if $p_{i_1}=\ldots=p_{i_m}$, then $m\leq n$. This stack is isomorphic to the smooth quotient stack $[\bP \oH^0(\bP^1,\cO(2n))^{ss}/\PGL_2]$, where the $\PGL_2$-action is linearized via the line bundle $\cO(1)$, and it admits a projective good moduli space \cite{MFK}*{Proposition 4.1}.

The good moduli space $\cC_{2n}^{\git}\to C_{2n}^{\git}$ is properly stable, but it is not a coarse moduli space due to the presence of strictly semistable points in $\bP \oH^0(\bP^1,\cO(2n))^{ss}$: the point in the stack corresponding to $\bP^1$ and $D=n\cdot 0 + n\cdot \infty$ corresponds to the unique closed orbit of the strictly semistable locus.

Consider the stack $\Ctilde_{2n}^{\CY}$ whose objects are pairs $(P\to S, D)$ where:
\begin{enumerate}
    \item the morphism $P\to S$ is a twisted conic \cite{DLV}*{\S 1.1}, i.e. a proper, flat morphism whose geometric fibers are rational twisted curves having at most one node, and such node have automorphism group isomorphic to $\bmu_2$;
    \item the divisor $D\subset P$ is contained in the schematic locus of $P$, is finite of degree $2n$ over $S$, and on the geometric fibers $D=p_1+
    \ldots +p_{2n}$ has degree at least $n$ on each irreducible component, and if $p_{i_1}=\ldots =p_{i_m}$, then $m\leq n$.
\end{enumerate}
We aim at proving the following.
\begin{Teo}\label{thm:CY extended blow up}
    The following hold true:
    \begin{enumerate}
        \item there is a morphism $\pi\colon \Ctilde^{\CY}_{2n}\to\cC^{\git}_{2n}$ which realizes the first stack as an extended weighted blow-up of the second;
        \item the stack $\Ctilde^{\CY}_{2n}$ contains a proper Deligne--Mumford stack $\overline{C}_{2n}^{\CY}$ admitting a projective coarse moduli space.
    \end{enumerate}
\end{Teo}
In particular, one can compactify the space of maps to $\cC_{2n}^\git$ by considering maps to $\Ctilde^{\CY}_{2n}$ as the enlargement of \Cref{thm_intro_simplified}.
So while $\cC_{2n}^\git$ does not contain an open
substack which is Deligne--Mumford, one can choose an enlargement of $\cC_{2n}^\git$ 
(namely, $\Ctilde^{\CY}_{2n}$) which is itself a moduli space, and which contains an open substack which is proper and Deligne--Mumford.

Recall that there exists a smooth algebraic stack $\cC^{\CY}_{2n}$ whose objects are the same as the ones of $\Ctilde_{2n}^{\CY}$ albeit the node of the singular conics is not twisted \cite{ascher2023moduli}*{Definition 16.7, Lemma 16.8}. This stack has a natural forgetful morphism $\cC^{\CY}_{2n} \to \cC$, where $\cC$ is the stack of rational curves having at most one node.
Let $\cR$ be the stack of twisted conics \cite{DLV}*{Definition 1.1}, which also has a forgetful morphism $\cR\to \cC$. Then the following is straightforward.
\begin{Lemma}\label{lm:cartesian}
    We have $\Ctilde^{\CY}_{2n} \simeq \cC^{\CY}_{2n} \times_{\cC} \cR$.
\end{Lemma}
We will also need the following.
\begin{Lemma}\label{lm:root}
    The stack $\cR$ is an order $2$ root stack of $\cC$ along the smooth divisor of singular curves.
\end{Lemma}
\begin{proof}
    Observe that $\cC$ is isomorphic to the quotient stack $[\bP (\oH^0(\bP^2,\cO(2) \smallsetminus \Delta_2)/\PGL_3]$, where $\Delta_2$ is the locus of quadrics of rank $\leq 1$. Then the substack $\cC_1$ of singular curves is isomorphic to $[(\Delta_1\smallsetminus\Delta_2)/\PGL_3]$, which is a smooth divisor.
    
    Let $\Ctilde$ denote the root stack of $\cC$ along $\cC_1$. In order to define a morphism $\cR\to \Ctilde$, consider a twisted conic $P\to S$ and let $\overline{P}\to S$ be its coarse space, where $S=\spec(A)$ is smooth. By picking an atlas of $\cR$ and up to further shrinking it, we can assume that the locus of points with singular fibers is the vanishing locus of $f\in A$, where $f$ is not a zero divisor.
    
    \'{E}tale locally around the twisted node, the twisted conic is isomorphic to $\spec(A[u,v]/(uv-f)/\bmu_2]$, where the action is given by $u\mapsto (-u)$ and $v\mapsto (-v)$. The coarse space is instead isomorphic to $\spec(A[x,y]/(xy-f^2))$, and the coarse moduli space morphism sends $x\mapsto u^2$ and $v\mapsto v^2$. Let $S_1\subset S$ be the locus of points whose fiber in $\overline{P}\to S$ is singular: then we have just showed that $S_1$ has a root of order $2$, locally defined by the vanishing locus of $f$. In this way, we have defined a morphism $\cR \to \widetilde{\cC}$.

    On the other hand, given a family of rational curves $\overline{P}\to S=\spec(A)$, and $T=V(f)\subset S$ a divisor such that $2T=S_1$, we have that \'etale locally $\overline{P}$ looks like $\spec(A[x,y]/(xy-f^2))$, hence we can replace each of these \'etale local charts with the stack $[\spec(A[u,v]/(uv-f))/\bmu_2]$, and glue them back together, thus obtaining a twisted conic $P\to S$. This defines a morphism $\widetilde{\cC}\to \cR$, and it is straightforward to check that the two morphisms of stacks that we defined are one the inverse of the other.
\end{proof}
A combination of \Cref{lm:cartesian} and \Cref{lm:root} immediately gives us the following.
\begin{Prop}\label{prop:Ctilde is root}
    The stack $\widetilde{\cC}^{\CY}_{2n}$ is a $\bmu_2$-root stack of $\cC^{\CY}_{2n}$ along the divisor of singular curves.
\end{Prop}
\begin{Cor}\label{cor:Ctilde has gms}
    The stack $\widetilde{\cC}^{\CY}_{2n}$ is an algebraic stack of finite type and it admits a projective good moduli space, which is isomorphic to $\cC^{\git}_{2n}$.
\end{Cor}
\begin{proof}
    This follows from \Cref{prop:Ctilde is root} and the fact that $\cC^{\CY}_{2n}$ is an algebraic stack of finite type with projective good moduli space isomorphic to $C^{\git}_{2n}$ \cite{ascher2023moduli}*{Theorem 16.9}.
\end{proof}
We are ready to prove the result stated at the beginning of this section. We begin with the following Remark, where we explain the main ideas of the proof of \Cref{thm:CY extended blow up}.
\begin{Remark}\label{remark_why_we_need_twisted_nodes}
    Consider a one parameter family $\pi:(\cC,\Delta)\to \spec(R)$ in $\cC^{\CY}_{2n}$, and assume that the generic fiber of $\pi$ is smooth and the special one is singular. In other terms, the family $\cC\to \spec(R)$ is a degeneration of $\bP^1$ to a nodal union of two $\bP^1$s. To construct a map $\cC^{\CY}_{2n}\to \cC^{\git}_{2n}$, one might want to:
    \begin{enumerate}
        \item first blow-up the singular point on the central fiber of $\cC$ to get $B\to \cC$; let $E$ be the exceptional divisor.
\item Contract the proper transforms of the two irreducible components of the central fiber of $\cC$ in $B$, so that the push forward of the proper transform of $\Delta$ intersects the pushforward of $E$ in two points. 
    \end{enumerate}However, as it is, this approach does not work as the central fiber of $B\to \spec(R)$ is not reduced when, for example, the total space $\cC$ is smooth. On the other hand, a local computation shows that one can arrange $E$ to be reduced whenever the nodal point $x$ on the central fiber of $\cC$ is an $A_{2n+1}$ singularity. So we would like to only consider families when $x$ is an $A_{2n+1}$ singularity. The idea is to force the central fiber to be a \textit{twisted} curve with a $\bmu_2$-stabilizer at the node, so that when we take the corresponding coarse moduli space, the nodal point is indeed always an $A_{2n+1}$ singularity.
\end{Remark}
We will also need the following technical lemma.
\begin{Lemma}\label{lemma_one_can_check_rep_at_poli_pts}
    Let $f:\cX\to \cY$ be a $S$-complete morphism, with $\cX$ quasicompact. Let $p$ be a point in $\cX$ and $x$ a closed $k$-point in the closure of $p$. If $f$ is representable at $x$, then it is representable at $p$.
\end{Lemma}
If $\cX$ admits a separated good moduli space, then $f$ is automatically $S$-complete. So to check that a morphism of algebraic stacks with a separated good moduli space is representable, it suffices to check that it is representable at the closed points of the domain. The following proof was suggested to us by Andres Fernandez Herrero, improving an argument we had in a previous version of this paper.

\begin{proof}It suffices to show that the relative inertia $I_{\cX/\cY}$ is trivial at $p$. Consider $i:\spec(A)\to \cY$ a smooth morphism such that the map $x\to \cX\to \cY$ lifts to $x\to \spec(A)\to \cY$. Let $g\colon\cX':=\cX\times_{\cY}\spec(A) \to \cX$ be the induced map. As the relative inertia commutes with smooth base change, we have that $I_{\cX'/\spec(A)}$ is trivial at $x'$ for every $x'\in g^{-1}(x)$, and it suffices to check that there is a point $p'$ in the fiber of $p$ such that $I_{\cX'/\spec(A)}$ is trivial at $p'$. We can choose $x'$ and $p'$ so that $x'$ is a closed point and is in the closure of $p'$.  

Observe that $\cX'$ is S-complete, since being S-complete is stable under base change from \cite{AHH}*{Remark 3.40} and $\spec(A)$ is S-complete. Then $\cX'$ is locally linearly reductive from \cite{alper2019cartan} and from \cite{luna}. Then $\cX'$ has unpunctured inertia from \cite{AHH}*{Theorem 5.3}.

Consider a morphism $g:\spec(A)\to \cX'$, where $A$ is a Henselian local ring, and assume that $g$ sends the generic point of $\spec(A)$ to $p'$ and the closed one to $x'$. Consider the pull-back of the inertia $g^*\cI_{\cX'}$ on $\spec(A)$, which is an affine group scheme over $\spec(A)$. Our goal is to show that it is the trivial group scheme. By upper-semicontinuity of fiber dimension for group schemes, since $\cI_{\cX'}$ is trivial at $x'$ then $g^*\cI_{\cX'}\to \spec(A)$ is quasifinite. Since we are on a field of characteristic 0, the group
scheme is unramified (that is, the fibers of $g^*\cI_{\cX'}\to \spec(A)$ are \'etale). Then from the structure of affine unramified morphisms over strictly Henselian rings \cite{stacks-project}*{\href{https://stacks.math.columbia.edu/tag/00UY}{Tag 00UY}}, we have that $g^*\cI_{\cX'}=\cI_1\sqcup \cI_2$ where $\cI_1$ is a disjoint union of finitely many closed
subschemes of $\spec(A)$ and the image of $\cI_2\to \spec(A)$ does not contain the closed point of $\spec(A)$. By the definition of unpuncuted inertia and since $\cX'$ has unpuncuted inertia, we have $\cI_2=\emptyset$. Since the special fiber of $g^*\cI_{\cX'}$ is trivial, $\cI_1$ must consist of a single closed $Z \to \spec(A)$. Observe that $Z\cong \spec(R)$ (it is the identity section), so we have that $g^*\cI_{\cX'}$ is trivial as desired.
\end{proof}

\begin{Remark}
    It might be tempting to try to prove \Cref{lemma_one_can_check_rep_at_poli_pts} by proving that the locus where $f$ is representable is open. This is false: the action of $\Gm\rtimes \bmu_2$ on $\bA^2$, where $\Gm$ acts with weights $(1,-1)$ and $\bmu_2$ swaps the two axis, has $\{(x,y):xy=0\}\smallsetminus\{(0,0)\}$ as  locus with trivial stabilizers. So the locus in $[\bA^2/\Gm\rtimes \bmu_2]$ where $[\bA^2/\Gm\rtimes \bmu_2]\to \spec(k)$ is representable is not open. Similarly, the stack $\cX:=[\bA^2\smallsetminus\{(0,0)\}/\Gm\rtimes \bmu_2]$ is Deligne--Mumford, non-separated, and the locus where $\cX\to\spec(k)$ is representable is not open. 
\end{Remark}

\begin{proof}[Proof of \Cref{thm:CY extended blow up}]We divide the proof in several steps.

    \noindent
    \textbf{Step 1}: we first need to construct a morphism $\widetilde{\cC}^{\CY}_{2n} \to \cC^{\git}_{2n}$. Let $P\to S$ be a twisted curve with $S=\spec(A)$, and $\overline{P}\to S$ its coarse space. As already done before, we can assume that the locus of points with a singular fiber is the vanishing locus of $f\in A$, where $f$ is not a zero divisor.

    Set $\widetilde{P}\to S$ to be the blow-up of $P$ along the relative singular locus, with exceptional divisor $E$. We claim (1) that $\overline{P}\to S$ has reduced fibers and that (2) there is a canonical $S$-morphism $g\colon\overline{P} \to P'$ with $(P'\to S)$ a smooth rational curve, obtained by contracting the irreducible components of the singular fibers which are not $E$.

    Assuming (1), consider the line bundle $\cF=\omega_{\widetilde{P}/S}^\vee(-E)$: by construction this has degree two on the smooth fibers. Let $R_{\overline{s}} \cup E_{\overline{s}} \cup R'_{\overline{s}}$ be a singular fiber; then the restriction of $\cF$ has multidegree $(0,2,0)$. Let $g\colon P\to P'$ be the restriction onto the image of the $S$-morphism induced by $\cF$; one can verify that $g$ has the properties described in (2). Indeed, $\oH^1(\cF_s)=0$ for every $s\in S$, so the formation of $\Proj(\bigoplus_{n\ge 0}\pi_*(\cF^{\otimes n}))$ commutes with base change from cohomology and base change. Then we can check the desired statement fiberwise, where one can check it easily.

    To prove (1), observe that up to passing to an \'{e}tale-local neighbourhood of the relative singular locus, we can substitute $\overline{P}$ with $\spec(A[x,y]/(xy-f^2))$; we must have such a local description, because $\overline{P}$ is the coarse space of a curve with a twisted node.
    A quick computation shows that the blow-up along $(x,y,f)$ of $\spec(A[x,y]/(xy-f^2))$ is covered by the following three charts:
    \begin{itemize}
        \item[(i)] $\spec(A[u,v]/(uv-1))$ with $x=uf$ and $y=vf$;
        \item[(ii)]$\spec(A[w,y]/(f-wy)$ with $x=w^2y$;
        \item [(iii)] $\spec(A[w',x]/(f-w'x))$ with $y=w'^2x$.
    \end{itemize}
    In particular, the fiber over $S_1$ is reduced.

    Given $(P\to S, D)$ an object of $\Ctilde^{\CY}_{2n}$, we can therefore consider the smooth rational curve $P'\to S$ together with the divisor $g(D)$. It is immediate to check that the latter is still a divisor, it is finite over $S$ of degree $2n$, and if there is a geometric fiber such that $D_s=p_1+\ldots +p_{2n}$ and $p_{i_1}=\ldots=p_{i_m}$, then $m\leq n$. In other words, the pair $(P'\to S, g(D))$ is an object of $\cC^{\git}_{2n}$. In this way, we have defined $\pi\colon \Ctilde^{\CY}_{2n} \to \cC^{\git}_{2n}$.

    \noindent
    \textbf{Step 2}: in order to apply the criterion given in \Cref{cor:criterion 2}, we first observe that there is a schematically dense, open embedding $\cC^{\git}_{2n}\to \cC^{\CY}_{2n}$ \cite{ascher2023moduli}*{Lemma 16.8} and that $\pi$ induces an isomorphism of good moduli spaces. As the image of this embedding does not intersect the divisor of singular curves, we can lift this embedding to $\Ctilde_{2n}^{\CY}$, as the latter is a root stack over $\cC^{\CY}_{2n}$ by \Cref{lm:root}. It is still schematically-dense.

    We know from \Cref{cor:Ctilde has gms} that $\Ctilde^{\CY}_{2n}$ has a projective good moduli space, which is normal because $\Ctilde^{\CY}_{2n}$ is a root stack over $\cC^{\CY}_{2n}$ (see \Cref{prop:Ctilde is root}) and the latter is smooth stack \cite{ascher2023moduli}*{Lemma 16.8}. The morphism $\pi$ induces a morphism of good moduli spaces $C_{2n}^{\CY}\to C_{2n}^{\git}$, which has a section given by the morphism induced by $i$, which is an isomorphism (this follows from how we defined $i$ and \cite{ascher2023moduli}*{Theorem 16.9}). It follows that the morphism $\pi$ induces an isomorphism of good moduli spaces.

    \noindent
    \textbf{Step 3}: to apply \Cref{cor:criterion 2} we need to show that the morphism $\xi:\Ctilde^{\CY}_{2n} \to \cC^{\git}_{2n}\times [\bA^1/\Gm]$, given by $\pi$ and the Cartier divisor $\widetilde{\cE}\subset \Ctilde^{\CY}_{2n}$ of singular curves, is representable. From how $\pi$ is constructed, it suffices to check that $\xi$ is representable at the points of $\widetilde{\cC}_{2n}^\CY$ which correspond to singular curves (namely, $\widetilde{\cE}$). Let $\pi_i$ the two projections of $\cC^{\git}_{2n}\times [\bA^1/\Gm]$. As there is a unique polystable point $p\in \widetilde{\cE}$, namely the curve where the support of the divisor consists of two points, it from \Cref{lemma_one_can_check_rep_at_poli_pts} it suffices to check that $\xi$ is representable at this point. We recall a few results from \cite{DLV}*{\S 1.5} on $\Aut_{\Ctilde^{\CY}_{2n}}(p)$.
    
     Let $(C,E)$ be the pair in $\widetilde{\cE}$ corresponding to $p$, namely a nodal twisted cubic $C$ with a divisor $E$ consisting of two smooth points on $C$, one per irreducible component. If we denote by $D$ the nodal union of two $\bP^1$s attached at the point $[0,1]$, there is an action of $\bmu_2$ on $D$ which preserves the irreducible components and on each irreducible component it acts as $[a,b]\mapsto [-a,b]$. The stack quotient $[D/\bmu_2]$ has three stacky points, two of which are smooth. There is a map $[D/\bmu_2]\to C$ which rigidifies $[D/\bmu_2]$ at the two smooth stacky points; the corresponding points on $C$ are the support of the divisor $E$, and we still denote by $E$ the two smooth points on $D$ fixed by $\bmu_2$. To summarize:
     \begin{itemize}
         \item $D$ is the nodal union of two $\bP^1$, and $E $ are two smooth points on $D$, one per irreducible component,
         \item $\bmu_2$ acts on $D$, non-trivially on each irreducible component, fixing $E$, and
         \item $C$ is the quotient $[D/\bmu_2]$ rigidified at $[E/\bmu_2]$.
     \end{itemize}Then:
    \begin{enumerate}
\item $\Aut(D,E)\cong \Gm^2\rtimes \bmu_2$ with where $\Gm$ acts by scaling each branch, and $\bmu_2$ swaps the two branches. 
        \item There is a map $\Aut(D,E)\to \Aut([D/\bmu_2],[E/\bmu_2])$ that is surjective with kernel $(-1,-1,\Id)$, and an isomorphism $\Aut([D/\bmu_2],[E/\bmu_2])\cong\Aut(C,E)$ \cite{DLV}*{Pag. 10}.
        \item The action of $\Aut_{\widetilde{\cC}^{\CY}_{2n}}(p)=\Aut(C,E)$ on the fiber of $\cO_{\widetilde{\cC}^{\CY}_{2n}}(\widetilde{\cE})$ at $p$ (i.e., on $\Aut_{[\bA^1/\Gm]}(\pi_2\circ \xi (p))$) is computed via the composition $\Gm^2\rtimes \bmu_2\cong \Aut(D,E)\to \Gm$, $(a,b,\sigma)\mapsto ab$ \cite{DLV}*{Proposition A.4 and pag. 12}.
\item The map $\widetilde{\cC}^{\CY}_{2n}\to \cC^\git_{2n}$ sends the twisted conic $C$ to the $\bP^1$ obtained by performing the blow-up of the reduced stacky point of $C$ and taking the coarse moduli space of the exceptional divisor. In a neighbourhood of the node, the blow-up is the $\bmu_2$-quotient of $\Proj(k[x,y,X_0,X_1]/(X_0y-X_1x,xy))$. We study the chart where $X_0\neq 0$; the other chart is similar. In this chart the blow-up is $\spec(k[x,\frac{X_1}{X_0}]/x^2\frac{X_1}{X_0})$ and the action of $\bmu_2$ sends $x\mapsto -x$.
Hence the coarse moduli space is $\spec(k[x^2,\frac{X_1}{X_0}]/x^2\frac{X_1}{X_0})$. The exceptional divisor is $\spec(k[\frac{X_1}{X_0}])$. So the composition \[\Gm^2\rtimes \bmu_2\cong \Aut(D,E)\to \Aut(C,E)=\Aut_{\widetilde{\cC}^{\CY}_{2n}}(p)\xrightarrow{\pi_1\circ \xi} \Aut_{\cC^{\git}_{2n}}(\pi_1\circ \xi(p))\text{ sends }\]\[(a,b,\sigma)\mapsto \left(\frac{a}{b},\sigma\right).\]
    \end{enumerate}
Then from points (3) and (4), the composition
\begin{align*}
\Gm^2\rtimes \bmu_2\cong \Aut(D,E)\to &\Aut_{\widetilde{\cC}^{\CY}_{2n}}(p)\to \Aut_{\cC^{\git}_{2n}}(\pi_1\circ\xi(p))\times \Aut_{[\bA^1/\Gm]}(\pi_2\circ\xi(p))\\&\text{sends } (a,b,\sigma)\mapsto \left(\left(\frac{a}{b},\sigma\right),ab\right).
\end{align*}
The kernel of this map is generated by $(-1,-1,\Id)$ so from (2) the map 
\[\Aut(C,E)\cong \Aut_{\widetilde{\cC}^{\CY}_{2n}}(p)\to \Aut_{\cC^{\git}_{2n}}(\pi_1\circ\xi(p))\times \Aut_{[\bA^1/\Gm]}(\pi_2\circ\xi(p))\]
is injective, as desired.  

    \noindent
    \textbf{Step 4}: we are only left with proving that $\Ctilde_{2n}^{\CY}$ contains a proper Deligne--Mumford stack with a coarse moduli space. Consider $\cU_{2n}$ the complement of the closed substack whose objects are pairs $(P\to S, D)$ where on each geometric fiber the restriction of $D$ is supported on only two points.
    Observe finally that $\cU_{2n}$ is Deligne--Mumford, and it has a projective coarse moduli space. Indeed, it is a root stack of the KSBA-moduli space compactifying $(\bP^2,(\frac{1}{n}+\epsilon )D)$  where $D\subseteq \bP^1$ is a smooth divisor of degree $2d$.
\end{proof}

\subsubsection{Moduli of pointed rational curves}\label{subsubsection_ordered_pts}
Let $\Mcal_{0,2n}$ denote the moduli stack of smooth rational curves with $2n$ markings, for $n>1$, which is a quasi-projective variety.

One possible modular compactification of $\Mcal_{0,2n}$ can be constructed as follows: consider the scheme \[\bP {\rm H}^0(\bP^1,\cO(1))^{\times 2n} \simeq (\bP^1)^{\times 2n}\]parametrizing $2n$-uples of points on $\bP^1$. There is a natural $\PGL_2$-action on this variety, which can be linearized via the line bundle $L=\cO(2)\boxtimes\ldots\boxtimes\cO(2)$. 

Set $\Mbar_{0,2n}^\git=[(\bP^1)^{\times 2n}(L)^{ss}/\PGL_2]$: then this is an algebraic stack, with good moduli space the GIT quotient $\overline{M}_{0,2n}^{\git}:=(\bP^1)^{\times 2n}(L)^{ss}/\!\!/\PGL_2$. 

The semistable locus of this linearized $\PGL_2$-action on $(\bP^1)^{\times n}$ consists of the configurations with at most $n$ markings collapsing together, and the stable points are those with at most $n-1$ points collapsing together. The closed orbits of the strictly semistable locus are the semistable configurations with the markings supported on only two points (see the proof of \cite{MFK}*{Proposition 4.1}).

Therefore, the stack $\Mbar_{0,2n}^{\git}$ is a compactification of $\Mcal_{0,2n}$, but it does not have a proper, Deligne--Mumford substack with a coarse moduli space, so we cannot apply the theory of stable quasimaps of \S \ref{sec:qmaps}. 

Consider the stack $\Mtilde_{0,2n}^{\Hass}$ whose objects are tuples $(P\to S,p_1,\ldots,p_{2n})$ where:
\begin{itemize}
    \item the morphism $P\to S$ is a twisted conic and $p_i\colon S \to P$ are sections landing in the smooth locus of $P\to S$;
    \item on every geometric point $s$, each irreducible component of $P_s$ contains at least $n$ sections, and no more than $n$ sections can coincide.
\end{itemize}
We omit the definition of the morphisms.
\begin{Prop}\label{prop:hass to git}
    The following hold true:
    \begin{enumerate}
        \item there is a morphism $\pi\colon \Mtilde^{\Hass}_{0,2n}\to\Mbar^{\git}_{0,2n}$ which realizes the first stack as an extended weighted blow-up of the second;
        \item the stack $\Mtilde^{\Hass}_{0,2n}$ contains a proper Deligne--Mumford stack $\Mbar_{0,2n}^{\Hass}$ admitting a projective coarse moduli space.
    \end{enumerate}
\end{Prop}
The Deligne--Mumford stack $\Mbar_{0,2n}^{\Hass}$ appearing in \Cref{prop:hass to git} is the moduli stack of stable weighted $2n$-marked genus zero curves with weights $(\frac{1+\epsilon}{n},\ldots,\frac{1+\epsilon}{n})$ constructed in \cite{Hassett_w}.
We will need the following lemma.
\begin{Lemma}\label{lm:stein of coh aff is gms}
    Let $f\colon \cX \to Y$ be a cohomologically affine morphism. Then $\pi\colon \cX\to X=\spec_Y(f_*\cO_{\cX})$ is a good moduli space.
\end{Lemma}
\begin{proof}
    Consider the canonical factorization $\cX\overset{\pi}{\to}X \overset{g}{\to} Y$. The functor $g_*\colon \operatorname{QCoh}(X) \to \operatorname{QCoh}(Y)$ is an equivalence of categories, because $g$ is affine. We deduce that $\pi_*$ is exact. We conclude by observing that by construction $\pi_*\cO_{\cX}\simeq \cO_X$.
\end{proof}

\begin{proof}[Proof of \Cref{prop:hass to git}]
    As before, we divide the proof in several steps.

    \noindent
    \textbf{Step 1}: the morphism $\pi\colon \Mtilde^{\Hass}_{0,2n}\to\Mbar^{\git}_{0,2n}$ is constructed exactly as in the first step of the proof \Cref{thm:CY extended blow up}, with the only difference that one has to take the image of the sections $p_i\colon S\to P$ instead of the image of the divisor $D\subset P$.     
    
    \noindent
    \textbf{Step 2}: we want to apply \Cref{cor:criterion 2}. First observe that there is a schematically-dense open embedding $i\colon \Mbar^{\git}_{0,2n}\to \Mtilde^{\Hass}_{0,2n}$, defined similarly as in the proof of \Cref{thm:CY extended blow up}. It is easy to check that $\pi\circ i = \id$.
    
    Next, we show that $\Mtilde_{0,2n}^{\Hass}$ admits a good moduli space and that $\pi\colon\Mtilde_{0,2n}^{\Hass}\to\Mtilde^{\git}_{0,2n}$ induces an isomorphism at the level of good moduli spaces.

    Observe that there is a natural morphism $\Mtilde^{\Hass}_{0,2n} \to \Ctilde^{\CY}_{2n}$, which is representable and finite: indeed, it is induced by taking the quotient with respect to the natural action of the symmetric group on the labels of the sections.

    As $\Ctilde^{\CY}_{2n} \to C^{\git}_{2n}$ is cohomologically affine by \Cref{thm:CY extended blow up}, so it is the composition $f\colon \Mtilde^{\Hass}_{0,2n} \to C^{\git}_{2n}$. Define $X={\spec}_{C^{\git}_{2n}}(f_*\cO_{\Mtilde^{\Hass}_{0,2n}})$. Applying \Cref{lm:stein of coh aff is gms}, we deduce that $\Mtilde^{\Hass}_{0,2n} \to X$ is a good moduli space.

    Let $\overline{M}_{0,2n}^{\git}$ be the moduli space of $\Mbar^{\git}_{0,2n}$: then the composition $\Mtilde^{\Hass}_{0,2n} \to \Mbar^{\git}_{0,2n} \to \overline{M}_{0,2n}^{\git}$ induces a morphism $f\colon X \to \overline{M}_{0,2n}^{\git}$ which can be easily checked to be birational and bijective. As the target is a normal variety and both the domain and the target are projective, by Zariski's main theorem we deduce that $f$ is an isomorphism.

    \noindent
    \textbf{Step 3}: we show that there is a representable morphism $\Mtilde^{\Hass}_{0,2n} \to \Mbar^{\git}_{0,2n}\times [\bA^1/\Gm]$. The morphism to $[\bA^1/\Gm]$ is induced by the Cartier divisor of singular twisted conics, so that we have a commutative diagram
    \[
    \begin{tikzcd}
        \Mtilde^{\Hass}_{0,2n} \ar[r] \ar[d] & \Mbar^{\git}_{0,2n}\times [\bA^1/\Gm] \ar[d]\\
        \Ctilde^{\CY}_{2n} \ar[r] & \cC^{\git}_{2n}\times [\bA^1/\Gm].
    \end{tikzcd}
    \]
    As the vertical arrows are representable, and as we have seen in the proof of \Cref{thm:CY extended blow up} that the bottom horizontal arrow is representable, we deduce that also the top horizontal arrow is representable.

    \noindent
    \textbf{Step 4}: we can apply the criterion given in \Cref{cor:criterion 2}. To conclude, observe that $\Mtilde^{\Hass}_{0,2n}$ contains the Hassett moduli stack $\Mbar^{\Hass}_{0,2n}$ of weighted $2n$-marked rational curves, with weights $(\frac{1+\epsilon}{n},\ldots,\frac{1+\epsilon}{n})$ constructed in \cite{Hassett_w}.
\end{proof}
\begin{Remark}
    It follows from the results in this appendix that one can compactify the moduli space of maps $C\to \cC_{2n}^{\git}$ as in \Cref{thm_intro_simplified} by considering the enlargement $\cC_{2n}^{\git}\subseteq \widetilde{\cC}^\CY_{2n}$ and via \Cref{teo_intro_cX_contains_open_proper_dm}. In particular, the locus $\cU$ parametrizing maps $\bP^1\to \cC_{2n}^{\git}$ which send the generic point of $\bP^1$ to a strictly stable point, and such that the composition $\bP^1\to\cC^\git_{2n}\to C^\git_{2n}$ with the good moduli space is finite, can be compactified $\cU\subseteq \overline{\cU}$ in a way such that the boundary parametrizes maps $\cC\to \widetilde{\cC}^\git_{2n}$ from a rational twisted curve, and which are stable in the sense of \Cref{def_stable_qmap}. A similar picture goes through if we consider instead $\overline{\cM}^\git_{0,2n}$.
\end{Remark}

\bibliographystyle{amsalpha}
\bibliography{bibliography}

\end{document}